\documentclass[a4paper,11pt,reqno,noindent]{amsart}
\usepackage[centertags]{amsmath}
\usepackage{amsfonts,amssymb,amsthm,dsfont,cases,amscd,esint,enumerate}
\usepackage[T1]{fontenc}
\usepackage[english]{babel}
\usepackage[utf8]{inputenc}
\usepackage{newlfont}
\usepackage{color}
\usepackage[body={15cm,21.5cm},centering]{geometry} 
\usepackage{fancyhdr}
\theoremstyle{plain}
\newtheorem{theor0}{Theorem}[section]
\newenvironment{theor}
  {\pushQED{\qed}\begin{theor0}}
  {\popQED\end{theor0}}
\newtheorem{lem0}[theor0]{Lemma}

\newtheorem{prop0}[theor0]{Proposition}
\newenvironment{prop}
  {\pushQED{\qed}\begin{prop0}}
  {\popQED\end{prop0}}
\newtheorem{cor0}[theor0]{Corollary}

\theoremstyle{definition}
\newtheorem{rems0}[theor0]{Remarks}

\newtheorem{rem0}[theor0]{Remark}
\newenvironment{rem}
  {\pushQED{\qed}\begin{rem0}}
  {\popQED\end{rem0}}

\theoremstyle{plain}
\newtheorem{as0}[theor0]{Assumption}
\newenvironment{as}
  {\pushQED{\qed}\begin{as0}}
  {\popQED\end{as0}}
\newtheorem*{asn0*}{\assumptionnumber}
  \providecommand{\assumptionnumber}{}
  \makeatletter
  \newenvironment{asn0}[2]
   {\renewcommand{\assumptionnumber}{Condition \!(#1) {\normalfont--- #2}}
    \begin{asn0*}
    \protected@edef\@currentlabel{{\normalfont(#1)}}}
   {\end{asn0*}}
  \makeatother

\mathchardef\emptyset="001F
\numberwithin{equation}{section}

\newcommand{\eps}{\varepsilon}
\newcommand{\dv}{{\operatorname{div}}}

\newcommand{\N}{\mathbb N}
\newcommand{\T}{\mathbb T}
\newcommand{\Sp}{\mathbb S}
\renewcommand{\S}{\mathbb S}
\newcommand{\e}{\varepsilon}

\newcommand{\Pc}{\mathcal{P}}
\newcommand{\Pe}{\operatorname{Pe}}
\newcommand{\act}{\operatorname{act}}
\newcommand{\tr}{\operatorname{tr}}
\newcommand{\ro}{\operatorname{ro}}
\newcommand{\fl}{\operatorname{fl}}

\newcommand{\Cc}{\mathcal{C}}
\newcommand{\calC}{\mathcal{C}}
\newcommand{\Rey}{\operatorname{Re}}
\newcommand{\Wi}{\operatorname{Wi}}

\newcommand{\dd}{\mathrm{d} }
\newcommand{\R}{\mathbb R}
\newcommand{\Z}{\mathbb Z}

\newcommand{\loc}{{\operatorname{loc}}}
\newcommand{\Id}{\operatorname{Id}}

\newcommand{\D}{\operatorname{D}}
\newcommand{\Ld}{\operatorname{L}}
\newcommand{\Div}{{\operatorname{div}}}
\newcommand{\Tr}{{\operatorname{tr}}}
\newcommand{\step}[1]{\noindent \textit{Step} #1.}
\newcommand{\substep}[1]{\noindent \textit{Substep} #1.}

\usepackage[colorlinks,citecolor=black,urlcolor=black]{hyperref}

\title[Hydrodynamic limit of multiscale models for suspensions]{Hydrodynamic limit of multiscale viscoelastic models for rigid particle 
suspensions}
\author{
Mitia Duerinckx, \
Lucas Ertzbischoff,\\
Alexandre Girodroux-Lavigne, \
Richard M. H\"ofer}
\address[Mitia Duerinckx]{Universit\'e Libre de Bruxelles, D\'epartement de Math\'ematique, 1050~Brussels, Belgium}
\email{mitia.duerinckx@ulb.be}
\address[Lucas Ertzbischoff - \textit{corresponding author}]{CEREMADE, CNRS, Université Paris-Dauphine, PSL Research University, 75016 Paris, France}
\email{ertzbischoff@ceremade.dauphine.fr}
\address[Alexandre Girodroux-Lavigne]{University of Regensburg, Faculty of Mathematics, Regensburg, Germany}
\email{alexandre.girodroux-lavigne@mathematik.uni-regensburg.de}
\address[Richard M. H\"ofer]{University of Regensburg, Faculty of Mathematics, Regensburg, Germany}
\email{richard.hoefer@ur.de}

\begin{document}
\selectlanguage{english}

\begin{abstract}
We study the multiscale viscoelastic Doi model for suspensions of Brownian rigid rod-like particles, as well as its generalization by Saintillan and Shelley for self-propelled particles.
{We consider the regime of a small Weissenberg number, which corresponds to a fast rotational diffusion compared to the fluid velocity gradient, and we analyze the resulting hydrodynamic approximation.} More precisely, we show the asymptotic validity of macroscopic nonlinear viscoelastic models, in form of so-called ordered fluid models,
as an expansion in the Weissenberg number.
The result holds for zero Reynolds number in 3D and for arbitrary Reynolds number in~2D. Along the way, we establish several new well-posedness and regularity results for nonlinear fluid models, which may be of independent interest.
\end{abstract}

\maketitle
\setcounter{tocdepth}{1}
\tableofcontents

\section{Introduction}
\begingroup\allowdisplaybreaks
\subsection{General overview}
Suspensions of rigid particles in Stokesian fluids are ubiquitous both in nature and in applications and are known to 
typically
exhibit non-Newtonian behaviors.
These systems can be described on different scales: either by macroscopic non-Newtonian fluid models,
or by so-called multiscale kinetic models, or else on the microscopic scale as suspended particles moving in a fluid flow.
{Let us briefly describe these different levels of physical description:}
\begin{enumerate}[$\bullet$]
\item {\it Macroscopic non-Newtonian fluid models:}\\
An incompressible fluid flow is generally described by the Navier--Stokes equations,
\begin{equation}\label{eq:NS00}
\left\{\begin{array}{l}
\Rey(\partial_t +u\cdot\nabla)u-\Div(\sigma) + \nabla p\,=\,h,\\[1mm]
\Div(u)=0,
\end{array}\right.
\end{equation}
where $u:\R^+\times\R^d\to\R^d$ is the fluid velocity field,
where the source term $h$ accounts for internal forces and where $\sigma$ is the deviatoric fluid stress (a trace-free symmetric matrix field).
For Newtonian fluids, the stress is linear in the strain rate, that is,
\begin{equation}\label{eq:Newt}
\sigma = 2  \D(u),
\end{equation}
where the strain rate $\D(u) = \frac12(\nabla u + (\nabla u)^T)$ is defined as the symmetric gradient of the fluid velocity. In contrast, non-Newtonian fluids are characterized by more complex constitutive laws relating the stress $\sigma$ to the strain rate~$\D(u)$, describing various possible types of nonlinear and memory effects.
In {\it macroscopic} non-Newtonian models,
these laws are assumed to take some explicit form and are
usually fitted phenomenologically to experimental rheological measurements.
In  Reiner--Rivlin and in generalized Newtonian fluid models, the stress $\sigma$ is simply taken to be a local function of~$\D(u)$, meaning that only nonlinear effects are retained.
To further describe memory effects of non-Newtonian fluids, such as viscoelastic properties, more realistic models rather relate $\sigma$ and $\D(u)$ via integral or differential equations, which aim to take into account the dependence of the stress on the fluid deformation history. Such models that are frame invariant are generically called {\it simple fluids}, of which the celebrated Oldroyd--B and {FENE--P} models are particular cases. In the fast relaxation limit, simple fluid models reduce to the hierarchy of so-called ordered fluid models.
{We refer to Section~\ref{sec:standard-2nd-order} for a detailed account of such models.}
\smallskip\item {\it Multiscale kinetic models for suspensions:}\\
So-called {\it multiscale} or {\it micro-macro} models go one step away from the pure phenomenology, towards a microscopically more accurate description of particle suspensions.
More precisely, the macroscopic fluid equation~\eqref{eq:NS00} gets coupled via the stress $\sigma$ to a kinetic equation describing
the evolution of the suspended solid phase.
The latter is modeled by a particle density function $f(t,x,n) \in \R^+$ at time $t$, where~$x$ is the position of particles and~$n$ is their `state': for instance, $n \in \R^d$ may describe the relative position of endpoints of elongated
particles {(hence $n\in\Sp^{d-1}$ in case of rigid suspended particles as will be considered in this work).}
The evolution of the particle density~$f$ is then modeled by a Fokker--Planck equation describing transport with the fluid and diffusive effects. Finally, the coupling to the macroscopic fluid equation is expressed through an explicit constitutive law $\sigma = \sigma(f,\D(u))$, which is typically derived from formal microscopic considerations in dilute regime.
Such models describe how the microscopic state of the particles adapts collectively to local fluid deformations and how the macroscopic fluid flow gets itself effectively impacted.
Popular models include the kinetic FENE and Hookean dumbbell models for dilute suspensions of flexible polymers, the so-called Doi model for suspensions of Brownian rigid rod-like particles,
and the Doi--Saintillan--Shelley for corresponding active (self-propelled) particles.
{The Doi and Doi--Saintillan--Shelley models are displayed in the next section, see~\eqref{eq:eps-kin00}--\eqref{eq:def-sigma2-00}, and we refer to Section~\ref{sec:nondimensionalization} for modelling details.}
\smallskip\item {\it Microscopic models:}\\
At the particle scale, we can formulate a fully detailed hydrodynamic model describing the motion of suspended particles in a background fluid flow. The fluid phase is then described by the Navier--Stokes equations with no-slip conditions at the boundary of the particles and is coupled to the particle dynamics driven by Newton's equations of motion. {This leads to a highly complex dynamics, for which we refer for instance to~\cite{Hillairet-Sabbagh} and references therein; we skip over the details as we will not discuss this type of models further in this work.}
\end{enumerate}
From the modeling perspective, while microscopic models are certainly impractical due to their complexity and the huge number of particles in real-life systems, one can argue that multiscale kinetic models are more satisfactory than macroscopic fluid models as they retain some information from the fluid-particle coupling and therefore better reveal mechanisms leading to non-Newtonian behavior.
Yet, macroscopic models appeal through their simpler description: in particular, they are much more accessible for numerical simulations and proner to comparison with experimental rheological measurements.

Most previous works on particle suspensions have aimed either to study properties of macroscopic non-Newtonian fluid models or multiscale kinetic models, or else to derive multiscale kinetic models rigorously from microscopic particle dynamics, which has indeed attracted considerable interest in recent years; see Section~\ref{sec:previous.results} below for references.
In the present contribution, we fill the gap in the micro-macro understanding of non-Newtonian effects of particle suspensions by further studying the derivation of macroscopic fluid models from multiscale kinetic models in suitable regimes.

In some exceptional cases, the kinetic equation describing the state of the particles in multiscale kinetic models can be integrated out and leads to a closed macroscopic equation for the stress $\sigma$ in terms of the strain rate $\D(u)$. This is for instance the case for the kinetic Hookean dumbbell model, which is well-known to be formally equivalent to the Oldroyd--B macroscopic fluid model (see e.g.\@ the recent rigorous analysis in~\cite{DebiecSuli23}).
In general, however, no exact macroscopic closure is available and we can only hope for perturbative closures to be valid in suitable asymptotic regimes.

In fact, in the regime of weak non-Newtonian effects, a fairly large class of non-Newtonian macroscopic fluid models
is believed to be well approximated by {a special family of models, called {\it ordered fluid} models. More precisely, the latter} are expected to be good approximations for all viscoelastic
fluids {in the fading memory regime}, that is, in the regime when the elastic time-dependent effects due to suspended particles in the fluid have an inherent relaxation timescale that is much shorter than the overall timescale of the fluid flow.
This ratio of timescales is the so-called Weissenberg number $\Wi$. In other words, ordered fluid models arise formally as expansions of {any viscoelastic
fluid model} at small~$\Wi$ (sometimes referred to as the retarded motion expansion). A first-order fluid is a Newtonian fluid, a second-order fluid is a non-Newtonian fluid where effects of order $O(\Wi^2)$ are neglected, a third-order fluid amounts to neglecting effects of order $O(\Wi^3)$, etc.

\medskip
In the present work, we focus on the Doi model, which is a multiscale kinetic model for suspensions of Brownian rigid rod-like particles, and we further consider its generalization by Saintillan and Shelley for active (self-propelled) particles.
We rigorously analyze {the hydrodynamic limit of these models in the small-$\Wi$ regime,} which was extensively studied on a formal level in the physics literature in the 1970s~\cite{HinchLeal72, Brenner74},
and we confirm the asymptotic validity of ordered fluid models in this setting.
As a prerequisite, the justification of the asymptotic expansion requires a careful study both of the Doi--Saintillan--Shelley model and of ordered fluid models: in particular, we establish new well-posedness and regularity results for these nonlinear viscoelastic fluid models, which may be of independent interest. Moreover, as the Doi--Saintillan--Shelley model
can itself be derived from a microscopic hydrodynamic model (at least formally), see Section~\ref{sec:previous.results},
our derivation of macroscopic ordered fluid models comes together with explicit expressions for rheological parameters in terms of microscopic characteristics of the underlying particle suspension, {which are checked to fit qualitatively with experiments. In short, while ordered fluid models have often been dismissed due to ill-posedness and instability issues, we restore their credibility by showing that they are perfectly well-posedness in some suitable perturbative sense, that they can be rigorously derived from multiscale kinetic models for suspensions, and that they predict relevant non-Newtonian effects observed in rheological experiments --- while more reputable and simpler formal closures like the Oldroyd--B model are found to be inaccurate.}

\subsection{Informal statement of main results}
We start from the following dimensionless Doi--Saintillan--Shelley model describing suspensions of (active or passive) Brownian rigid rod-like particles in a Stokesian fluid
\begin{equation}\label{eq:eps-kin00}
\left\{\begin{array}{l}
\Rey(\partial_t+u_\e\cdot\nabla)u_\e-\Delta u_\e+\nabla p_\e=h+ \frac{1}{\e}\Div(\sigma_1[f_\e])+\Div(\sigma_2[f_\e,\nabla u_\e]),\\[2mm]
\partial_t f_\e+\Div_x\big((u_\e+U_0n)f_\e\big)+\Div_n\big(\pi_n^\bot(\nabla u_\e)nf_\e\big)=\tfrac1 \Pe \Delta_xf_\e+\tfrac1\e\Delta_nf_\e,\\[2mm]
\Div(u_\e)=0,
\end{array}\right.
\end{equation}
In this system, the unknowns are the fluid velocity field $u\colon\R^+\times\T^d\to\R^d$ and the particle position and orientation density function $f\colon\R^+\times\T^d\times\Sp^{d-1}\to\R^+$, while the elastic and viscous stresses are defined as 
\begin{eqnarray}
\sigma_1[f]&:=&\lambda\theta\int_{\Sp^{d-1}}\big(n\otimes n-\tfrac1d\Id\big)f(\cdot,n)\,dn,\\[2mm]
\sigma_2[f,\nabla u]&:=&\lambda\int_{\Sp^{d-1}}(n\otimes n)(\nabla u)(n\otimes n)\,f\,dn.\label{eq:def-sigma2-00}
\end{eqnarray}
The parameter $\e:=\Wi>0$ stands for the Weissenberg number, {$\Rey \ge0$ is the} Reynolds number, $\Pe>0$ is the so-called P\'eclet number, and the source term $h$ is taken to be smooth and accounts for internal forces.
The state variable $n$ in the kinetic Fokker--Planck equation for $f$ describes the orientation of rigid particles on the unit sphere, hence $n\in \Sp^{d-1}$: we write~$\Delta_n$ and~$\Div_n$ for the Laplace--Beltrami operator and the divergence on the sphere~$\Sp^{d-1}$,
and we also use the short-hand notation $$\pi_n^\bot:=\Id-n\otimes n$$ for the orthogonal projection onto $n^\bot$.
{We refer to Section~\ref{sec:nondimensionalization} for references and modelling details on this Doi--Saintillan--Shelley system~\eqref{eq:eps-kin00}, together with the non-dimensionalization procedure and the description of the different parameters $U_0,\lambda,\theta\in\R$.}
Note that the so-called Doi model for passive suspensions is recovered for the special choice $U_0=0$ and~$\theta=6$.
We consider this system for simplicity in a finite box $\T^d=[0,1)^d$ with periodic boundary conditions, and we consider space dimension~$d=2$ or~$3$.
The system is complemented by initial conditions
\begin{equation}\label{eq:eps-kin-bis}
\left\{\begin{array}{l}
f_\e|_{t=0}=f^\circ_\e,\\[1mm]
u_\e|_{t=0}=u^\circ~~\text{if $\Rey\ne0$},
\end{array}\right.
\end{equation}
and for uniqueness we require
\begin{align*}
\int_{\T^d} u_\e=0\quad\text{if $\Rey=0$}.
\end{align*}
{In this setting, we shall study the hydrodynamic limit as $\e\downarrow0$ in the regime with $\lambda\le1$ and with
$\lambda\theta\ll1$ small enough, and our main result is the following asymptotic validity of second-order fluid models.
We refer to Theorem~\ref{th:main} in Section~\ref{sec:main.result} for a more detailed statement}, including the precise well-preparedness requirement for kinetic initial data, as well as the explicit expression for the effective coefficients $\mu_0,\eta_1,\gamma_1,\gamma_2$.
{New results on the well-posedness of the kinetic model~\eqref{eq:eps-kin00} and of the second-order fluid model~\eqref{eq:nonstand-REtensor-intro} are also postponed to Section~\ref{sec:main.results}.}

\begin{theor}[Informal statement of the main result]\label{thm:informal}
Consider either the Stokes case $\Rey=0$ with $d\le3$, or the Navier--Stokes case $\Rey\ne0$ with $d=2$.
Given an initial particle density $f_\e^\circ \in C^\infty\cap\Pc(\T^d\times\Sp^{d-1})$ {that is well-prepared in a sense that will be clarified later (see Assumption~\ref{as:well-prepared} in Section~\ref{sec:formal.expansion}),} and given also an initial fluid velocity $u^\circ \in C^\infty(\T^d)^d$ with $\Div(u^\circ)=0$ in the Navier--Stokes case,
consider a weak global solution $(u_\e,f_\e)$ of the Cauchy problem for the Doi--Saintillan--Shelley model~\eqref{eq:eps-kin00}. For all $T > 0$,
{there is a constant $C<\infty$ such that the following holds: provided that
\begin{equation}\label{eq:smallness-cond}
\e\le\tfrac1C\qquad\text{and}\qquad \lambda\theta (1+\Pe)\le\tfrac1C,
\end{equation}
the fluid velocity $u_\e$ satisfies
\begin{eqnarray*}
\begin{array}{rlll}
\|\nabla(u_\e-\bar u_\e)\|_{\Ld^2 (0,T; \Ld^2(\T^d)^{d\times d})}&\le&C\e^2,&\quad\text{if~~$\Rey=0$},\\[1mm]
\|u_\e-\bar u_\e\|_{\Ld^\infty (0,T; \Ld^2(\T^d)^d)}+\|\nabla(u_\e-\bar u_\e)\|_{\Ld^2 (0,T; \Ld^2(\T^d)^{d\times d})}&\le&C\e^2,&\quad\text{if~~$\Rey\ne0$},
\end{array}
\end{eqnarray*}
and the particle spatial density $\rho_\e:=\fint_{\Sp^{d-1}}f_\e(\cdot,n)\,\dd n$ satisfies
\[\| \rho_\e - \bar \rho_\e\|_{\Ld^\infty (0,T; \Ld^2(\T^d))} + \|\nabla( \rho_\e - \bar \rho_\e)\|_{\Ld^2 (0,T;\Ld^2(\T^d)^d)} \,\le\,C\e^2,\]
where} $(\bar u_\e, \bar \rho_\e)$ solves the following (non-standard) second-order fluid equation
\begin{align}\label{eq:second-order-fluid-intro-eqn}
\left\{\begin{array}{l}
\Rey(\partial_t+\bar u_\e\cdot\nabla)\bar u_\e-\Div(\bar\sigma_\e)+\nabla\bar p_\e\,=\,h+O(\e^2),\\[1mm]
(\partial_t+\bar u_\e\cdot\nabla)\bar\rho_\e-(\tfrac 1 \Pe+\e\mu_0)\Delta\bar\rho_\e=O(\e^2),\\[1mm]
\bar\sigma_\e=(1+\eta_1\bar\rho_\e)A_1(\bar u_\e)+\e\gamma_1A_2'(\bar u_\e,\bar\rho_\e)+\e\gamma_2\bar\rho_\e A_1(\bar u_\e)^2,\\[1mm]
\Div(\bar u_\e) = 0,
\end{array}\right.
\end{align}
for some explicit coefficients $\mu_0,\eta_1,\gamma_1,\gamma_2 \in \R$,
where 
\begin{align} \label{strain.rate}
    A_1(u):=2\D(u)=(\nabla u)^T + \nabla u
\end{align} 
is the strain rate
and where $A_2'$ is a (non-standard) inhomogeneous, diffusive version of the second Rivlin-Ericksen tensor, defined by
\begin{equation}\label{eq:nonstand-REtensor-intro}
A_2'(\rho,u)\,:=\,\big(\partial_t-{\tfrac 1 \Pe}\Delta+ u\cdot\nabla\big)\big(\rho A_1 (u)\big)+\rho\big((\nabla u)^TA_1 (u)+A_1 (u)(\nabla u)\big).\qedhere
\end{equation}
\end{theor}

\begin{rem} \label{rem:main.informal}
A few comments are in order:
\begin{enumerate}[(a)]
\item {\it Limitation regarding Reynolds number:}\\Due to classical regularity issues for the Navier--Stokes equations, we limit ourselves to the Stokes case $\Rey = 0$ for the 3D model, and we only study the Navier--Stokes case $\Rey \ne0$ in 2D. In contrast, the 3D Navier--Stokes case is much more complicated and is not discussed in this work as the well-posedness of the system is then still open.
\smallskip\item {\it Limitation regarding P\'eclet number:}\\We take into account a non-vanishing spatial diffusion $O(\frac1\Pe)$ in the Doi--Saintillan--Shelley model~\eqref{eq:eps-kin00}, which differs {from the setting usually considered in applications where $\Pe=\infty$: a nontrivial spatial diffusion $\frac1\Pe>0$} is actually needed in this work for technical well-posedness reasons.
{In fact, in view of~\eqref{eq:smallness-cond}, we could consider a small spatial diffusion provided that the particle volume fraction $\lambda$ is small enough.}
Due to the spatial diffusion, the structure of ordered fluid equations needs to be slightly adapted, leading to the non-standard definition of Rivlin--Ericksen tensors in~\eqref{eq:nonstand-REtensor-intro}. {We refer to Sections~\ref{sec:nonstandard-ordered-Pe} and~\ref{sec:ordered.nonstandard} for full details on this modification of ordered fluid models at $\Pe<\infty$. Note that it is reminiscent of the version of the Oldroyd--B model with stress diffusion that is often considered both for analytical and numerical studies~\cite{renardyThomases21,DebiecSuli23}}.
\smallskip\item {\it Well-preparedness:}\\The well-preparedness assumption for the initial kinetic data will be clarified later on in Assumption~\ref{as:well-prepared}. Informally, it amounts to assuming that initial data are locally at equilibrium with respect to the dynamics of orientations and are perturbatively compatible with the formal $\e$-expansion. It allows to avoid initial boundary layers.
\smallskip\item {\it Well-posedness of the models:}\\
Along with the above asymptotic expansion result, we prove well-posedness both for the Doi--Saintillan--Shelley model and for the second-order fluid equation:\vspace{0.1cm}
\begin{enumerate}[---]
\item The existence of global weak solutions for the Doi--Saintillan--Shelley model~\eqref{eq:eps-kin00} is postponed to Section~\ref{sec:well-posed-DSS}. We further obtain a new weak-strong uniqueness principle, which implies some stability that is at the very heart of the above result.
\vspace{0.1cm}\item The second-order fluid equation~\eqref{eq:second-order-fluid-intro-eqn} is well-known to be ill-posed whenever~\mbox{$\gamma_1<0$} (which is indeed the case for relevant effective parameters): {this is particularly clear at $\Pe<\infty$ as the equation then behaves like a backward heat equation (see also e.g.~\cite{Galdi-93} at $\Pe=\infty$).}
Yet, we can define several well-posed notions of approximate solutions that only satisfy~\eqref{eq:second-order-fluid-intro-eqn} up to higher-order~$O(\e^2)$ errors.
This discussion is postponed to Section~\ref{sec:formal.expansion} and Appendix~\ref{app.Bous}, where {we present two approaches to fix this issue: First,} we introduce the notion of approximate \textit{hierarchical solutions}, {which naturally appear as small-$\e$ expansions; see Proposition~\ref{prop:hier-inhom}. Second,} by means of a Boussinesq-type perturbative rearrangement, we additionally {provide a reformulation} of second-order fluid equations {in terms of a closed} well-posed system, {which may be more desirable in particular for stability issues;} see Proposition~\ref{prop:Bous}.
{This allows us to overcome the usual objections about the ill-posedness and instability of ordered fluid equations, thus restoring their legitimity.}
\end{enumerate}
\smallskip\item {\it Non-Newtonian features in the hydrodynamic limit:}\\
The explicit expression for the effective second-order fluid coefficients $\mu_0,\eta_1,\gamma_1,\gamma_2$ is postponed to Section~\ref{sec:main.result}, and we further refer to Section~\ref{secrem:non-newtonian actif} for the description of rheological properties of the obtained macroscopic model. The expressions for the coefficients agree with those computed in~\cite{HinchLeal72, Brenner74} in the case of passive suspensions, and they qualitatively match experimental data and formal predictions on active suspensions.
{In particular, for passive particles, we recover non-zero normal-stress differences, with first and second normal-stress coefficients that are positive and negative, respectively, exactly as observed in experiments. Note that for passive suspensions this actually rules out the validity of the Oldroyd--B model, which would rather predict a zero second normal-stress coefficient (see Remark~\ref{rem:Oldroyd-lala}).}
\smallskip\item {\it Extensions:}\\A similar result could be obtained with the same approach when starting from the co-rotational kinetic FENE model for elastic polymers (see e.g.~\cite{Lions-Masmoudi-07} for a review of this system). For conciseness, we do not repeat our analysis in that setting and leave the adaptation to the reader.
\qedhere
\end{enumerate}
\end{rem}\medskip

\subsection{Previous results}\label{sec:previous.results}
We briefly review previous results, both rigorous or not, related to the multiscale description of particle suspensions and related systems:

\begin{enumerate}[$\bullet$]
\item {\it Macroscopic non-Newtonian rheology.}\\
{Since the early works that have led to the formal derivation of the Doi or Doi--Saintillan--Shelley models in the physics literature, see e.g.~\cite{Shima40, KuhnKuhn45, RisemannKirkwood50, Saito51},
the underlying goal has been to actually derive macroscopic non-Newtonian features, typically} via the explicit calculations of stresses in specific flows such as simple shear.
This has led for instance to asymptotic formulas for the shear viscosity in simple shear flow for very large or very small Weissenberg number~$\Wi$. Normal-stress differences (see~\eqref{definition normal stress} below) have also been computed at small $\Wi$ by Giesekus~\cite{Giesekus62}, showing that the elastic stress $\sigma_1$ does not contribute to the second normal-stress difference but that the viscous stress does. This was in contradiction with Weissenberg's original conjecture that all real-world fluids must have vanishing second normal-stress difference, a conjecture that was later falsified also through experiments.
Hinch and Leal~\cite{HinchLeal71,HinchLeal72} systematically computed expansions for the stress in simple shear flow both at small and at large $\Wi$. They also noticed that their findings at small $\Wi$ agree with second-order fluid models, but they did not investigate whether this holds in more general fluid flows.
{Similar results have been obtained by Brenner \cite{Brenner74} who extended the expansion to third-order fluids.}

\smallskip
\item{\it Validity of formal closures.}\\
Although no exact macroscopic closure is available for the Doi and Doi--Saintillan--Shelley models, formal approximate {moment closures have been proposed in~\cite[Chapter 8.7]{DoiEdwards88}. Such moment closures are  studied  more systematically both analytically and numerically in \cite{HelzelTzavaras16, HelzelTzavaras17, DahmHelzel22} with a special focus on sedimenting particles.}
A particular instance is the Oldroyd--B model, which is well-known to be remarkably an exact macroscopic closure for the kinetic Hookean dumbbell model (see e.g.\@ the recent rigorous analysis in~\cite{DebiecSuli23} in the stress-diffusive case).
For the Doi model, as an exact closure is not available, we may wonder whether macroscopic closures are at least asymptotically valid in some scaling regimes: in fact, we will see that the validity of the Oldroyd--B closure already fails at order $O(\Wi)$ in the small-$\Wi$ regime (see Remark~\ref{rem:Oldroyd-lala}).
The situation is very similar for kinetic FENE models for elastic polymers:
no exact macroscopic closure holds, but the so-called FENE--P model is still a very popular formal approximate closure (see e.g.~\cite{BirdDotsonJohnson80}, where the name FENE--P is attributed in reference to an earlier paper of Peterlin~\cite{Peterlin66}).

\smallskip
\item {\it Small-$\Wi$ regime and hydrodynamic limits.}\\
{At small $\Wi$, the Doi--Saintillan--Shelley model~\eqref{eq:eps-kin00} undergoes a strong diffusion in orientation (cf.\@ factor $\frac1\e=\frac1{\Wi}$ in front of $\Delta_n$): to leading order, the orientations of the particles relax instantaneously to the steady state, which corresponds to isotropic orientations. While this leading order amounts to a trivial Newtonian behavior, next-order $O(\Wi)$ corrections encode nontrivial non-Newtonian effects where the stress starts to depend on the local fluid deformation and its history.
The formal perturbative expansion in powers of $\Wi$ is comparable to the Hilbert expansion method in the Boltzmann theory~\cite{Hilbert-17,Caflisch-80,Golse-05,Saint-Raymond-09}, where one looks for a solution of the  Boltzmann equation as a formal power series in terms of the Knudsen number $\mathrm{Kn} \ll 1$ and where the leading-order approximation simply leads to compressible Euler equations. From the closely related Chapman--Enskog asymptotics expansion, one can (formally) obtain compressible Navier--Stokes equations as a $O(\mathrm{Kn})$ correction to the compressible Euler system, up to $O(\mathrm{Kn}^2)$ errors (see e.g.~\cite[Section~5.2]{Golse-05} or~\cite[Section~2.2.2]{Saint-Raymond-09}).
In a similar way, for the Doi model, second-order fluid equations are obtained in this work as a $O(\Wi)$ correction to the Stokes equations, up to $O(\Wi^2)$ errors.}

Our results can be compared with corresponding results for the Doi--Onsager model for liquid crystals, which indeed shares some similarities with the Doi model that we consider in this work. In~\cite{E-Zhang-06,Wang-Zhang-Zhang-15}, the macroscopic Ericksen--Leslie system is derived from the Doi--Onsager model by means of a Hilbert expansion.
Note however that in that case the leading term in the expansion already yields a non-trivial system, so that higher-order corrections are not investigated in~\cite{E-Zhang-06,Wang-Zhang-Zhang-15}.
We also mention recent related work on hydrodynamic limits for alignment models~\cite{Degond-Motsch-08,DFMA-17,DFMAT-19,DegondFrouvelleLiu21} and for flocking models~\cite{karper2015hydrodynamic,kang2015asymptotic,figalli2018rigorous},
{as well as a preliminary result for kinetic FENE models for elastic polymers~\cite{LemaouPicassoDegond02}.}

Note that hydrodynamic limits have been investigated also for various other kinetic models for particle suspensions in different settings:
for instance, in the context of sedimentation for small inertial particles, let us mention the inertialess limits studied in~\cite{Jab00, Hofer18InertialessLimit, HanKwan-Michel,Ertzbischoff-Boussinesq}, as well as the high-diffusion limit in velocity investigated in~\cite{GoudonJabinVasseur04a,GoudonJabinVasseur04b,Mellet-Vasseur,Su-Yao}.

In a different direction, we also mention recent work~\cite{HohShellStab2010, AlbrittonOhm22, CotiZelatiDietertGerardVaret22, CotiZelatiDietertGerardVaret24}, where the stability and mixing properties of the Doi--Saintillan--Shelley model have been investigated {(neglecting however the viscous stress~$\sigma_2$)}.
\end{enumerate}

\subsection{Structure of the article}
{The article is split into five main sections, in addition to three appendices:
\begin{enumerate}[$\bullet$]
\item In Section~\ref{sec:main.results}, we state our main results. We start with the well-posedness of the Doi--Saintillan--Shelley system, cf.~Proposition~\ref{prop:well-posedness}, and the perturbative well-posedness of second-order fluid equations, cf.~Proposition~\ref{prop:hier-inhom}, before turning to the detailed statement of the hydrodynamic approximation result that we anticipated in Theorem~\ref{thm:informal}, cf.~Theorem~\ref{th:main}. In Section~\ref{sec:Hilbert-formal}, we also include a brief account of the formal $\e$-expansion that motivates our different results.
\smallskip\item Section~\ref{sec:expansion} is devoted to the proof of the hydrodynamic approximation result.
\smallskip\item In Section~\ref{sec:Doi-physics}, we comment on the different models that we consider in this work and on the physical content of our results. More precisely, in Section~\ref{sec:nondimensionalization}, we motivate the Doi--Saintillan--Shelley model~\eqref{eq:eps-kin00} that is our starting point in this work and we describe the signification of the different parameters. Next, in Section~\ref{sec:standard-2nd-order}, we recall the usual definition of ordered fluid models and their non-Newtonian features, before turning to their non-standard modifications at finite P\'eclet number and for inhomogeneous suspensions, leading to the specific equations~\eqref{eq:second-order-fluid-intro-eqn} that we derive in this work. We also compare the non-Newtonian features of the derived models with experimental data.
\smallskip\item In Section~\ref{sec:Conclusion}, we include a few concluding remarks on our main results.
\smallskip\item In Appendix~\ref{app:well-posedness}, we give the proof of Proposition~\ref{prop:well-posedness} on the well-posedness of the Doi--Saintillan--Shelley system~\eqref{eq:eps-kin00}.
\smallskip\item In Appendix~\ref{app.Bous}, we  give the proof of Proposition~\ref{prop:hier-inhom} on the perturbative well-posedness of second-order fluid equations, and we further develop alternative notions of perturbative solutions.
\smallskip\item In Appendix~\ref{app:3rd-order}, we derive the next-order $\e$-expansion of the Doi--Saintillan--Shelley theory, which leads to some physically important corrections to second-order fluid models.
\end{enumerate}}

\subsection*{Notations}
We summarize the main notations that we use in this work:
\begin{enumerate}[$\bullet$]
\item We denote by $C\ge1$ any constant that only depends on the dimension $d$ and possibly on other controlled quantities to be specified. We use the notation $\lesssim$ for $\le C\times$ up to such a multiplicative constant $C$. We write $\ll$ (resp.~$\gg$) for $\le C\times$ (resp.~$\ge C\times$) up to a sufficiently large multiplicative constant $C$. When needed, we add subscripts to indicate dependence on other parameters.
\smallskip\item For a vector field $u$ and a matrix field $S$, we set $(\nabla u)_{ij}:=\nabla_j u_i$, $S^T_{ij}:=S_{ji}$, $\D( u):=\frac12(\nabla u+(\nabla u)^T)$, and $\Div(S)_i:=\nabla_j S_{ij}$ (we systematically use Einstein's summation convention on repeated indices).
\smallskip\item We denote by $\dd n$ the (not normalized) Lebesgue measure on the $(d-1)$-dimensional unit sphere $\Sp^{d-1}$,
and we denote its area by $\omega_d := |\S_{d-1}|$.
Differential operators with a subscript $n$ (such as $\mathrm{div}_n$ and $\Delta_n$) refer to differential operators on $\Sp^{d-1}$, endowed with the natural Riemannian metric.
\smallskip\item For $n \in \Sp^{d-1}$, we denote by $\pi_n^\bot:=\Id-n\otimes n$ the orthogonal projection on $n^\bot$.
\smallskip\item We let {$\langle g\rangle(x):=\fint_{\Sp^{d-1}}g(x,n)\,\dd n:=\frac1{\omega_d}\int_{\Sp^{d-1}}g(x,n)\,\dd n$} be the angular averaging of a function $g$ on~\mbox{$\T^d\times\Sp^{d-1}$}. We also use the short-hand notation $P_1^\bot g:=g-\langle g\rangle$.
\smallskip\item We denote by $H^k(\T^d)$ (resp.\@ $H^k(\T^d \times \S^{d-1})$) the standard $\Ld^2$ Sobolev spaces for functions depending on $x \in \T^d$ (resp.\@ $(x,n) \in \T^d \times \S^{d-1}$), and we use the notation $\|\cdot\|_{H^k_x}$ (resp.~$\|\cdot\|_{H^k_{x,n}}$) for the corresponding norms.
For time-dependent functions, given a Banach space $X$ and $t>0$, the norm of $\Ld^p(0,t;X)$ is denoted by $\|\cdot\|_{\Ld_t^p X}$.
\smallskip\item The space of probability measures on $\T^d$ (resp.\@ on $\T^d \times \S^{d-1}$) is denoted by $\Pc(\T^d)$ (resp.\@ by~$\Pc(\T^d \times \S^{d-1})$).
\end{enumerate}

%%%%%%%%%%%%%%%%%%%%%%%%%
%%%%%%%%%%%%%%%%%%%%%%%%%

\section{Main results} \label{sec:main.results}
{In this section, we describe our main results in more details. We start with our well-posedness results for the Doi--Saintillan--Shelley model and for second-order fluid equations, before giving a more complete statement of our hydrodynamic limit result anticipated in Theorem~\ref{thm:informal}. For notational simplicity, we shall consider $\Rey\in\{0,1\}$.}

\subsection{Well-posedness of the Doi--Saintillan--Shelley model}\label{sec:well-posed-DSS}
{Although the Doi--Saintillan--Shelley system~\eqref{eq:eps-kin00} is an instance of the more general class of Fokker--Planck--Navier--Stokes systems, we emphasize two main peculiarities:}
\begin{enumerate}[---]
\item we include the contribution of the viscous stress $\sigma_2$, which arises from the rigidity of underlying suspended particles on the microscale and effectively modifies the solvent viscosity;
\smallskip\item we also include the effect of particle swimming via~$U_0$,
thus creating local changes in the spatial density $\rho_\e=\fint_{\Sp^{d-1}}f_\e(\cdot,n)\,\dd n$, which no longer remains constant in general.
\end{enumerate}
In contrast, most previous works have focused on corresponding models without viscous stress $\sigma_2$ and without particle swimming $U_0=0$.
In that simplified setting, for the 3D Stokes case, the existence of global entropy solutions was proven in~\cite{otto2008continuity}, and the global well-posedness of smooth solutions was proven in~\cite{constantin2010global} (without translational diffusion, $\Pe=\infty$).
In the Navier--Stokes case, corresponding well-posedness results were obtained for instance in~\cite{constantin2008global}.
The model including the viscous stress $\sigma_2[f_\e,\nabla u_\e]$ but without particle swimming $U_0=0$ was first studied in~\cite{Lions-Masmoudi-07}, where the existence of global weak solutions was proven for the Navier--Stokes case in 2D and 3D (without translational diffusion, $\Pe=\infty$). We also refer to~\cite{La19} for the global well-posedness of smooth solutions in the~2D Navier--Stokes case.
Particle swimming $U_0\ne0$ was first considered in~\cite{chen2013global}, where the authors studied the corresponding model without viscous stress $\sigma_2$ and proved the global existence of weak entropy solutions both for the Stokes and Navier--Stokes cases in 2D and 3D, as well as the existence of energy solutions for the Stokes case and their uniqueness in 2D.

{Unsurprisingly, building on similar ideas as in previous work on related systems,} we can actually establish the existence of global energy solutions for the full model~\eqref{eq:eps-kin00} in the 2D and 3D Stokes cases, as well as in the 2D Navier--Stokes case, and we further obtain a weak-strong uniqueness principle.
To our knowledge, this is the first result where both the viscous stress and the swimming forces are included at the same time.
The proof is postponed to Appendix~\ref{app:well-posedness}.

\begin{prop}\label{prop:well-posedness}
Consider either the Stokes case $\Rey=0$ with $d\le3$, or the Navier--Stokes case $\Rey=1$ with $d=2$.
Given $\e>0$, given $h\in\Ld^\infty_\loc(\R^+;\Ld^2(\T^d)^d)$, given an initial condition $f^\circ_\e\in\Ld^2\cap\Pc(\T^d \times\Sp^{d-1})$, and given also $u^\circ\in \Ld^2(\T^d)^d$ with $\Div(u^\circ)=0$ in the Navier--Stokes case, the Cauchy problem~\eqref{eq:eps-kin00}--\eqref{eq:eps-kin-bis} admits a global weak solution $(u_\e,f_\e)$ with:
\begin{enumerate}[(i)]
\item{in the Stokes case $\Rey=0$, $d\le3$,}
\begin{eqnarray*}
u_\e&\in&\Ld^{\infty}_\loc(\R^+;H^1(\T^d)^d),\\[2mm]
f_\e&\in&\Ld^\infty_\loc(\R^+;\Ld^2\cap\Pc(\T^d\times\Sp^{d-1}))\cap\Ld^2_\loc(\R^+;H^1(\T^d\times\Sp^{d-1}));
\end{eqnarray*}
\item{in the Navier--Stokes case $\Rey=1$, $d=2$,}
\begin{eqnarray*}
u_\e&\in&\Ld^{\infty}_\loc(\R^+;\Ld^2(\T^2)^2)\cap \Ld^{2}_\loc(\R^+;H^1(\T^2)^2),\\[2mm]
f_\e&\in&\Ld^\infty_\loc(\R^+;\Ld^2\cap\Pc(\T^2\times\Sp^{1}))\cap\Ld^2_\loc(\R^+;H^1(\T^2\times\Sp^{1})).
\end{eqnarray*}
\end{enumerate}
In both cases, a weak-strong uniqueness principle further holds: if $(u_\e,f_\e)$ and $(u'_\e,f'_\e)$ are two such global weak solutions with identical initial conditions, {and if $(u'_\e,f'_\e)$ has the following additional regularity,
\begin{eqnarray*}
 u'_\e &\in& \Ld^2_{\mathrm{loc}}(\R^+; W^{1, \infty}(\T^d)^d), \\
 f'_\e  &\in& \Ld^{\infty}_{\mathrm{loc}}(\R^+; \Ld^{\infty}(\T^d \times \mathbb{S}^{d-1})),
\end{eqnarray*}
then we have $(u_\e,f_\e)=(u'_\e,f'_\e)$.}
\end{prop}

\subsection{Formal Hilbert expansion}\label{sec:Hilbert-formal}
{Before continuing with the statement of our main results, let us perform the formal Hilbert expansion of the Doi--Saintillan--Shelley model~\eqref{eq:eps-kin00} as $\e\downarrow0$. This will shed light both on the required well-preparedness assumptions and on our perturbative well-posedness result for the ordered fluid equation~\eqref{eq:second-order-fluid-intro-eqn} that we aim to derive.
Let $(u_\e,f_\e)$ be a global weak solution of~\eqref{eq:eps-kin00}. Given the fast relaxation of particle orientations due to the strong rotational diffusion, we} naturally split the particle density as
\[f_\e(x,n)\,=\,\rho_\e(x)+g_\e(x,n),\qquad\rho_\e\,:=\,\langle f_\e\rangle\,:=\,\fint_{\Sp^{d-1}}f_\e(\cdot,n)\,\dd n\,\in\,\tfrac1{\omega_d}\Pc(\T^d),\]
where $\rho_\e$ stands for the spatial density and where~$\langle g_\e\rangle=0$.
{In order to avoid initial boundary layers,} we then make the following well-preparedness assumption: we assume 
\begin{equation}\label{eq:loc-equi-no-wellprepared}
\rho_\e|_{t=0}=\rho^\circ\in\tfrac1{\omega_d} \Pc(\T^d),\qquad g_\e|_{t=0}=O(\e),
\end{equation}
meaning that we start from an initial density that is to leading order at equilibrium with respect to the strong rotational diffusion in~\eqref{eq:eps-kin00}.
In this setting, we shall analyze the asymptotic behavior of the solution $(u_\e,f_\e)$ and derive a hydrodynamic approximation in the spirit of Hilbert's expansion method in the Boltzmann theory~\cite{Hilbert-17,Caflisch-80,Golse-05,Saint-Raymond-09}.
We start from the ansatz 
\begin{align}
\begin{split}\label{eq:formal-exp}
u_\e&=u_0+\e u_1+\e^2u_2+\ldots,\\
\rho_\e&=\rho_0+\e \rho_1+\e^2\rho_2+\ldots,\\
g_\e&=g_0 + \e g_1+\e^2g_2+\e^3g_3+\ldots,
\end{split}
\end{align}
with $\rho_0\in\frac1{\omega_d}\Pc(\T^d)$, with $\int_{\T^d}\rho_k=0$ for $k\ge1$, and with $\langle g_k\rangle=0$ for all $k\ge0$.
Inserting it into the system~\eqref{eq:eps-kin00}, and identifying powers of $\e$, 
we are led formally to the following hierarchy of coupled equations:
\begin{enumerate}[$\bullet$]
\item \emph{Order $O(\e^{-1})$:} The equation for the particle density yields $\Delta_n g_0 = 0$, and therefore, as by definition $\langle g_0 \rangle =0$, we must have
\[g_0=0.\]
On the other hand, the fluid equation yields $\Div_x ( \sigma_1 [\rho_0 + g_0]) = 0$, which is then automatically satisfied as $g_0=0$ and as $\rho_0$ does not depend on $n$.
\smallskip\item \emph{Order $O(\e^{0})$:} The triplet $(u_0,\rho_0,g_1)$ satisfies
\begin{equation}
\left\{\begin{array}{l}\label{eq :(u_0,f_1)}
\Rey(\partial_t+u_0\cdot\nabla)u_0-\Delta u_0+\nabla p_0=h+\Div(\sigma_1[g_1])+\Div(\sigma_2[\rho_0,\nabla u_0]),\\[1mm]
\Delta_n g_1 = U_0n\cdot\nabla_x\rho_0+\Div_n\big(\pi_n^{\perp} (\nabla u_0) n\rho_0\big),\\[1mm]
(\partial_t-\tfrac1 \Pe \Delta+u_0\cdot\nabla)\rho_0=0,\\[1mm]
\Div(u_0)=0,~~\langle g_1\rangle=0,\\[1mm]
u_0|_{t=0}=u^\circ~~\text{if $\Rey\ne0$}, \ \ \int_{\T^d} u_0 = 0~~\text{if $\Rey=0$},\\[1mm]
\rho_0|_{t=0}=\rho^\circ.
\end{array}\right.
\end{equation}
\item \emph{Order $O(\e^{1})$:} The triplet $(u_1,\rho_1,g_2)$ satisfies
\begin{equation}
\left\{\begin{array}{l}\label{eq :(u_1,f_2)}
\Rey\big((\partial_t+u_0\cdot\nabla)u_1+(u_1\cdot\nabla)u_0\big)-\Delta u_1+\nabla p_1\\[1mm]
\hspace{2cm}\,=\, \Div(\sigma_1[g_2])+\Div\big(\sigma_2[\rho_0,\nabla u_1]+\sigma_2[\rho_1+g_1,\nabla u_0]\big),\\[1mm]
\Delta_ng_2=(\partial_t-\tfrac1 \Pe \Delta_x+u_0\cdot\nabla_x)g_1+P_1^\bot\big(U_0n\cdot\nabla_x(\rho_1+g_1)\big)\\[1mm]
\hspace{4cm}+\Div_n\big(\pi_n^\bot(\nabla u_1)n\rho_0+\pi_n^\bot(\nabla u_0)n(\rho_1+g_1)\big),\\[1mm]
(\partial_t-\tfrac1 \Pe \Delta_x+u_0\cdot\nabla_x)\rho_1+u_1\cdot\nabla_x\rho_0+\langle U_0n\cdot\nabla_xg_1\rangle=0,\\[1mm]
\Div(u_1)=0,~~\langle g_2 \rangle = 0,\\[1mm]
u_1|_{t=0}=0~~\text{if $\Rey\ne0$}, \ \ \int_{\T^d}u_1=0~~\text{if $\Rey=0$},\\[1mm]
\rho_1|_{t=0}=0.
\end{array}\right.
\end{equation}
\item \emph{Order $O(\e^{2})$:} The triplet $(u_2,\rho_2,f_3)$ satisfies
\begin{equation}
\left\{\begin{array}{l}\label{eq :(u_2,f_3)}
\Rey\big((\partial_t+u_0\cdot\nabla)u_2+(u_1\cdot\nabla)u_1+(u_2\cdot\nabla)u_0\big)-\Delta u_2+\nabla p_2\\[1mm]
\hspace{1cm}= \,\Div(\sigma_1[g_3])+\Div\big(\sigma_2[\rho_0,\nabla u_2]+\sigma_2[\rho_1+g_1,\nabla u_1]+\sigma_2[\rho_2+g_2,\nabla u_0]\big),\\[1mm]
\Delta_ng_3=(\partial_t-\tfrac1 \Pe \Delta_x+u_0\cdot\nabla_x)g_2+u_1\cdot\nabla_xg_1+P_1^\bot\big(U_0n\cdot\nabla_x(\rho_2+g_2)\big)\\[1mm]
\hspace{1.5cm} \, + \, \Div_n\big(\pi_n^\bot(\nabla u_2)n\rho_0+\pi_n^\bot(\nabla u_1)n(\rho_1+g_1)+\pi_n^\bot(\nabla u_0)n(\rho_2+g_2)\big),\\[1mm]
(\partial_t-\tfrac1 \Pe \Delta_x+u_0\cdot\nabla_x)\rho_2+u_1\cdot\nabla_x\rho_1+u_2\cdot\nabla_x\rho_0+\langle U_0n\cdot\nabla_xg_2\rangle=0,\\[1mm]
\Div(u_2)=0,~~\langle g_3 \rangle = 0,\\[1mm]
u_2|_{t=0}=0~~\text{if $\Rey\ne0$}, \ \ \int_{\T^d}u_2=0~~\text{if $\Rey=0$},\\[1mm]
\rho_2|_{t=0}=0.
\end{array}\right.
\end{equation}
\end{enumerate}
{We could pursue this hierarchy to higher orders but shall truncate it to second order for shortness in this work --- and making it rigorous is the main technical aspect of the proof of our result.
To ensure the well-posedness of this hierarchy, it is important to notice that it is actually triangular as we can eliminate $g_1,g_2,g_3$ (explicitly!) in terms of the velocity fields $u_0,u_1,u_2$ and of the spatial densities $\rho_0,\rho_1,\rho_2$, cf.~Proposition~\ref{prop:first-terms}.
This leads to a simpler hierarchy of equations for $(u_0,\rho_0)$, $(u_1,\rho_1)$, $(u_2,\rho_2)$.}

Given the above hierarchy, we understand that the well-preparedness condition~\eqref{eq:loc-equi-no-wellprepared} for initial data needs to be further strengthened to indeed avoid initial boundary layers: more precisely, we need to assume that the initial condition $f_\e|_{t=0}=f_\e^\circ$ is compatible with the above hierarchy, meaning that $g_\e$ coincides initially with the expansion $\e g_1+\e^2g_2+\e^3g_3+\ldots$ To accuracy $O(\e^2)$, this leads to the following well-preparedness assumption. Note that similar issues are well known for higher-order hydrodynamic expansions in the Boltzmann theory, see e.g.~\cite{Caflisch-80,seiji1983euler,Lachowicz-87}.

\begin{as}[Well-preparedness]\label{as:well-prepared}
Let $h\in C^\infty(\R^+\times\T^d)^d$, let $\rho^\circ\in H^s(\T^d)\cap\frac1{\omega_d}\Pc(\T^d)$ for some $s\gg1$, and let also $u^\circ\in H^s(\T^d)^d$ with $\Div(u^\circ)=0$ in the Navier--Stokes case.
We assume that the initial condition $f_\e|_{t=0}=f_\e^\circ$ for the Doi--Saintillan--Shelley system~\eqref{eq:eps-kin00} is well-prepared in the following sense:
{there is a constant $C_0<\infty$ such that for all~$0<\e\le\frac{1}{C_0}$,
decomposing}
\[f_\e^\circ(x,n)\,=\,\rho_\e^\circ(x)+g_\e^\circ(x,n),\qquad\rho_\e^\circ\,:=\,\langle f_\e^\circ\rangle\,=\,\fint_{\Sp^{d-1}}f_\e^\circ(\cdot,n)\,\dd n\,\in\,\tfrac1{\omega_d}\Pc(\T^d),\]
we have
\begin{equation*}
\rho_\e^\circ=\rho^\circ
\qquad\text{and}\qquad
\e^{\frac{1}{2}}\|g_\e^\circ-(\e g_1+\e^2 g_2)|_{t=0}
\|_{\Ld^2_{x,n}}\le\,C_0\e^3,
\end{equation*}
where $g_1$ and $g_2$ are the solutions of the hierarchical equations~\eqref{eq :(u_0,f_1)} and~\eqref{eq :(u_1,f_2)} with initial data $\rho^\circ$ in the Stokes case and $(u^\circ,\rho^\circ)$ in the Navier--Stokes case.
\end{as}

\begin{rem}This well-preparedness assumption is compatible with the positivity \mbox{$f_\e|_{t=0} \ge 0$} for~$\e$ small enough, which is necessary to ensure the well-posedness of the Doi--Saintillan--Shelley system~\eqref{eq:eps-kin00}, cf.~Proposition~\ref{prop:well-posedness}.
\end{rem}

\subsection{Well-posedness of ordered fluid models: hierarchical solutions} \label{sec:formal.expansion}
{We turn to our well-posedness result for the ordered fluid equation~\eqref{eq:second-order-fluid-intro-eqn}. Strictly speaking, this equation is well known to be ill-posed whenever $\gamma_1<0$: this is particularly clear at $\Pe<\infty$ as the equation then behaves like a backward heat equation (see also~\cite{Galdi-93} at $\Pe=\infty$, as well as~\cite{Dunn-Rajagopal-95} for a discussion from the physics perspective).
Yet, we can devise a notion of (say, second-order) {\it approximate hierarchical solutions} that satisfy the desired equation up to $O(\e^2)$ and are indeed well-posed.
The definition of hierarchical solutions naturally arises from the formal $\e$-expansion of solutions: equivalently, it coincides with $(u_0+\e u_1,\rho_0+\e\rho_1)$ where $(u_0,\rho_0),(u_1,\rho_1)$ satisfy the above hierarchy~\eqref{eq :(u_0,f_1)}--\eqref{eq :(u_1,f_2)} obtained from the Doi--Saintillan--Shelley model up to explicitly eliminating $g_1,g_2$.
More precisely, this is stated in the following proposition.} The proof is postponed to Appendix~\ref{app.Bous}, {as well as an alternative definition of well-posed perturbative solutions rather inspired by the Boussinesq theory for water waves.} Note that the regularity theory for the system~\eqref{eq:v0rho0} below is quite delicate in the 3D Stokes case: a particularly careful stepwise argument is needed to first cover low-regularity situations.

\begin{prop}[Approximate hierarchical solutions]\label{prop:hier-inhom}
Consider equation~\eqref{eq:second-order-fluid-intro-eqn} with parameters $\eta_0,\Pe>0$, $\eta_1,\mu_0\ge0$, and
$\gamma_1,\gamma_2\in\R$.
If $(u_0,\rho_0),(u_1,\rho_1)$ are smooth solutions of the following two auxiliary systems,
\begin{equation} \label{eq:v0rho0}
\left\{\begin{array}{l} 
\Rey(\partial_t+u_0\cdot\nabla)u_0-\Div\big(2(\eta_0+\eta_1\rho_0)\D(u_0)\big)+\nabla p_0=h,\\[1mm]
(\partial_t-\tfrac1\Pe\Delta+u_0\cdot\nabla)\rho_0=0,\\[1mm]
\Div(u_0)=0,\\[1mm]
u_0|_{t=0}=u^\circ\quad\text{if $\Rey\ne0$},\quad\int_{\T^d}u_0=0\quad\text{if $\Rey=0$},\\[1mm]
\rho_0|_{t=0}=\rho^\circ,
\end{array}\right.
\end{equation}
\begin{equation} \label{eq:v1rho1}
\left\{\begin{array}{l} 
\Rey\big((\partial_t+u_0\cdot\nabla)u_1+(u_1\cdot\nabla)u_0\big)
-\Div\big(2(\eta_0+\eta_1\rho_0)\D(u_1)\big)
+\nabla p_1\\[0,5mm]
\hspace{3.5cm}=\Div\big(2\eta_1\rho_1\D(u_0)+\gamma_1 A_2'(\rho_0,u_0)+\gamma_2(2\D(u_0))^2\big),\\[1mm]
(\partial_t-\tfrac1\Pe\Delta+u_0\cdot\nabla)\rho_1=\mu_0\Delta\rho_0-u_1\cdot\nabla\rho_0,\\[1mm]
\Div(u_1)=0,\\[1mm]
u_1|_{t=0}=0\quad\text{if $\Rey\ne0$},\quad\int_{\T^d}u_1=0\quad\text{if $\Rey=0$},\\[1mm]
\rho_1|_{t=0}=0,
\end{array}\right.
\end{equation}
then the superposition $(\bar u_\e,\bar\rho_\e)=(u_0+\e u_1,\rho_0+\e\rho_1)$ satisfies equation~\eqref{eq:second-order-fluid-intro-eqn} with some controlled remainder $O(\e^2)$ and with initial condition $\bar \rho_\e|_{t=0}=\rho^\circ$ (and $\bar u_\e|_{t=0}=u^\circ$ if~\mbox{$\Rey\ne0$}). For the well-posedness of~\eqref{eq:v0rho0} and~\eqref{eq:v1rho1}, we separately consider the Stokes and Navier--Stokes cases:
{\begin{enumerate}[(i)]
\item \emph{Stokes case $\Rey=0$, $d\le3$:}\\
Given $\rho^\circ \in \Ld^2(\T^d)\cap\frac1{\omega_d}\Pc(\T^d)$ and $h\in \Ld^{\infty}_{\loc}(\R^+;H^{-1}(\T^d)^d)$, there exists a unique global solution~$(u_0,\rho_0)$ of~\eqref{eq:v0rho0} with
\begin{eqnarray*}
u_0&\in&\Ld^\infty_\loc(\R^+;H^{1}(\T^d)^d),\\
\rho_0 &\in& \Ld^\infty_\loc\big(\R^+;\Ld^2(\T^d)\cap\tfrac1{\omega_d}\Pc(\T^d)\big)\cap \Ld^2_\loc(\R^+;H^{1}(\T^d)).
\end{eqnarray*}
Moreover, for all integers $s\ge\frac d2+1$, provided that~$\rho^\circ \in H^s(\T^d)$ and that $h$ belongs to $\Ld^{\infty}_{\loc}(\R^+; H^{s-1} ( \R^+ \times \T^d)^d)$, this solution further satisfies
\begin{eqnarray*}
u_0&\in&\Ld^\infty_\loc(\R^+;H^{s+1}(\T^d)^d),\\
\rho_0 &\in& \Ld^\infty_\loc(\R^+;H^s(\T^d))\cap\Ld^2_\loc(\R^+;H^{s+1}(\T^d)).
\end{eqnarray*}
In that case, if furthermore $h \in W^{1,\infty}_{\loc}(\R^+;H^{s-3}(\T^d)^d)$, there exists a unique global solution $(u_1,\rho_1)$ of~\eqref{eq:v1rho1} with
\begin{eqnarray*}
u_1&\in&\Ld^\infty_\loc(\R^+;H^{s-1}(\T^d)^d),\\
\rho_1 &\in& \Ld^\infty_\loc(\R^+;H^{s-2}(\T^d))\cap\Ld^2_\loc(\R^+;H^{s-1}(\T^d)).
\end{eqnarray*}
\item \emph{Navier--Stokes case $\Rey=1$, $d=2$:}\\
Given {$\rho^\circ\in\Ld^2(\T^2)\cap\frac1{2\pi}\Pc(\T^2)$,} $u^\circ \in \Ld^2(\T^2)^2$ which satisfies $\Div(u^\circ)=0$, and $h\in \Ld^{2}_{\loc}(\R^+; H^{-1} ( \T^2)^2)$, there exists a unique global solution $(u_0,\rho_0)$ of~\eqref{eq:v0rho0} with
\begin{eqnarray*}
u_0&\in&\Ld^\infty_\loc(\R^+;\Ld^2(\T^2)^2)\cap\Ld^2_\loc(\R^+;H^1(\T^2)^2),\\
\rho_0 &\in& {\Ld^\infty_\loc\big(\R^+;\Ld^2(\T^2)\cap\tfrac1{2\pi}\Pc(\T^2)\big)\cap\Ld^2_\loc(\R^+;H^1(\T^2)).}
\end{eqnarray*}
Moreover, for all integers $s\ge2$, provided that $\rho^\circ \in H^s(\T^2)$, $u^\circ \in H^s(\T^2)^2$, and $h\in \Ld^{2}_{\loc}(\R^+; H^{s-1} ( \T^2)^2)$, this solution further satisfies
\begin{eqnarray*}
u_0&\in&\Ld^\infty_\loc(\R^+;H^{s}(\T^2)^2)\cap\Ld^2_\loc(\R^+;H^{s+1}(\T^2)^2),\\
\rho_0 &\in& \Ld^\infty_\loc(\R^+;H^s(\T^2))\cap\Ld^2_\loc(\R^+;H^{s+1}(\T^2)).
\end{eqnarray*}
In that case, there exists a unique global solution $(u_1,\rho_1)$ of~\eqref{eq:v1rho1} with
\begin{align*}
u_1&~~\in~~\Ld^\infty_\loc(\R^+;H^{s-2}(\T^2)^2)\cap\Ld^2_\loc(\R^+;H^{s-1}(\T^2)^2),\\
\rho_1 &~~\in~~ \Ld^\infty_\loc(\R^+;H^{s-2}(\T^2))\cap\Ld^2_\loc(\R^+;H^{s-1}(\T^2)).\qedhere
\end{align*}
\end{enumerate}}
\end{prop}

\subsection{Hydrodynamic limit} \label{sec:main.result}
{At last, we turn to the precise statement of our main result, that is, the rigorous hydrodynamic approximation of the Doi--Saintillan--Shelley theory in the small-$\Wi$ regime, {as already anticipated in Theorem~\ref{thm:informal}.} For shortness, we focus on the {\it first-order} approximation and the emergence of {\it second-order} fluid equations: we shall derive~\eqref{eq:second-order-fluid-intro-eqn} with explicit coefficients
\begin{align}\label{eq:parameters-2nd}
\begin{split}
\eta_0= 1, \qquad \eta_1= \lambda\tfrac{(\theta+2)\omega_d}{2d(d+2)}, \qquad \mu_0 = \tfrac{1}{d(d-1)}U_0^2 , \\[1mm]
\gamma_1 = -\lambda\theta\tfrac{\omega_d}{4d^2(d+2)} , \qquad  \gamma_2 = \lambda\tfrac{\omega_d}{2d^2(d+4)}(\theta+\tfrac{2d}{d+2}),
\end{split}
\end{align}
where we recall that $\lambda,\theta,U_0$ are parameters from the Doi--Saintillan--Shelley system~\eqref{eq:eps-kin00}.
{The proof is given in Section~\ref{sec:expansion} and a physical interpretation of parameters~\eqref{eq:parameters-2nd} is postponed to Section~\ref{secrem:non-newtonian actif}.
Note that the same analysis could be pursued to arbitrary order: we refer in particular to Section~\ref{secrem:third-order} for the next-order description and its physical implications.}

\begin{theor}[Small-$\Wi$ expansion]\label{th:main}
Let $h\in C^\infty(\R^+\times\T^d)^d$, let $\rho^\circ\in H^s(\T^d)\cap\frac1{\omega_d}\Pc(\T^d)$ for some~$s\gg1$,
and let also $u^\circ \in H^s(\T^d)^d$ with $\Div(u^\circ)=0$ in the Navier--Stokes case.
Denote by $(u_\e, f_\e)$ the global solution of the Doi--Saintillan--Shelley model~\eqref{eq:eps-kin00} as given by Proposition~\ref{prop:well-posedness} with initial condition $f_\e|_{t=0}=f_\e^\circ\in\Ld^2\cap\Pc(\T^d\times\S^{d-1})$
satisfying the well-preparedness of Assumption~\ref{as:well-prepared}.
There is some universal constant $C_0<\infty$ such that the following holds provided that
\begin{equation}\label{eq:smallness-cond-RE}
\lambda\theta (1+\Pe)\|\rho^\circ\|_{\Ld^\infty_x} \, \le \, \tfrac1{C_0}.
\end{equation}
\begin{enumerate}[(i)]
\item \emph{Stokes case $\Rey=0$, $d\le3$:}\\
Let $(\bar u_\e, \bar \rho_\e)$ be the unique global approximate hierarchical solution of~\eqref{eq:second-order-fluid-intro-eqn} as given by Proposition~\ref{prop:hier-inhom}(i) with explicit coefficients~\eqref{eq:parameters-2nd}.
Then we have for all $t \ge0$,
\begin{eqnarray*}
\|\nabla(u_\e-\bar u_\e)\|_{\Ld^2_t\Ld^2_x}&\lesssim&\e^2,\\
\| \rho_\e - \bar \rho_\e\|_{\Ld^\infty_t\Ld^2_x} + \|\nabla( \rho_\e - \bar \rho_\e)\|_{\Ld^2_t\Ld^2_x} &\lesssim& \e^2,
\end{eqnarray*}
where multiplicative constants depend on $\Pe$ and on an upper bound on $t,\lambda,U_0,$ and on controlled norms of $h$ and $\rho^\circ$.
\smallskip
\item \emph{Navier--Stokes case $\Rey=0$, $d=2$:}\\
Let $(\bar u_\e, \bar \rho_\e)$ be the unique global approximate hierarchical solution of~\eqref{eq:second-order-fluid-intro-eqn} as given by Proposition~\ref{prop:hier-inhom}(ii) with explicit coefficients~\eqref{eq:parameters-2nd}.
Then we have for all $t \ge0$,
\begin{align*}
\|u_\e-\bar u_\e\|_{\Ld^\infty_t\Ld^2_x}+\|\nabla(u_\e-\bar u_\e)\|_{\Ld^2_t\Ld^2_x}
&~\lesssim~\e^2,\\
\| \rho_\e - \bar \rho_\e\|_{\Ld^\infty_t\Ld^2_x} + \|\nabla( \rho_\e - \bar \rho_\e)\|_{\Ld^2_t\Ld^2_x}
&~\lesssim~ \e^2,
\end{align*}
where multiplicative constants depend on $\Pe$ and on an upper bound on $t,\lambda,U_0,$ and on controlled norms of $h, \rho^\circ,$ and $u^\circ$.
\qedhere
\end{enumerate}
\end{theor}

%%%%%%%%%%%%%%%%%%%%%%%%
%%%%%%%%%%%%%%%%%%%%%%%%

\section{Proof of the main result}\label{sec:expansion}

{This section is devoted to the proof of Theorem~\ref{th:main}, which amounts to the justification of the formal $\e$-expansion described in Section~\ref{sec:Hilbert-formal}.
In a nutshell, the proof relies on the weak-strong principle that we have established for the Doi--Saintillan--Shelley model, cf.\@ Proposition~\ref{prop:well-posedness}, which we now adapt to derive a stability result for approximate solutions obtained by truncating the formal $\e$-expansion.
Although not surprising in itself, this argument requires a number of non-trivial computations and careful developments.

We split the proof into three main parts and start by checking the well-posedness and the propagation of regularity for the hierarchy~\eqref{eq :(u_0,f_1)}--\eqref{eq :(u_2,f_3)} formally derived as asymptotic $\e$-expansion of the Doi--Saintillan--Shelley system.
This is stated as the following proposition. We emphasize that we were not able to find references for some of the systems that appear in this hierarchy, and we believe that some of these new well-posedness results are of independent interest (see in particular the system~\eqref{eq :(u_0,f_1)-Re0-rewr} below).
As already noticed above, the hierarchy~\eqref{eq :(u_0,f_1)}--\eqref{eq :(u_2,f_3)} is triangular as we can eliminate $g_1,g_2,g_3$ in terms of velocity fields and spatial densities, and this elimination is actually made explicitly below, cf.~\eqref{f_1.explicit}--\eqref{f_2.explicit}.}
We focus on the Stokes case $\Rey=0$, while the 2D Navier--Stokes case follows up to straightforward adaptations and is omitted for shortness.
The proof is given in Section~\ref{sec:lem:first-terms}. 

{\begin{prop}[Well-posedness of hierarchy]
\label{prop:first-terms}
Consider the Stokes case $\Rey=0$, $d\le3$.
\begin{enumerate}[(i)]
\item \emph{Well-posedness for $(u_0,\rho_0,g_1)$:}\\
Given integer $s \geq \tfrac{d}{2}+1$, $\rho^\circ \in H^s(\T^d)\cap\frac1{\omega_d} \Pc(\T^d)$, and $h \in \Ld^\infty_{\loc}(\R^+;H^{s-1}(\T^d)^d)$, there exists a unique solution $(u_0,\rho_0,g_1)$ of the Cauchy problem~\eqref{eq :(u_0,f_1)} with
\begin{eqnarray*}
    u_0 &\in& \Ld^\infty_\loc(\R^+;H^{s+1}(\T^d)^d), \\
    \rho_0 &\in& \Ld^\infty_\loc\big(\R^+; H^s(\T^d)\cap\tfrac1{\omega_d}\Pc(\T^d)\big) \cap\Ld^2_\loc(\R^+;H^{s+1}(\T^d)),\\
    g_1 &\in& \Ld^\infty_\loc(\R^+;H^{s-1}(\T^d\times \Sp^{d-1})),
\end{eqnarray*}
and $g_1$ is given by the explicit formula
\begin{align}\label{f_1.explicit}
g_1(\cdot,n)
= -\tfrac{1}{d-1}U_0n\cdot\nabla\rho_0+\tfrac12(n \otimes n):\rho_0\D(u_0).
\end{align}
Moreover, for all $r\ge0$, provided that $s$ is chosen large enough and that $h\in W^{2,\infty}(\R^+;H^{s-1}(\T^d)^d)$, we also have
\begin{eqnarray*}
\partial_t u_0, \,\partial_t^2 u_0 &\in& \Ld^\infty_{\loc}(\R^+; H^{r}(\T^d)^d),\\
\partial_t \rho_0, \,\partial_t^2 \rho_0 &\in& \Ld^\infty_{\loc}(\R^+; H^{r}(\T^d)).
\end{eqnarray*}
\item \emph{Well-posedness for $(u_1,\rho_1,g_2)$:}\\
Given $s\ge0$, provided that the solution $(u_0,\rho_0)$ of item~(i) is such that
\begin{eqnarray*}
u_0&\in&\Ld^\infty_\loc(\R^+;H^{s+3}(\T^d)^d)\cap W^{1,\infty}_\loc(\R^+;H^{s+1}(\T^d)^d),\\
\rho_0&\in&\Ld^\infty_\loc(\R^+;H^{s+2}(\T^d))\cap W^{1,\infty}_\loc(\R^+;H^s(\T^d)),
\end{eqnarray*}
there exists a unique solution $(u_1,\rho_1,g_2)$ of the Cauchy problem~\eqref{eq :(u_1,f_2)} with
\begin{eqnarray*}
u_1&\in& \Ld^\infty_\loc(\R^+;H^{s+1}(\T^d)^d), \\
\rho_1 &\in& \Ld^\infty_\loc(\R^+; H^s(\T^d)) \cap\Ld^2_\loc(\R^+;H^{s+1}(\T^d)),\\
g_2 &\in& \Ld^\infty_\loc(\R^+;H^s(\T^d\times \Sp^{d-1})),
\end{eqnarray*} 
and $g_2$ is given by the explicit formula
\begingroup\allowdisplaybreaks
\begin{eqnarray}\label{f_2.explicit}
\qquad g_2(\cdot,n) &=&
-\tfrac1{d-1}U_0 n \cdot\nabla\rho_1
-\tfrac1{2d}(n \otimes n-\tfrac1d\Id):\Big(\tfrac14A_2'(u_0,\rho_0)-\rho_0\D(u_0)^2\nonumber\\
&&\hspace{4cm} -d\rho_0\D( u_1)-d\rho_1\D(u_0) -\tfrac{1}{d-1}U_0^2\nabla^2\rho_0\Big)\nonumber\\
&& -\tfrac{1}{3(d-1)}U_0(\nabla\rho_0)\cdot\Big(n\big((n\otimes n):\D(u_0)\big)
+\tfrac{4}{d-1}\D(u_0)n\Big)\nonumber\\
&& -\tfrac1{6(d+1)}U_0\rho_0\,\Div_x\Big(n\big((n \otimes n):\D(u_0)\big)
+\tfrac4{d-1}\D(u_0)n\Big)\nonumber\\
&& +\tfrac{1}8\rho_0\Big(\big((n\otimes n):\D(u_0)\big)^2-\tfrac2{d(d+2)}\tr(\D(u_0)^2)\Big),
\end{eqnarray}
\endgroup
in terms of the (non-standard) Rivlin--Ericksen tensor $A_2'$ defined in \eqref{2ndordered.fluid-inhom-ex/def-A2'}. Moreover, for all $r\ge0$, provided that $s$ is chosen large enough, we also have
\begin{eqnarray*}
\partial_t u_1 &\in& \Ld^\infty_{\loc}(\R^+; H^{r}(\T^d)^d),\\
\partial_t \rho_1 &\in& \Ld^\infty_{\loc}(\R^+; H^{r}(\T^d)).
\end{eqnarray*}
\item \emph{Well-posedness for $(u_2,\rho_2,g_3)$:}\\
Given $s\ge0$, provided that the solutions $(u_0,\rho_0)$ and $(u_1,\rho_1)$ of items~(i) and~(ii) are such that
\begin{eqnarray*}
(u_0,\rho_0)&\in&W^{2,\infty}_\loc(\R^+;H^{r}(\T^d)^{d+1}),\\
(u_1,\rho_1)&\in&W^{1,\infty}_\loc(\R^+;H^{r}(\T^d)^{d+1}),
\end{eqnarray*}
for some $r$ large enough, then there exists a unique solution $(u_2,\rho_2,g_3)$ of the Cauchy problem \eqref{eq :(u_2,f_3)} with
\begin{align*}
u_2 &~~\in~~ \Ld^\infty_\loc(\R^+;H^{s+1}(\T^d)^d),\\
\rho_2 &~~\in~~ \Ld^\infty_\loc(\R^+; H^s(\T^d))\cap\Ld^2_\loc(\R^+; H^{s+1}(\T^d)),\\
g_3 &~~\in~~ \Ld^2_\loc(\R^+;H^{s-1}(\T^d\times \Sp^{d-1})).
\end{align*}
\end{enumerate}
Associated to all the above well-posedness results are estimates of the corresponding norms of the solutions in terms of all the parameters and of the controlled norms of the data.
\end{prop}
}

With the above construction of the hierarchy $\{u_n,\rho_n,g_{n+1}\}_{n\ge0}$,
we can now turn to the justification of the formal expansion~\eqref{eq:formal-exp}.
The proof is given in Section~\ref{sec:prop:main}.
{Note that the well-preparedness assumption~\eqref{eq:well-prepared^4} below for initial data is one order stronger than in Assumption~\ref{as:well-prepared}: indeed, while our main result in Section~\ref{sec:main.result} focusses on $O(\e)$ effects, only deriving second-order fluid models, the present result further describes~$O(\e^2)$ effects and therefore requires this strengthened well-preparedness condition. Although not needed for the purposes of our main result, the present next-order analysis is included to illustrate how the $\e$-expansion can be pursued to arbitrary order without additional mathematical difficulties, then leading to higher-order fluid models; we refer to Appendix~\ref{app:3rd-order} for the corresponding derivation of third-order fluid models.}

\begin{prop}[Error estimates for $\e$-expansion]\label{prop:main}
Let $h\in C^\infty(\R^+\times\T^d)^d$, let $\rho^\circ\in H^s(\T^d)\cap\frac1{\omega_d}\Pc(\T^d)$ for some~$s\gg1$,
and let also $u^\circ \in H^s(\T^d)^d$ with $\Div(u^\circ)=0$ in the Navier--Stokes case.
Denote by $(u_\e,f_\e)$ the solution of the Doi--Saintillan--Shelley model~\eqref{eq:eps-kin00} as given by Proposition~\ref{prop:well-posedness},
and assume that the initial condition $f_\e|_{t=0}=f_\e^\circ$ is well-prepared in the following sense: {there is a constant $C_0<\infty$ such that for all $0<\e\le\frac1{C_0}$, decomposing $f_\e^\circ=\rho_\e^\circ+g_\e^\circ$ with $\rho_\e^\circ:=\langle f_\e^\circ\rangle$, we have
\begin{equation}
{\rho_\e^\circ\,=\,\rho^\circ}
\qquad\text{and}\qquad
\e^\frac12\|g_\e^\circ-(\e g_1+\e^2 g_2+\e^3g_3)|_{t=0}\|_{\Ld^2_{x,n}}\,\le\,C_0\e^4,\label{eq:well-prepared^4}
\end{equation}
where $g_1,g_2,g_3$ are the solutions of hierarchical equations~\eqref{eq :(u_0,f_1)}--\eqref{eq :(u_2,f_3)}
as defined in Proposition~\ref{prop:first-terms} with data~$h,u^\circ,\rho^\circ$.
There is some universal constant $C_1<\infty$ such that the following holds provided that}
\begin{equation*}
\lambda\theta (1+\Pe)\|\rho^\circ\|_{\Ld^\infty_{x}} \, \le \, \tfrac1{C_1}.
\end{equation*}
\begin{enumerate}[(i)]
\item \emph{Stokes case $\Rey=0$, $d\le3$:}\\
For all $t\ge0$, we have
\begin{eqnarray*}
\hspace{1cm}\e^\frac12\|\nabla (u_\e-u_0-\e u_1-\e^2u_2)\|_{\Ld^\infty_t\Ld^2_{x}}
+\|\nabla (u_\e-u_0-\e u_1-\e^2u_2)\|_{\Ld^2_t\Ld^2_x}&\!\!\le\!\!&\Cc(t)\e^3,\\
\hspace{1cm}\e^\frac12\|g_\e-\e g_1-\e^2g_2-\e^3 g_3\|_{\Ld^\infty_t\Ld^2_{x,n}}\!+\|\nabla_n(g_\e-\e g_1-\e^2g_2-\e^3 g_3)\|_{\Ld^2_t\Ld^2_{x,n}}&\!\!\le\!\!&\Cc(t)\e^4,\\
\hspace{1cm}\| \rho_\e - \rho_0 - \e \rho_1 - \e^2 \rho_2 \|_{\Ld^\infty_t \Ld^2_x} + \|  \nabla(\rho_\e - \rho_0 - \e \rho_1 - \e^2 \rho_2)\|_{\Ld^2_t\Ld^2_x} & \!\!\le\!\!& \Cc(t) \e^3,
\end{eqnarray*}
{provided that $\e\Cc(t)\le1$,}
where the multiplicative constant $\Cc(t)$ depends on $
\Pe$ and on an upper bound on
\begin{gather*}
\qquad t,~{C_0,~C_1,~\lambda,~U_0,}\quad
\|(\nabla u_0,\nabla u_1,\nabla u_2)\|_{\Ld^\infty_t\Ld^\infty_x},\quad
\|(\rho_0, \rho_1, \rho_2)\|_{\Ld^\infty_t\Ld^\infty_x},\\
\qquad
\|(g_1,g_2,g_3)\|_{\Ld^\infty_tW^{1,\infty}_{x,n}},\quad
\|(\partial_t-\Delta_x)g_3\|_{\Ld^\infty_t\Ld^2_{x,n}}.
\end{gather*}
\item \emph{Navier--Stokes case $\Rey=0$, $d=2$:}\\
For all $t\ge0$, we have
\begin{eqnarray*}
\hspace{1cm}\|u_\e-u_0-\e u_1-\e^2u_2\|_{\Ld^\infty_t\Ld^2_x}+\|\nabla (u_\e-u_0-\e u_1-\e^2u_2)\|_{\Ld^2_t\Ld^2_x}&\!\!\le\!\!&\Cc(t)\e^3,\\
\hspace{1cm}\e^{\frac12}\|g_\e-\e g_1-\e^2g_2-\e^3 g_3\|_{\Ld^\infty_t\Ld^2_{x,n}}\!
+
\|\nabla_n(g_\e-\e g_1-\e^2g_2-\e^3 g_3)\|_{\Ld^2_t\Ld^2_{x,n}}&\!\!\le\!\!&\Cc(t)\e^4,\\
\hspace{1cm}\| \rho_\e - \rho_0 - \e \rho_1 - \e^2 \rho_2 \|_{\Ld^\infty_t \Ld^2_x}  + \|  \nabla(\rho_\e - \rho_0 - \e \rho_1 - \e^2 \rho_2)\|_{\Ld^2_t\Ld^2_x} & \!\!\le\!\! & \Cc(t) \e^3,
\end{eqnarray*}
{provided that $\e\Cc(t)\le1$,}
where $\Cc(t)$ now depends $
\Pe$ and on an upper bound on
\begin{gather*}
\qquad
t,~{C_0,~C_1,~\lambda,~U_0,}\quad
\|(u_0,u_1,u_2)\|_{\Ld^\infty_tW^{1,\infty}_x},
\quad\|(\rho_0, \rho_1, \rho_2)\|_{\Ld^\infty_t\Ld^\infty_x},\\
\qquad
\|(g_1,g_2,g_3)\|_{\Ld^\infty_tW^{1,\infty}_{x,n}},
\quad\|(\partial_t-\Delta_x)g_3\|_{\Ld^\infty_t\Ld^2_{x,n}}.
\qedhere
\end{gather*}
\end{enumerate}
\end{prop}

Finally, to conclude the proof of Theorem~\ref{th:main}, it remains to identify the equation satisfied by the truncated $\e$-expansion $(\bar u_\e,\bar \rho_\e):=(u_0+\e u_1,\rho_0 + \e \rho_1)$: we show that it precisely coincides with {the so-called approximate hierarchical solution of the second-order fluid model~\eqref{eq:second-order-fluid-intro-eqn} defined in Proposition~\ref{prop:hier-inhom}} for some suitable choice of parameters.
The proof is given in Section~\ref{seclem:Bous}.

\begin{prop}[From hierarchy to second-order fluids]\label{lem:Bous}
Given the solutions $(u_0,\rho_0,g_1)$ and $(u_1,\rho_1,g_2)$ of the hierarchy~\eqref{eq :(u_0,f_1)}--\eqref{eq :(u_1,f_2)}, as constructed in Proposition~\ref{prop:first-terms} in the Stokes case,
the superposition $(\bar u_\e,\bar\rho_\e):=(u_0+\e u_1,\rho_0+\e\rho_1)$ coincides with the unique hierarchical solution of the second-order fluid model~\eqref{eq:second-order-fluid-intro-eqn} in the sense of Proposition~\ref{prop:hier-inhom}, with coefficients $\eta_0, \eta_1, \mu_0, \gamma_1,\gamma_2$ given explicitly by~\eqref{eq:parameters-2nd}.
\end{prop}

The combination of Propositions~\ref{prop:main} and~\ref{lem:Bous} completes the proof of Theorem~\ref{th:main}.
Note however that we only appeal to Assumption~\ref{as:well-prepared} in the statement of Theorem~\ref{th:main}, which is one order weaker than
the well-preparedness assumption~\eqref{eq:well-prepared^4} required in Proposition~\ref{prop:main}.
Indeed, as explained, we focus in Theorem~\ref{th:main} on $O(\e)$ effects, only deriving the second-order fluid model, while in Proposition~\ref{prop:main} we took care to further describe~$O(\e^2)$ effects. For the purposes of Theorem~\ref{th:main},
the well-preparedness assumption~\eqref{eq:well-prepared^4} can therefore simply be replaced by Assumption~\ref{as:well-prepared}.

\subsection{Computational tools for spherical calculus}
{In this section, we briefly recall several computational tools that will be used throughout this work to compute derivatives and integrals on the sphere.}
First, we recall that the Laplace--Beltrami operator $\Delta_n$ on the sphere $\Sp^{d-1}$ ($d \geq 2$) can be computed as follows: given a smooth function $g:\Sp^{d-1}\to\R$, we can extend it to $\R^d\setminus\{0\}$ by setting $G(x):=g \big(\tfrac{x}{|x|} \big)$, and we then have
\begin{equation}\label{eq:form-Delta-n-comput}
\Delta_ng\,=\,(\Delta_xG)|_{\Sp^{d-1}}.
\end{equation}
In particular, we can compute in this way
\begin{eqnarray}
\Delta_n(n_i)&=& (1-d)n_i,\label{eq:comput-triangle}\\
\Delta_n(n_in_j)&=&2\delta_{ij}-2dn_in_j,\nonumber\\
\Delta_n(n_in_jn_k)&=&2(\delta_{ij}n_k+\delta_{ik}n_j+\delta_{jk}n_i)
-3(d+1)n_in_jn_k,\nonumber\\
\hspace{-1cm}\Delta_n(n_in_jn_kn_l)&=&2(\delta_{ij}n_kn_l+\delta_{kj}n_in_l+\delta_{ki}n_jn_l+\delta_{il}n_jn_k+\delta_{lk}n_in_j+\delta_{jl}n_in_k)\nonumber\\
&&-4(d+2)n_in_jn_kn_l,\nonumber
\end{eqnarray}
and so on for higher-order polynomials.
These formulas can be used to explicitly invert $\Delta_n$ on mean-zero polynomial expressions: for any trace-free symmetric matrix $A\in\R^{d\times d}$, we find for instance,
\begin{eqnarray}
\Delta_n^{-1}(n)&=& - \tfrac 1 {d-1}  n,\label{Laplace.1}\\
\Delta_n^{-1}\big(n \otimes n:A\big)&=& - \tfrac 1 {2d}  \,n \otimes n:A,\label{Laplace.2}\\
\Delta_n^{-1}\Big(n(n\otimes n:A)\Big)&=&- \tfrac 1 {3 (d+1)} \Big( n( n \otimes n:A)  + \tfrac{4}{d-1} A n\Big),\label{Laplace.3}\\
\Delta_n^{-1}\left((n \otimes n:A)^2 - \tfrac 2 {d(d+2)} \tr(A^2)\right) &=&- \tfrac  1 {4 (d+2)} \Big( (n \otimes n:A)^2+\tfrac 4 {d }n \otimes n:A^2\hspace{1cm}\nonumber \\
&&\hspace{2.5cm}- \tfrac 2d(\tfrac1{d+2}+\tfrac2d)\, \tr (A^2)\Big).\label{Laplace.4}
\end{eqnarray}
Henceforth, the pseudo-inverse $\Delta_n^{-1}$ is chosen to be defined as an operator from mean-zero fields to mean-zero fields.

We also recall that the divergence of functions on $\Sp^{d-1}$ can be computed similarly as the Laplace--Beltrami operator~\eqref{eq:form-Delta-n-comput} by an extension procedure: for any trace-free matrix $A$ and any smooth function $g:\Sp^{d-1}\to\R$, we find for instance
\begin{align} \label{spherical.divergence}
\Div_n(\pi_n^\bot A ng)\,=\,\nabla_ng\cdot An-d\,(n\otimes n): Ag.
\end{align}

Finally, we further note that the above differential formulas~\eqref{eq:comput-triangle} imply by direct integration the following elementary integral identities for polynomial expressions on the sphere,
\begin{eqnarray}
\int_{\Sp^{d-1}}n_in_j\,\dd n&=&\tfrac{\omega_d}d\delta_{ij},\label{eq:explicit-integral-n}\\
\int_{\Sp^{d-1}}n_in_jn_kn_l\,\dd n&=&\tfrac{\omega_d}{d(d+2)}\big(\delta_{ij}\delta_{kl}+\delta_{kj}\delta_{il}+\delta_{ki}\delta_{jl}\big),\nonumber\\
\int_{\Sp^{d-1}}n_in_jn_kn_ln_mn_p\,\dd n&=&
\tfrac{\omega_d}{d(d+2)(d+4)}\big(\delta_{ij}\delta_{kl}\delta_{mp}+\ldots\big),\nonumber
\end{eqnarray}
and so on,
where we recall the notation $\omega_d=|\Sp^{d-1}|$.
These identities imply in particular, for any trace-free symmetric matrices $A,B\in\R^{d\times d}$,
\begin{eqnarray*}
\int_{\Sp^{d-1}} \big(n \otimes n - \tfrac1d{\Id}\big)\, \dd n &=&0,\\
\int_{\Sp^{d-1}} \big(n \otimes n-\tfrac1d\Id\big)\,(n \otimes n:A) \, \dd n &=&\tfrac{2\omega_d}{d(d+2)}A,\\
\int_{\Sp^{d-1}} \big(n \otimes n-\tfrac1d\Id\big)\,(n \otimes n:A)^2 \, \dd n
&=&\tfrac{8\omega_d}{d(d+2)(d+4)}\big(A^2-\tfrac1d\Tr(A^2)\Id\big).
\end{eqnarray*}

\subsection{Proof of Proposition~\ref{prop:first-terms}}\label{sec:lem:first-terms}
{We focus on items~(i) and~(ii), while the argument for item~(iii) easily follows by adapting that for item~(ii).}
We split the proof into three main steps.

\medskip
\step1 Explicit formulas for $g_1,g_2$.\\
With the above identities~\eqref{Laplace.1}--\eqref{Laplace.4}, we can explicitly solve the successive equations for $g_1$ and $g_2$ in~\eqref{eq :(u_0,f_1)}--\eqref{eq :(u_1,f_2)}. Note that these equations can be solved for fixed~$(x,t)$, thus treating $\rho_0,\rho_1,\D(u_0),\D(u_1)$ and their derivatives as parameters.
We start with the computation of~$g_1$.
Using~\eqref{spherical.divergence}, the defining equation for $g_1$ in~\eqref{eq :(u_0,f_1)} can be rewritten as
\begin{align*}
\Delta_ng_1
= U_0n\cdot\nabla\rho_0- d\,(n \otimes n):\rho_0\D(u_0).
\end{align*}
Hence, by~\eqref{Laplace.1} and~\eqref{Laplace.2},
the explicit form~\eqref{f_1.explicit} for $g_1$ follows.

We turn to the computation of $g_2$. Using~\eqref{spherical.divergence}, inserting the explicit form~\eqref{f_1.explicit} for~$g_1$,
and noting in particular that $\nabla_ng_1=-\tfrac{U_0}{d-1}\pi_n^\bot\nabla\rho_0+\rho_0\,\pi_n^\bot\D(u_0)n$,
the defining equation for $g_2$ in~\eqref{eq :(u_1,f_2)} can be rewritten as
\begin{eqnarray*}
\Delta_n g_2 &=& (\partial_t- \tfrac{1}{\Pe} \Delta_x + u_0 \cdot \nabla_x)g_1+P_1^\bot\big( U_0 n \cdot \nabla_x (\rho_1+g_1)\big)\\
&&-d\,(n\otimes n):\rho_0\D( u_1) -d\,(n\otimes n):(\rho_1+g_1)\D(u_0)+\nabla_ng_1\cdot (\nabla u_0)n\\
&=& (\partial_t- \tfrac{1}{\Pe}\Delta_x + u_0 \cdot \nabla_x)\big(-\tfrac{1}{d-1}U_0n\cdot\nabla\rho_0+\tfrac12(n \otimes n):\rho_0\D(u_0)\big)\\
&&+U_0 n \cdot \big(\nabla\rho_1-\tfrac{1}{d-1}(\nabla u_0)^T\nabla\rho_0\big)+\rho_0(n\otimes n):(\D(u_0)(\nabla u_0))\\
&&-(n\otimes n-\tfrac1d\Id):\big(d\rho_0\D( u_1)+d\rho_1\D(u_0)
+\tfrac{1}{d-1}U_0^2\nabla^2\rho_0\big)\\
&&+\tfrac{d+1}{d-1}U_0(n\cdot\nabla\rho_0)\big((n\otimes n):\D(u_0)\big)+\tfrac12U_0\rho_0 n\cdot\nabla_x\big((n \otimes n):\D(u_0)\big)\\
&&-\tfrac{d+2}2\rho_0\,\big((n\otimes n):\D(u_0)\big)^2.
\end{eqnarray*}
Recalling the following definition of Rivlin--Ericksen tensor,
\[A_2'(u_0,\rho_0)\,=\,\big(\partial_t- \tfrac{1}{\Pe}\Delta_x + u_0 \cdot \nabla_x\big)(2\rho_0\!\D(u_0))+(\nabla u_0)^T(2\rho_0\!\D(u_0))+(2\rho_0\!\D(u_0))(\nabla u_0),\]
and noting that $\tr(A_2'(u_0,\rho_0))=4\rho_0\tr(\D(u_0)^2)$, we can reformulate the above as
\begin{eqnarray*}
\Delta_n g_2 &=&-\tfrac{1}{d-1}U_0n\cdot\Big((\partial_t- \tfrac{1}{\Pe}\Delta + u_0 \cdot \nabla)\nabla\rho_0+(\nabla u_0)^T\nabla\rho_0\Big)+U_0 n \cdot\nabla\rho_1\\
&&+(n \otimes n-\tfrac1d\Id):\Big(\tfrac14A_2'(u_0,\rho_0)-d\rho_0\D( u_1)-d\rho_1\D(u_0)
-\tfrac{1}{d-1}U_0^2\nabla^2\rho_0\Big)\\
&&+\tfrac{d+1}{d-1}U_0(n\cdot\nabla\rho_0)\big((n\otimes n):\D(u_0)\big)+\tfrac12U_0\rho_0 n\cdot\nabla_x\big((n \otimes n):\D(u_0)\big)\\
&&-\tfrac{d+2}2\rho_0\Big(\big((n\otimes n):\D(u_0)\big)^2
-\tfrac2{d(d+2)}\tr(\D(u_0)^2)\Big).
\end{eqnarray*}
The equation in \eqref{eq :(u_0,f_1)} for $\rho_0$ entails that the first right-hand side term vanishes as
\begin{equation}\label{cancel-u_0}
(\partial_t- \tfrac{1}{\Pe}\Delta + u_0 \cdot \nabla)\nabla\rho_0+(\nabla u_0)^T\nabla\rho_0\,=\,\nabla(\partial_t-\tfrac{1}{\Pe}\Delta + u_0 \cdot \nabla)\rho_0\,=\,0.
\end{equation}
Now appealing to~\eqref{Laplace.2},
the second term reads
\begin{multline*}
\lefteqn{(n \otimes n-\tfrac1d\Id):\Big(\tfrac14A_2'(u_0,\rho_0)-d\rho_0\D( u_1)-d\rho_1\D(u_0)
-\tfrac{1}{d-1}U_0^2\nabla^2\rho_0\Big)}\\
\, =\, - \tfrac{1}{2d} \Delta_n \bigg( (n \otimes n):\Big(\tfrac14A_2'(u_0,\rho_0)-d\rho_0\D( u_1)-d\rho_1\D(u_0)
-\tfrac{1}{d-1}U_0^2\nabla^2\rho_0 \Big) \Id \bigg).
\end{multline*}
Using \eqref{Laplace.3}, we further get for the third term,
\begin{multline*}
\tfrac{d+1}{d-1}U_0(n\cdot\nabla\rho_0)\big((n\otimes n):\D(u_0)\big) \\
= - \tfrac{U_0}{3(d-1)}\Delta_n \bigg((\nabla\rho_0)\cdot\Big(n\big((n\otimes n):\D(u_0)\big)
+\tfrac{4}{d-1}\D(u_0)n\Big) \bigg),
\end{multline*}
and for the fourth term,
\begin{align*}
    \tfrac12U_0\rho_0 n\cdot\nabla_x\big((n \otimes n):\D(u_0)\big) = -\tfrac{U_0}{6(d+1)} \Delta_n \bigg(\rho_0\,\Div_x\Big(n\big((n \otimes n):\D(u_0)\big)\bigg).
\end{align*}
Finally, using \eqref{Laplace.4} for the last term, the explicit form~\eqref{f_2.explicit} for $g_2$ follows.

\medskip
\step2 {Proof of item~(i)}.\\
In view of the explicit solution for the equation for $g_1$ obtained in Step~1, cf.~\eqref{f_1.explicit},
the elastic stress $\sigma_1[g_1]$ defined in~\eqref{eq:eps-kin00} takes the form
\begin{equation*}
\sigma_1[g_1]
\,=\,\tfrac12\lambda\theta\int_{\Sp^{d-1}}(n\otimes n)\D(u_0)(n \otimes n)\rho_0\,\dd n\,=\,\tfrac12\theta\,\sigma_2[\rho_0,\nabla u_0],
\end{equation*}
where we have used that integrals of monomials of odd degree on $\Sp^{d-1}$ vanish.
By~\eqref{eq:explicit-integral-n}, this actually means
\begin{equation*}
\sigma_1[g_1]
\,=\,\tfrac{\omega_d}{d(d+2)}\lambda\theta\D(u_0)\rho_0\,=\,\tfrac12\theta\,\sigma_2[\rho_0,\nabla u_0].
\end{equation*}
In the Stokes case $\Rey=0$, the system~\eqref{eq :(u_0,f_1)} then becomes
\begin{equation}
\left\{\begin{array}{l}\label{eq :(u_0,f_1)-Re0-rewr}
-\Div\big((1+c_0\rho_0)2\!\D(u_0)\big)+\nabla p_0=h,\\[1mm]
(\partial_t-\tfrac1 \Pe \Delta+u_0\cdot\nabla)\rho_0=0,\\[1mm]
\Div(u_0)=0,\quad
\int_{\T^d} u_0 = 0,\\[1mm]
\rho_0|_{t=0}=\rho^\circ,
\end{array}\right.
\end{equation}
where we have set for shortness $c_0:=\lambda\tfrac{\omega_d(\theta+2)}{2d(d+2)}$.
Surprisingly, we could not find a reference for this natural system.
We establish the well-posedness and the propagation of regularity for this system, which we believe to be of independent interest, and we split the proof into five further substeps. Note that the propagation of regularity requires some care: we need to treat low and high regularity separately in several steps.

\medskip
\substep{2.1}
Well-posedness of energy solutions for~\eqref{eq :(u_0,f_1)-Re0-rewr}.\\
On the one hand, using the incompressibility constraint, the energy identity for the transport-diffusion equation for $\rho_0$ reads
\begin{equation}\label{eq:rho0-L2-energy}
\|{\rho}_0 \|_{\Ld^\infty_t\Ld^2_x}^2 +  \tfrac{2}{\Pe} \|\nabla {\rho}_0 \|_{\Ld^2_t\Ld^2_x}^2 = \|\rho^\circ\|_{\Ld^2(\T^d)}^2.
\end{equation}
On the other hand, for $\rho_0\ge0$, the energy identity for ${u}_0$ takes the form
\begin{align*}
 \|\nabla {u}_0\|_{\Ld^2(\T^d)}^2\,\le\,2\int_{\T^d}(1+c_0\rho_0)|\!\D(u_0)|^2\,=\, \int_{\T^d}hu_0,
\end{align*}
and thus, using $\int_{\T^d}u_0=0$,
\begin{equation}\label{eq:u0-L2-energy}
\|\nabla{u}_0\|_{\Ld^\infty_t\Ld^2_x} \, \le \, \|h\|_{H^{-1}(\T^d)}.
\end{equation}
By a standard fixed-point approach,
in view of these a priori estimates,
given an initial condition $\rho^\circ\in \Ld^2(\T^d)$ with $\rho^\circ\ge0$ and given $h\in \Ld^\infty_\loc(\R^+;H^{-1}(\T^d)^d)$, we easily check that the system~\eqref{eq :(u_0,f_1)-Re0-rewr} is globally well-posed with $u_0\in \Ld^\infty(\R^+;H^1(\T^d)^d)$ and $\rho_0\in\Ld^\infty(\R^+;\Ld^2(\T^d))\cap\Ld^2(\R^+;H^1(\T^d)^d)$ with $\rho_0\ge0$; we skip the detail.

\medskip
\substep{2.2} $H^1$ regularity for $(\nabla u_0,\rho_0)$: provided that the initial condition further satisfies $\rho^\circ \in H^1(\T^d)$, and provided that $h\in \Ld^\infty_\loc(\R^+;\Ld^2(\T^d)^d)$, we show that $\rho_0 \in \Ld^\infty_\loc(\R^+;H^1(\T^d))\cap\Ld^2_\loc(\R^+;H^2(\T^d))$ and that $u_0 \in \Ld^1_\loc(\R^+;H^2(\T^2)^2)$, with
\begin{align}
\| \nabla \rho_0 \|_{\Ld^\infty_t\Ld^2_x} + \|\nabla^2 \rho_0\|_{\Ld^2_t\Ld^2_x}
\,\le\, C(t,h,\rho^\circ),
\label{eq:nabrho0-L2-energy}\\
\left\{\begin{array}{lll}
\|\Delta u_0\|_{\Ld^2_t\Ld^2_x}
\,\lesssim\, C(t,h,\rho^\circ)&:&\text{if $d=2$},\\
\|\Delta u_0\|_{\Ld^1_t\Ld^2_x}
\,\lesssim\, C(t,h,\rho^\circ)&:&\text{if $d=3$},
\end{array}\right.
\label{eq:bnd-Deltau0-l21t}
\end{align}
where henceforth $C(t,h,\rho^\circ)$ stands for a constant further depending on an upper bound on~$t,c_0$, and on the controlled norms of the data $h$ and $\rho^\circ$.

We start with the proof of~\eqref{eq:nabrho0-L2-energy}.
Testing the equation for $\rho_0$ with $\Delta \rho_0$, and using the incompressibility of $u_0$, we find
\begin{eqnarray*}
\tfrac12\tfrac{\mathrm{d}}{\mathrm{d}t} \| \nabla \rho_0 \|_{\Ld^2_x}^2 +  \tfrac{1}{\Pe} \|\Delta \rho_0\|_{\Ld^2_x}^2
&=&\int_{\T^d} (u_0 \cdot \nabla \rho_0) \Delta \rho_0\\
&=&-\int_{\T^d} \nabla u_0:(\nabla \rho_0\otimes\nabla \rho_0)\\
&\le&\|\nabla u_0\|_{\Ld^2_x}\|\nabla \rho_0\|_{\Ld^4_x}^2.
\end{eqnarray*}
In dimension $d<4$, we can appeal to the Gagliardo--Nirenberg interpolation inequality
\[\|\nabla\rho_0\|_{\Ld^4_x} \,\lesssim \, \|\nabla\rho_0\|_{\Ld^2_x}^{1-\frac{d}{4}} \|\nabla^2\rho_0\|_{\Ld^2_x}^{\frac{d}{4}}.\]
The above then becomes, further using the elliptic estimate $\|\nabla^2 \rho_0\|_{\Ld^2_x}^2 \lesssim \|\Delta \rho_0\|_{\Ld^2_x}^2$,
\begin{equation*}
\tfrac12\tfrac{\mathrm{d}}{\mathrm{d}t} \| \nabla \rho_0 \|_{\Ld^2_x}^2 +  \tfrac{1}{\Pe} \|\Delta \rho_0\|_{\Ld^2_x}^2
\, \lesssim  \, \|\nabla u_0\|_{\Ld^2_x} \|\nabla \rho_0\|_{\Ld^2_x}^{2(1-\frac{d}{4})} \|\Delta\rho_0\|_{\Ld^2_x}^{\frac{d}{2}},
\end{equation*}
where the last factor can now be absorbed in the left-hand side by Young's inequality, to the effect of
\begin{equation*}
\tfrac{\mathrm{d}}{\mathrm{d}t} \| \nabla \rho_0 \|_{\Ld^2_x}^2 +  \tfrac{1}{\Pe} \|\Delta \rho_0\|_{\Ld^2_x}^2
\, \lesssim  \, \|\nabla u_0\|_{\Ld^2_x}^{\frac4{4-d}} \|\nabla \rho_0\|_{\Ld^2_x}^{2}.
\end{equation*}
By Gr\"onwall's inequality and the a priori energy estimate~\eqref{eq:u0-L2-energy}, this precisely proves the claim~\eqref{eq:nabrho0-L2-energy}.

Next, we turn to the proof of the corresponding estimate~\eqref{eq:bnd-Deltau0-l21t} for $u_0$. Testing the equation for $u_0$ with $\Delta u_0$, we find
\begin{eqnarray*}
\|\Delta u_0\|_{\Ld^2_x}^2
&=&-\int_{\T^d} h \cdot \Delta u_0 - 2c_0\int_{\T^d} \Div(\rho_0\D(u_0)) \cdot \Delta u_0\\
&=&-\int_{\T^d} h \cdot \Delta u_0 - c_0\int_{\T^d} \rho_0|\Delta u_0|^2- 2c_0\int_{\T^d} (\D(u_0)\nabla\rho_0 )\cdot \Delta u_0,
\end{eqnarray*}
and thus, by Young's inequality,
\begin{equation*}
\|\Delta u_0\|_{\Ld^2_x}
\,\lesssim\,\|h\|_{\Ld^2_x}+c_0\|\!\D(u_0)\nabla\rho_0\|_{\Ld^2_x}.
\end{equation*}
In dimension $d<4$, appealing to the Gagliardo--Nirenberg interpolation inequality, we can estimate for all $2< p,q<\infty$ with $\frac1p+\frac1q=\frac12$ and $d\le p,q\le\frac{2d}{d-2}$,
\begin{eqnarray*}
\|\!\D(u_0)\nabla\rho_0\|_{\Ld^2_x}
&\le&\|\nabla u_0\|_{\Ld^p_x}\|\nabla\rho_0\|_{\Ld^q_x}\\
&\lesssim_{p,q}&\|\nabla u_0\|_{\Ld^2_x}^{1-\frac dq}\|\nabla^2 u_0\|_{\Ld^2_x}^{\frac dq}\|\nabla\rho_0\|_{\Ld^2_x}^{1-\frac dp}\|\nabla^2\rho_0\|_{\Ld^2_x}^{\frac dp}.
\end{eqnarray*}
Inserting the a priori estimates~\eqref{eq:u0-L2-energy} and~\eqref{eq:nabrho0-L2-energy}, combining with the above, and using Young's inequality to absorb the $H^2$ norm of $u_0$, we are led to
\begin{equation*}
\|\Delta u_0\|_{\Ld^2_x}
\,\lesssim_{p,q}\,\|h\|_{\Ld^2_x}+C(t,h)\|\nabla\rho^\circ\|_{\Ld^2_x}^{\frac{q}{q-d}-\frac{d}{2}\frac{q-2}{q-d}}\|\nabla^2\rho_0\|_{\Ld^2_x}^{\frac d2\frac{q-2}{q-d}}.
\end{equation*}
Choosing any $q<\infty$ if $d=2$, and choosing $q=\frac{2d}{d-2}$ if $d=3$, we get
\begin{equation}\label{eq:bnd-Deltau0-l21t-pre}
\|\Delta u_0\|_{\Ld^2_x}
\,\lesssim_\gamma\,\|h\|_{\Ld^2_x}+C(t,h)\|\nabla\rho^\circ\|_{\Ld^2_x}^{\gamma}\|\nabla^2\rho_0\|_{\Ld^2_x}^{\frac{2}{4-d}},
\end{equation}
where we can take any exponent $\gamma>0$ if $d=2$, and $\gamma=0$ if $d=3$.
Combined with~\eqref{eq:nabrho0-L2-energy}, this yields the claim~\eqref{eq:bnd-Deltau0-l21t}.

\medskip
\substep{2.3} $H^2$ regularity for $(\nabla u_0,\rho_0)$: provided that the initial condition further satisfies $\rho^\circ \in H^2(\T^d)$, and provided that $h\in\Ld^\infty_\loc(\R^+; H^1(\T^d)^d)$, we show that $\rho_0 \in \Ld^\infty_\loc(\R^+;H^2(\T^d))\cap\Ld^2_\loc(\R^+;H^3(\T^d))$ and that $u_0 \in \Ld^\infty_\loc(\R^+;H^2(\T^d)^d)\cap\Ld^1_\loc(\R^+;H^3(\T^d)^d)$, with
\begin{gather}
\|\nabla^2\rho_0\|_{\Ld^\infty_t\Ld^2_x}
+\|\nabla^3\rho_0\|_{\Ld^2_t\Ld^2_x}
+\|\nabla^2 u_0\|_{\Ld^\infty_t\Ld^2(\T^d)}
\,\le\,C(t,h,\rho^\circ),
\label{eq:nabrho0-H2Linfty-est}\\
\left\{\begin{array}{lll}
\|\nabla^3u_0\|_{\Ld^\infty_t\Ld^2_x}^2\,\le\, C(t,h,\rho^\circ)&:&\text{if $d=2$},\\
\|\nabla^3u_0\|_{\Ld^4_t\Ld^2_x}^2\,\le\, C(t,h,\rho^\circ)&:&\text{if $d=3$}.
\end{array}\right.
\label{eq:nabu0-H3Linfty-est}
\end{gather}
We start with the proof of~\eqref{eq:nabrho0-H2Linfty-est}.
Applying $\nabla^2$ to both sides of the equation for $\rho_0$, testing with $\nabla^2\rho_0$, and integrating by parts, we find
\begin{eqnarray*}
\tfrac12\tfrac{\dd}{\dd t}\|\nabla^2\rho_0\|_{\Ld^2_x}^2+\|\nabla^3\rho_0\|_{\Ld^2_x}^2
&=&-\int_{\T^d}(\nabla^2_{ij}\rho_0)\nabla^2_{ij}\big((u_0)_k\nabla_k\rho_0\big)\\
&=&\int_{\T^d}(\nabla^3\rho_0)_{ijk}(\nabla u_0)_{ki}(\nabla\rho_0)_j
+\int_{\T^d}(\nabla^3\rho_0)_{ijj}(\nabla u_0)_{ki}(\nabla\rho_0)_k,
\end{eqnarray*}
and thus, using Young's inequality to absorb the factors $\nabla^3\rho_0$ in the right-hand side,
\begin{equation*}
\partial_t\|\nabla^2\rho_0\|_{\Ld^2_x}^2+\|\nabla^3\rho_0\|_{\Ld^2_x}^2
\,\lesssim\,\|\nabla u_0 \nabla\rho_0\|_{\Ld^2_x}^2.
\end{equation*}
By the Gagliardo--Nirenberg interpolation inequality, we can estimate for all $2<p,q<\infty$ with $\frac1p+\frac1q=\frac12$ and $d\le p,q\le\frac{2d}{d-2}$,
\begin{eqnarray*}
\|\nabla u_0 \nabla\rho_0\|_{\Ld^2_x}
&\le&\|\nabla u_0\|_{\Ld^p_x}\|\nabla\rho_0\|_{\Ld^q_x}\\
&\le&\|\nabla u_0\|_{\Ld^2_x}^{1-\frac dq}
\|\nabla^2 u_0\|_{\Ld^2_x}^{\frac dq}
\|\nabla\rho_0\|_{\Ld^2_x}^{1-\frac dp}
\|\nabla^2\rho_0\|_{\Ld^2_x}^{\frac dp}.
\end{eqnarray*}
Inserting the a priori estimates~\eqref{eq:u0-L2-energy}, \eqref{eq:nabrho0-L2-energy}, and~\eqref{eq:bnd-Deltau0-l21t-pre}, and combining with the above, we are led to
\begin{equation*}
\tfrac12\tfrac{\mathrm{d}}{\mathrm{d}t} \|\nabla^2\rho_0\|_{\Ld^2_x}^2+\|\nabla^3\rho_0\|_{\Ld^2_x}^2
\,\le\,C(t,h)
\|\nabla\rho^\circ\|_{\Ld^2_x}^{2(1-\frac dp)}
\Big(\|\nabla^2\rho_0\|_{\Ld^2_x}^{\frac{2d}p}+\|\nabla^2\rho_0\|_{\Ld^2_x}^{\frac{2d}p+\frac{4d}{q(4-d)}}\Big).
\end{equation*}
Choosing any $q<\infty$ if $d=2$, and choosing $q=\frac{2d}{d-2}$ if $d=3$, this proves
\begin{equation*}
\tfrac12\tfrac{\mathrm{d}}{\mathrm{d}t} \|\nabla^2\rho_0\|_{\Ld^2_x}^2+\|\nabla^3\rho_0\|_{\Ld^2_x}^2
\,\le\,C(t,h,\rho^\circ)\times\left\{\begin{array}{lll}
\|\nabla^2\rho_0\|_{\Ld^2_x}^{2}+1&:&\text{if $d=2$,}\\
\|\nabla^2\rho_0\|_{\Ld^2_x}^2+\|\nabla^2\rho_0\|_{\Ld^2_x}^4&:&\text{if $d=3$.}
\end{array}\right.
\end{equation*}
By Gr\"onwall's inequality, this yields
\begin{multline*}
\|\nabla^2\rho_0\|_{\Ld^\infty_t\Ld^2_x}^2+\|\nabla^3\rho_0\|_{\Ld^2_t\Ld^2_x}^2\\
\,\le\,C(t,h)\|\nabla^2\rho^\circ\|_{\Ld^2_x}^{2}\times\left\{\begin{array}{lll}
1&:&\text{if $d=2$,}\\
\exp\big(C(t,h)\|\nabla^2\rho_0\|_{\Ld^2_t\Ld^2_x}^2\big)&:&\text{if $d=3$.}
\end{array}\right.
\end{multline*}
By the a priori estimate~\eqref{eq:nabrho0-L2-energy}, this proves the claimed estimate on $\rho_0$. Combined with~\eqref{eq:bnd-Deltau0-l21t-pre}, this concludes the proof of~\eqref{eq:nabrho0-H2Linfty-est}.

We turn to the proof of~\eqref{eq:nabu0-H3Linfty-est}.
Applying $\nabla^2$ to both sides of the equation for $u_0$, testing with $\nabla^2u_0$, integrating by parts, and using Young's inequality, we find
\begin{equation*}
\|\nabla^3 u_0\|_{\Ld^2_x}
\,\lesssim\,\|\nabla h\|_{\Ld^2_x}
+\|\nabla\rho_0\otimes \nabla^2u_0\|_{\Ld^2_x}
+\|\nabla^2\rho_0\otimes\nabla u_0\|_{\Ld^2_x},
\end{equation*}
and thus, using as above the Gagliardo--Nirenberg interpolation inequality and the a priori estimates~\eqref{eq:u0-L2-energy}, \eqref{eq:nabrho0-L2-energy}, and~\eqref{eq:nabrho0-H2Linfty-est}, we obtain
\begin{equation}\label{eq:nabu0-H3Linfty-est-prebetter}
\|\nabla^3 u_0\|_{\Ld^2_x}^2
\,\lesssim\,C(t,h,\rho^\circ)\Big(1+\|\nabla^3\rho_0\|_{\Ld^2_x}^{d-2}\Big).
\end{equation}
Combined with~\eqref{eq:nabrho0-H2Linfty-est}, this proves~\eqref{eq:nabu0-H3Linfty-est}.

\medskip
\substep{2.4} $H^s$ regularity for $(\nabla u_0,\rho_0)$: for all integers $s >2$, provided that the initial condition further satisfies $\rho^\circ \in H^s(\T^d)$, and provided that $h\in\Ld^\infty_\loc(\R^+;H^{s-1}(\T^d)^d)$, we show that $\rho_0 \in \Ld^\infty_\loc(\R^+;H^s(\T^d))\cap\Ld^2_\loc(\R^+;H^{s+1}(\T^d))$ and that $u_0 \in \Ld^\infty_\loc(\R^+;H^{s+1}(\T^d)^d)$.
More precisely, we shall prove for all integers $s >2$,
\begin{gather}
\|\rho_0\|_{\Ld^\infty_t H^s_x}+\|\rho_0\|_{\Ld^2_tH^{s+1}_x}\,\lesssim_s\,\|\rho^\circ\|_{H^s_x}\exp\big(C_s\|u_0\|_{\Ld^1_tH^s_x}\big),\label{eq:estim-rho0-Hs}\\
\|u_0\|_{\Ld^\infty_tH^{s+1}_x}\,\lesssim_s\,\|h\|_{\Ld^\infty_tH^{s-1}_x}+\|u_0\|_{\Ld^\infty_tH^s_x}\|\rho_0\|_{\Ld^\infty_tH^s_x},\label{eq:estim-u0-Hs}
\end{gather}
and we note that these estimates indeed yield the conclusion by a direct iteration, starting from the results
of Substep~2.3 for $s=3$.

Let $s>2$ be a fixed integer, which entails in particular $s>\frac d2+1$ as we consider dimension~$d<4$.
We start with the proof of~\eqref{eq:estim-rho0-Hs}.
Applying $\langle\nabla\rangle^s:=(1-\Delta_x)^{s/2}$ to both sides of the equation for $\rho_0$, and testing it with $\langle\nabla\rangle^s\rho_0$, we find
\[\tfrac12\tfrac{\mathrm{d}}{\mathrm{d}t} \|\langle\nabla\rangle^s\rho_0\|_{\Ld^2_x}^2+  \tfrac{1}{\Pe}\|\nabla\langle\nabla\rangle^s\rho_0\|_{\Ld^2_x}^2=-\int_{\T^d}(\langle\nabla\rangle^s\rho_0)\langle\nabla\rangle^s(u_0\cdot\nabla\rho_0),\]
and thus, using the incompressibility constraint, 
\begin{eqnarray}
 \tfrac12\tfrac{\mathrm{d}}{\mathrm{d}t} \|\langle\nabla\rangle^s\rho_0\|_{\Ld^2_x}^2+ \tfrac{1}{\Pe} \|\nabla\langle\nabla\rangle^s\rho_0\|_{\Ld^2_x}^2
&=&-\int_{\T^d}(\langle\nabla\rangle^s\rho_0)[\langle\nabla\rangle^s,u_0\cdot\nabla]\rho_0\nonumber\\
&\le&\|\langle\nabla\rangle^s\rho_0\|_{\Ld^2_x}\|[\langle\nabla\rangle^s,u_0\cdot\nabla]\rho_0\|_{\Ld^2_x}.\label{eq:estim-predsdeozuhqbfndzq}
\end{eqnarray}
To estimate the last factor, we appeal to the following form of the Kato--Ponce commutator estimate~\cite[Lemma~X1]{KatoPonce} (see also~\cite[Lemma~3.4]{majda2002vorticity}): for all $u\in C^\infty(\T^d)^d$ and $\rho\in C^\infty(\T^d)$, we have
\begin{equation*}
\|[\langle\nabla\rangle^s,u\cdot\nabla]\rho\|_{\Ld^2_x}\,\lesssim_s\,\|u\|_{H^{s}_x}\|\rho\|_{W^{1,\infty}_x}+\|u\|_{W^{1,\infty}_x}\|\rho\|_{H^{s}_x},
\end{equation*}
and thus, by the Sobolev embedding with $s>\frac d2+1$,
\begin{equation}\label{eq:KatoPonce-est-00}
\|[\langle\nabla\rangle^s,u\cdot\nabla]\rho\|_{\Ld^2_x}\,\lesssim_s\,\|u\|_{H^{s}_x}\|\rho\|_{H^s_x}.
\end{equation}
Using this to estimate the last factor in~\eqref{eq:estim-predsdeozuhqbfndzq}, we are led to
\begin{equation}\label{eq:estim-rho0-Hs-pre}
\tfrac12\tfrac{\mathrm{d}}{\mathrm{d}t} \|\rho_0\|_{H^s_x}^2+\|\nabla\rho_0\|_{H^s_x}^2
\,\lesssim_s\,\|\rho_0\|_{H^s_x}^2\|u_0\|_{H^{s}_x},
\end{equation}
and the claim~\eqref{eq:estim-rho0-Hs} follows by Gr\"onwall's inequality.

We turn to the proof of~\eqref{eq:estim-u0-Hs}. Applying $\langle\nabla\rangle^s$ to both sides of the equation for $u_0$, and testing it with $\langle\nabla\rangle^su_0$, we find
\begin{eqnarray*}
\|\nabla\langle\nabla\rangle^su_0\|_{\Ld^2_x}^2
&=&\int_{\T^d}\langle\nabla\rangle^sh\cdot\langle\nabla\rangle^su_0-2c_0\int_{\T^d}\langle\nabla\rangle^s\D(u_0):\langle\nabla\rangle^s(\rho_0\D(u_0))\\
&\le&\int_{\T^d}\langle\nabla\rangle^sh\cdot\langle\nabla\rangle^su_0-2c_0\int_{\T^d}\langle\nabla\rangle^s\D(u_0):[\langle\nabla\rangle^s,\rho_0]\D(u_0),
\end{eqnarray*}
and thus, by Young's inequality,
\begin{equation*}
\|u_0\|_{H^{s+1}_x}
\,\lesssim\,\|h\|_{H^{s-1}_x}+\|[\langle\nabla\rangle^s,\rho_0]\D(u_0)\|_{\Ld^2_x}.
\end{equation*}
To estimate the last factor, we now appeal to the following form of the Kato--Ponce commutator estimate, instead of~\eqref{eq:KatoPonce-est-00},
with $s>\frac d2+1$,
\begin{equation*}
\|[\langle\nabla\rangle^s,\rho]\D(u)\|_{\Ld^2_x}\,\lesssim_s\,\|u\|_{H^{s}_x}\|\rho\|_{H^s_x}.
\end{equation*}
Using this to estimate the last factor in the above, the claim~\eqref{eq:estim-u0-Hs} follows.

\medskip
\substep{2.5} Time regularity.\\
For all $s\ge0$, if $\rho_0\in\Ld^\infty_\loc(\R^+;H^{s+2}(\T^d))$ and $u_0\in\Ld^\infty_\loc(\R^+;H^{s+3}(\T^d)^d)$, then the equation for $\rho_0$ yields, by the Sobolev embedding in dimension $d<4$,
\[\partial_t\rho_0\,=\,\tfrac1{\Pe}\Delta\rho_0-u_0\cdot\nabla\rho_0\,\in\,\Ld^\infty_\loc(\R^+;H^s(\T^d)).\]
Next, taking the time-derivative of both sides of the equation for $u_0$, we find
\[-\Div\big((1+c_0\rho_0)2\!\D(\partial_tu_0)\big)+\nabla\partial_tp_0=c_0\Div\big((\partial_t\rho_0)2\!\D(u_0)\big)+\partial_th.\]
For all $s\ge0$, if $\rho_0\in\Ld^\infty_\loc(\R^+;H^{s+2}(\T^d))$ and $u_0\in\Ld^\infty_\loc(\R^+;H^{s+3}(\T^d)^d)$, we deduce by elliptic regularity that $\partial_tu_0\in\Ld^\infty_\loc(\R^+;H^{s+1}(\T^d)^d)$ provided $h\in W^{1,\infty}_\loc(\R^+;H^{s-1}(\T^d)^d)$.
Higher time regularity is obtained similarly by direct induction.

\medskip
\step3 {Proof of item~(ii)}.\\
In the Stokes case $\Rey=0$, the system~\eqref{eq :(u_1,f_2)} reads
\begin{equation}
\left\{\begin{array}{l}\label{eq :(u_1,f_2)-rewr}
-\Delta u_1+\nabla p_1\,=\, \Div(\sigma_1[g_2])+\Div\big(\sigma_2[\rho_0,\nabla u_1]+\sigma_2[\rho_1+g_1,\nabla u_0]\big),\\[1mm]
(\partial_t-\tfrac1 \Pe \Delta_x+u_0\cdot\nabla_x)\rho_1+u_1\cdot\nabla_x\rho_0+\langle U_0n\cdot\nabla_xg_1\rangle=0,\\[1mm]
\Div(u_1)=0,
~~\int_{\T^d}u_1=0,\\[1mm]
\rho_1|_{t=0}=0,
\end{array}\right.
\end{equation}
where $g_1,g_2$ are given by~\eqref{f_1.explicit} and~\eqref{f_2.explicit}, respectively, in view of the explicit computations in Step~1.
Note that the form of $g_1,g_2$ ensures that this is a linear system for $(u_1,\rho_1)$.
We therefore focus on the proof of a priori energy estimates, while well-posedness and regularity properties easily follow.
On the one hand, testing the equation for $\rho_1$ with $\rho_1$ itself, and using the incompressibility constraint, we find
\begin{eqnarray*}
\tfrac12\tfrac{\dd}{\dd t}\|\rho_1\|_{\Ld^2_x}^2+\tfrac1{\Pe}\|\nabla\rho_1\|_{\Ld^2_x}^2
&=&-\int_{\T^d}\rho_1u_1\cdot\nabla\rho_0-\int_{\T^d}\Big(\int_{\Sp^{d-1}}\rho_1 U_0n\cdot\nabla_xg_1(\cdot,n)\,\dd n\Big)\\
&\le&\|\rho_1\|_{\Ld^2_x}\Big(\|u_1\|_{\Ld^2_x}\|\nabla\rho_0\|_{\Ld^\infty_x}+|U_0|\|\nabla_xg_1\|_{\Ld^2_x}\Big),
\end{eqnarray*}
and thus, by the explicit formula~\eqref{f_1.explicit} for $g_1$ and by the Sobolev embedding,
\begin{equation}\label{eq:estim-rho1-energy}
\tfrac12\tfrac{\dd}{\dd t}\|\rho_1\|_{\Ld^2_x}^2+\tfrac1{\Pe}\|\nabla\rho_1\|_{\Ld^2_x}^2
\,\lesssim\,\|\rho_1\|_{\Ld^2_x}\big(1+\|u_1\|_{\Ld^2_x}\big)\big(1+\|\rho_0\|_{H^3_x}+\|u_0\|_{H^2_x}\big)^2.
\end{equation}
On the other hand, testing the equation for $u_1$ with $u_1$ itself, we find
\[\int_{\T^d}|\nabla u_1|^2=-\int_{\T^d}\nabla u_1:\sigma_1[g_2]-\int_{\T^d}\nabla u_1:\sigma_2[\rho_0,\nabla u_1]-\int_{\T^d}\nabla u_1:\sigma_2[\rho_1+g_1,\nabla u_0],\]
and thus, recalling the definition of the elastic and viscous stresses, inserting the explicit formulas~\eqref{f_1.explicit} and~\eqref{f_2.explicit} for $g_1,g_2$, collecting all quadratic terms in $\nabla u_1$ in the left-hand side, and noting that integrals of monomials of odd degree on $\Sp^{d-1}$ vanish and using again the incompressibility constraint,
\begin{multline*}
\int_{\T^d}|\nabla u_1|^2
+\lambda(1+\tfrac12\theta)\int_{\T^d}\Big(\int_{\Sp^{d-1}}(n\otimes n:\nabla u_1)^2\,\dd n\Big)\rho_0\\
\hspace{-3cm}\,=\,
\tfrac1{2d}\lambda\theta\int_{\T^d}\Big(\int_{\Sp^{d-1}}(n\otimes n:\nabla u_1)(n \otimes n-\tfrac1d\Id)\,\dd n\Big)\\
\hspace{3cm}:\Big(\tfrac14A_2'(u_0,\rho_0)-\rho_0\D(u_0)^2
-d\rho_1\D(u_0)
-\tfrac{1}{d-1}U_0^2\nabla^2\rho_0\Big)\\
-\tfrac{1}8\lambda\theta\int_{\T^d}\bigg(\int_{\Sp^{d-1}}(n\otimes n:\nabla u_1)(n\otimes n:\nabla u_0)^2\,\dd n\bigg)\rho_0
-\int_{\T^d}\nabla u_1:\sigma_2[\rho_1+g_1,\nabla u_0].
\end{multline*}
Noting that the second left-hand side term is nonnegative as $\rho_0\ge0$, we deduce by the Sobolev embedding,
\begin{multline*}
\|\nabla u_1\|_{\Ld^2_x}
\,\lesssim\,
\|\rho_1\|_{\Ld^2_x}\|u_0\|_{H^3_x}
+\|\partial_t\rho_0\|_{\Ld^2_x}\|u_0\|_{H^3_x}\\
+\|\rho_0\|_{H^2_x}\|\partial_tu_0\|_{H^1_x}
+\|\rho_0\|_{H^2_x}\big(1+\|u_0\|_{H^3_x}\big)^2.
\end{multline*}
Combining this with~\eqref{eq:estim-rho1-energy} and with the Poincar\'e and Gr\"onwall inequalities, we deduce that $u_1$ is controlled in $\Ld^\infty_\loc(\R^+;H^1(\T^d)^d)$ and that $\rho_1$ is controlled in $\Ld^\infty_\loc(\R^+;\Ld^2(\T^d))\cap \Ld^2_\loc(\R^+;H^1(\T^d))$ provided that $u_0\in\Ld^\infty_\loc(\R^+;H^3(\T^d)^d)\cap W^{1,\infty}_\loc(\R^+;H^1(\T^d)^d)$ and $\rho_0\in\Ld^\infty_\loc(\R^+;H^2(\T^d))\cap W^{1,\infty}_\loc(\R^+;\Ld^2(\T^d))$.
\qed

\begingroup\allowdisplaybreaks
\subsection{Proof of Proposition~\ref{prop:main}}\label{sec:prop:main}
We aim to estimate the remainder terms $u_{k,\e}, f_{k,\e}, \rho_{k,\e},g_{k,\e}$ in the $\e$-expansion of the solution $(u_\e,f_\e)$ of the Doi--Saintillan--Shelley system, as defined through
\begin{eqnarray}
   u_\e & =& \sum_{j=0}^{k-1} \e^j u_j  +\e^k u_{k,\e},\label{new_eq:remainder-def}\\
   \rho_\e  & =& \sum_{j=0}^{k-1} \e^j \rho_j  +  \e^k \rho_{k,\e},\nonumber\\ 
   g_\e & =&  f_\e - \rho_\e  ~~=~~  \sum_{j=1}^{k-1} \e^j g_j  + \e^k g_{k,\e},\nonumber \\
   f_{k,\e} & =& \rho_{k,\e} + g_{k,\e},\nonumber
\end{eqnarray}
for $k \le3$.
We split the proof into five steps. We shall constantly use the short-hand notation $\mathcal{C}=\mathcal{C}(t)$ for multiplicative constants as defined in the statement,
the value of which may change from line to line.

\medskip
\step1 Energy estimate on the remainder $u_{3,\e}$.\\
Comparing equations for $u_\e,u_0,u_1$, cf.~\eqref{eq:eps-kin00} and Proposition~\ref{prop:first-terms}, using the linearity of $\sigma_1$ and the bilinearity of $\sigma_2$,
the remainder term $u_{3,\e}=\e^{-3}(u_\e-u_0-\e u_1 - \e^2 u_2)$ as defined in~\eqref{new_eq:remainder-def} satisfies the linearized Navier--Stokes equation
\begin{multline*}
\Rey(\partial_t+u_\e\cdot\nabla)u_{3,\e}
-\Delta u_{3,\e}-\Div(\sigma_2[f_\e,\nabla u_{3,\e}])+\nabla p_{3,\e}\\
=-\Rey\big[(u_{3,\e}\cdot\nabla) (u_0+\e u_1+\e^2u_2)+ (u_2\cdot\nabla) u_1+ ((u_1+\e u_2)\cdot\nabla) u_2\big]\\
+\Div(\sigma_1[g_{4,\e}])+\Div(\sigma_2[f_{3,\e} ,\nabla (u_0+\e u_1+\e^2u_2)])\\
+\Div(\sigma_2[ \rho_2 + g_2,\nabla u_1])+\Div(\sigma_2[ \rho_1 + g_1 +\e(\rho_2 + g_2),\nabla u_2]).
\end{multline*}
Here, we have also used the fact that $\sigma_1[\tau]=0$ for any function $\tau$ depending only on $x$, which follows from the simple observation that $\int_{\Sp^{d-1}}(n\otimes n-\tfrac1d\Id)\,\dd n=0$.
Testing this equation with $u_{3,\e}$, using the incompressibility constraints, and inserting the definition of elastic and viscous stresses~$\sigma_1,\sigma_2$, cf.~\eqref{eq:eps-kin00}, we get the following energy identity,
\begin{eqnarray*}
\lefteqn{\Rey\tfrac12\tfrac{\dd}{\dd t}\|u_{3,\e}\|_{\Ld^2_x}^2+\|\nabla u_{3,\e}\|_{\Ld^2_x}^2+\lambda\iint_{\T^d\times\Sp^{d-1}}\big(n\otimes n:\nabla u_{3,\e}\big)^2f_\e}\\
&=&-\Rey\int_{\T^d}u_{3,\e}\otimes u_{3,\e}:\nabla(u_0+\e u_1+\e^2u_2)\\
& & -\Rey\int_{\T^d}u_{3,\e}\cdot\big((u_2\cdot\nabla)u_1+((u_1 + \e u_2)\cdot\nabla)u_2\big)\\
& &-\lambda\theta\int_{\T^d}\nabla u_{3,\e}:\Big(\int_{\Sp^{d-1}}(n\otimes n-\tfrac1d\Id)\,g_{4,\e}\,\dd n\Big)\\
& &-\lambda \int_{\T^d}\nabla u_{3,\e}:\Big(\int_{\Sp^{d-1}}(n\otimes n)(\nabla (u_0+\e u_1+\e^2u_2))(n\otimes n) f_{3,\e} \,\dd n\Big)\\
& &-\lambda\int_{\T^d}\nabla u_{3,\e}:\Big(\int_{\Sp^{d-1}}(n\otimes n)(\nabla u_1)(n\otimes n)(\rho_2+g_2)\,\dd n\Big)\\
& &-\lambda\int_{\T^d}\nabla u_{3,\e}:\Big(\int_{\Sp^{d-1}}(n\otimes n)(\nabla u_2)(n\otimes n)(\rho_1+g_1+\e( \rho_2+g_2) )\,\dd n\Big),
\end{eqnarray*}
where we note that the viscous stress $\sigma_2[f_\e,\nabla u_{3,\e}]$ has led to an additional dissipation term in the left-hand side.
Appealing to the Cauchy--Schwarz inequality and to Young's inequality, we may then deduce
\begin{align*} 
&\Rey\tfrac{\dd}{\dd t}\|u_{3,\e}\|_{\Ld^2_x}^2+\|\nabla u_{3,\e}\|_{\Ld^2_x}^2\\
\quad & \lesssim \quad \Rey \left(1+\|(\nabla u_0,\nabla u_1,\nabla u_2)\|_{\Ld^\infty_x} \right)\Big(\|u_{3,\e}\|_{\Ld^2_x}^2+\|(u_1,u_2)\|_{\Ld^2_x}^2\Big) + \lambda^2\theta^2\|g_{4,\e}\|_{\Ld^2_{x,n}}^2 \\
&   \qquad +\lambda^2\|(\nabla u_0,\nabla u_1,\nabla u_2)\|_{\Ld^\infty_x}^2\|(f_{3,\e},g_1,g_2, \rho_1,\rho_2)\|_{\Ld^2_{x,n}}^2.
\end{align*}
Recalling the short-hand notation $\calC$ for multiplicative constants, and further decomposing $f_{3,\e}=\rho_{3,\e} +g_3 + \e g_{4,\e}$, we get 
\begin{multline}\label{est.u_3,eps}
\Rey\tfrac{\dd}{\dd t}\|u_{3,\e}\|_{\Ld^2_x}^2+\|\nabla u_{3,\e}\|_{\Ld^2_x}^2 \\
\, \le \, \mathcal{C}+\mathcal{C}\Big(\Rey\|u_{3,\e}\|_{\Ld^2_x}^2+\|\rho_{3,\e}\|_{\Ld^2_{x}}^2 \Big)
+(\lambda^2\theta^2+\calC\e^2)\|g_{4,\e}\|_{\Ld^2_{x,n}}^2.
\end{multline}

\medskip
\step2 Energy estimate for $\rho_{3,\e}$.\\
The equation satisfied by $\rho_{3,\e}=\e^{-3}(\rho_\e-\rho_0-\e\rho_1-\e^2\rho_2)$ takes the form
\begin{multline*}
    \partial_t \rho_{3,\e} - \tfrac1 \Pe  \Delta_x \rho_{3,\e}
    \,=\, - u_\e \cdot \nabla \rho_{3,\e}
    - u_{3,\e} \cdot \nabla (\rho_0 + \e \rho_1 + \e^2 \rho_2)\\
    - u_2 \cdot \nabla \rho_1
    - (u_1 + \e u_2) \cdot \nabla \rho_2 
    + \langle U_0 n \cdot \nabla_x g_{3,\e}\rangle.
\end{multline*}
Testing it with $\rho_{3,\e}$ itself, we can use the incompressibility constraint and several integrations by parts to get
\begin{multline*} 
    \tfrac12\tfrac{\mathrm{d}}{\mathrm{d}t} \|\rho_{3,\e}\|_{\Ld^2_x}^2 + \tfrac1 \Pe
    \|\nabla \rho_{3,\e}\|_{\Ld^2_x}^2  \\
    \, \lesssim \,  \Big(
    \|(\rho_0,\e\rho_1,\e^2\rho_2)\|_{\Ld^\infty_{x}} \|u_{3,\e}\|_{\Ld^2_x}  +  \|(\rho_1, \rho_2)\|_{\Ld^\infty_x}\|(u_1,u_2)\|_{\Ld^2_x} 
    +\|g_{3,\e}\|_{\Ld^2_{x,n}} \Big)\|\nabla \rho_{3,\e}\|_{\Ld^2_x},
\end{multline*}
hence, after absorption in the spatial dissipation,
\begin{equation*}
     \tfrac{\mathrm{d}}{\mathrm{d}t} \|\rho_{3,\e}\|_{\Ld^2_x}^2 + 
    \tfrac1{\Pe}\|\nabla \rho_{3,\e}\|_{\Ld^2_x}^2
    \, \lesssim \, 
     \big(\Pe\|\rho_0\|_{\Ld^\infty_{x}}^2+\e^2\calC\big)\|u_{3,\e}\|_{\Ld^2_x}^2
    +\mathcal{C} \big(1+\|g_{3,\e}\|_{\Ld^2_{x,n}}^2\big).
\end{equation*}
Using the maximum principle
\begin{equation}\label{eq:max-pr-rho0}
\|\rho_0\|_{\Ld^\infty_{t,x}}\,\le\,\|\rho^\circ\|_{\Ld^\infty_{x}},
\end{equation}
which holds for a (regular) solution $\rho_0$ of~\eqref{eq :(u_0,f_1)} with initial data~$\rho^\circ$, we deduce
\begin{equation}\label{est.rho_3,eps}
\tfrac{\mathrm{d}}{\mathrm{d}t} \|\rho_{3,\e}\|_{\Ld^2_x}^2
+\tfrac1{\Pe}\|\nabla \rho_{3,\e}\|_{\Ld^2_x}^2
\,\lesssim\,
\big(\Pe\|\rho^\circ\|_{\Ld^\infty_{x}}^2+\e^2\calC\big)\|u_{3,\e}\|_{\Ld^2_x}^2
+\mathcal{C} \big(1+\|g_{3,\e}\|_{\Ld^2_{x,n}}^2\big).
\end{equation}

\medskip
\step3 Energy estimate for $g_{4,\e}$.\\
The equation satisfied by $g_{4,\e}=\e^{-4}(g_\e-\e g_1-\e^2g_2-\e^3g_3)$
takes the form
\begin{equation}\label{new_eq:eqn-g3eps}
\e\partial_tg_{4,\e}-\Delta_n g_{4,\e}-\tfrac \e \Pe\Delta_xg_{4,\e}\\
\, = \, \sum_{i=1}^7 \mathrm{T}_i,
\end{equation}
in terms of
\begin{eqnarray*}
T_1&:=&-(\partial_t - \tfrac1 \Pe \Delta_x)g_3 -\Div_x\big( u_1 g_2 + u_2(g_1 + \eps g_2)\big),\\
T_2&:=&- \Div_x(u_\e g_{3,\e}),\\
T_3&:=&\e\Div_x(u_{3,\e} (g_1 + \e g_2)),\\
T_4&:=&\Div_x P_1^\perp \big( U_0 n f_{3,\e} \big),\\
T_5&:=&\Div_n\big(\pi_n^\bot(\nabla u_\e)n f_{3,\e}\big) ,\\
T_6&:=&\Div_n\big(\pi_n^\bot(\nabla u_{3,\e})n (\rho_0 + \e f_1 + \e^2 f_2 ) \big),\\
T_7&:=&\Div_n\big(\pi_n^\bot( \nabla u_1)n  f_2 + \pi_n^\bot( \nabla u_2)n (f_1 + \e f_2)\big).
\end{eqnarray*}
Testing this equation with $g_{4,\e}$ itself, we separately analyze the effect of the seven different source terms. For ${T}_1,{T}_6$, and ${T}_7$, we integrate by parts, we use Poincar\'e's inequality on~$\mathbb{S}^{d-1}$, recalling $\langle g_{4,\e} \rangle=0$, and we use the maximum principle~\eqref{eq:max-pr-rho0}, to the effect of
\begin{eqnarray*}
    \iint_{\T^d\times\Sp^{d-1}}   {T}_1 g_{4,\e}& \lesssim & \mathcal{C} \Vert \nabla_n g_{4,\e} \Vert_{\Ld^2_{x,n}}, \\
    \iint_{\T^d\times\Sp^{d-1}} {T}_6 g_{4,\e} & \lesssim &\big(\|\rho^\circ\|_{\Ld_{x}^\infty}+\e\calC\big)  \Vert \nabla u_{3,\e} \Vert_{\Ld^2_x} \Vert \nabla_n g_{4,\e} \Vert_{\Ld^2_{x,n}}, \\
    \iint_{\T^d\times\Sp^{d-1}} {T}_7 g_{4,\e} & \lesssim &\mathcal{C} \Vert \nabla_n g_{4,\e} \Vert_{\Ld^2_{x,n}}.
\end{eqnarray*}
For the term ${T}_2$, we use the incompressibility constraint, we further decompose $g_{3,\e}=g_3+\e g_{4,\e}$, and we integrate by parts, to the effect of
\begin{eqnarray*}
    \iint_{\T^d\times\Sp^{d-1}}{T}_2 g_{4,\e}
    &=& - \iint_{\T^d\times\Sp^{d-1}} g_{4,\e}\,\dv(u_\e g_{3,\e})\\
    & =&\iint_{\T^d\times\Sp^{d-1}} g_{4,\e}  u_\e\cdot \nabla_x g_{3} + \e \iint_{\T^d\times\Sp^{d-1}} g_{4,\e} u_\e\cdot \nabla_x  g_{4,\e} \\
    &=&\iint_{\T^d\times\Sp^{d-1}} g_{4,\e}u_\e\cdot \nabla_x g_{3},
\end{eqnarray*}
which is then estimated by
\begin{equation*}
\iint_{\T^d\times\Sp^{d-1}} {T}_2 g_{4,\e}
\, \le \, \mathcal{C}\Vert u_\eps \Vert_{\Ld^2_x}  \Vert \nabla_n g_{4,\e} \Vert_{\Ld^2_{x,n}},
\end{equation*}
again using Poincar\'e's inequality. For the term ${T}_3$, we get
\begin{eqnarray*}
    \iint_{\T^d\times\Sp^{d-1}} {T}_3 g_{4,\e} &=& \e \iint_{\T^d\times\Sp^{d-1}} g_{4,\e}u_{3,\e} \cdot \nabla_x (g_1 + \e g_2) \\
    &\le&  \e \mathcal C \Vert u_{3,\e} \Vert_{\Ld^2_x} \Vert \nabla_n g_{4,\e}\Vert_{\Ld^2_{x,n}}.
\end{eqnarray*}
For the term ${T}_4$, decomposing $f_{3,\e}=\rho_{3,\e} +g_3 + \e g_{4,\e}$, using that~$U_0$ is a constant and that~$P_1^\bot g_{4,\e}=g_{4,\e}$, and integrating by parts, we find
\begin{align*}
     \iint_{\T^d\times\Sp^{d-1}} {T}_4 g_{4,\e} & \, = \,
     \iint_{\T^d\times\Sp^{d-1}} g_{4,\e}\dv_x P_1^\perp \big( U_0 n f_{3,\e} \big)\\
     & \, = \,  \iint_{\T^d\times\Sp^{d-1}} g_{4,\e}\dv_x  \big( U_0 n (\rho_{3,\e} +g_3) \big)+\e U_0\iint_{\T^d\times\Sp^{d-1}} g_{4,\e}\dv_x (  n g_{4,\e})\\
     & \,=\, \iint_{\T^d\times\Sp^{d-1}} g_{4,\e}\dv_x  \big( U_0 n ( \rho_{3,\e} + g_3)\big),
\end{align*}
which is then estimated by
\begin{align*}
   \iint_{\T^d\times\Sp^{d-1}} {T}_4 g_{4,\e} \dd x \dd n  \,  \lesssim \,  \big( \|\nabla \rho_{3,\e}\|_{\Ld^2_x}+ \mathcal{C} \big) \| \nabla_n g_{4,\e} \|_{\Ld^2_{x,n}}.
\end{align*}
Finally, decomposing $f_{3,\e}=\rho_{3,\e}+g_3+\e g_{4,\e}$ and $u_\e=u_0+\e u_{1,\e}$, and integrating by parts, we split the remaining term ${T}_5$ as follows,
\begin{eqnarray*}
\lefteqn{\iint_{\T^d\times\Sp^{d-1}} {T}_5g_{4,\e}}\\
    &=&\iint_{\T^d\times\Sp^{d-1}} g_{4,\e} \Div_n\big(\pi_n^\bot(\nabla u_\eps )n f_{3,\e}\big)\\
    &=& \iint_{\T^d\times\Sp^{d-1}} g_{4,\e} \Div_n\big(\pi_n^\bot(\nabla u_\e)n (\rho_{3,\e}  +g_3 + \e g_{4,\e} )\big)
    \\
    &=& \iint_{\T^d\times\Sp^{d-1}} \rho_{3,\e} g_{4,\e} \Div_n\big(\pi_n^\bot(\nabla u_\e)n\big) +  \iint_{\T^d\times\Sp^{d-1}} g_{4,\e} \Div_n\big(\pi_n^\bot(\nabla u_\e)n g_3\big)\\
    &&  + \tfrac \eps 2  \iint_{\T^d\times\Sp^{d-1}} | g_{4,\e}|^2 \dv_n(\pi_n^\bot(\nabla u_\e)n)\\
    &=&\iint_{\T^d\times\Sp^{d-1}}  \rho_{3,\e} g_{4,\e} \Div_n\big(\pi_n^\bot(\nabla u_0)n\big) 
    +  \e \iint_{\T^d\times\Sp^{d-1}} \rho_{3,\e} g_{4,\e} \Div_n\big(\pi_n^\bot(\nabla u_{1,\e})n\big)\\
    &&+\iint_{\T^d\times\Sp^{d-1}} g_{4,\e} \Div_n\big(\pi_n^\bot(\nabla u_\e)n g_3\big)
    + \tfrac \eps 2  \iint_{\T^d\times\Sp^{d-1}} | g_{4,\e}|^2 \dv_n(\pi_n^\bot(\nabla u_0)n)\\
    &&+ \tfrac {\eps^2} 2  \iint_{\T^d\times\Sp^{d-1}} | g_{4,\e}|^2 \dv_n(\pi_n^\bot(\nabla u_{1,\e})n),
\end{eqnarray*}
which is then estimated by
\begin{multline*}
\iint_{\T^d\times\Sp^{d-1}}  {T}_5 g_{4,\e}
\, \lesssim \,   \e \mathcal{C}\Vert   g_{4,\e} \Vert_{\Ld^2_{x,n}}^2+ \mathcal{C} \big( \Vert \nabla u_{\e} \Vert_{\Ld^{2}_{x}}+\Vert \rho_{3,\e} \Vert_{\Ld^{2}_{x}}   \big)\Vert  \nabla_n g_{4,\e} \Vert_{\Ld^2_{x,n}}\\
+ \e  \|\nabla u_{1,\e}\|_{\Ld^2_x} \big(\|g_{4,\e}\|_{\Ld^4_x \Ld^2_n}^2 + \|  \rho_{3,\e} \|_{\Ld^{4}_{x}} \|  g_{4,\e} \|_{\Ld^4_{x} \Ld^2_n} \big).
\end{multline*}
Testing the equation~\eqref{new_eq:eqn-g3eps} with $g_{4,\e}$ itself and using the above estimates on the different source terms, together with Young's inequality, we end up with
\begin{multline}\label{eq:estim-g4eps-energy-0}
\e\tfrac{\dd}{\dd t}\|g_{4,\e}\|_{\Ld^2_{x,n}}^2+\|\nabla_ng_{4,\e}\|_{\Ld^2_{x,n}}^2+\tfrac \e \Pe \|\nabla_x g_{4,\e}\|_{\Ld^2_{x,n}}^2\\
\, \lesssim\,
\mathcal{C}\big(1+\|u_\e\|_{H^1_x}^2+ \| \rho_{3,\e}\|_{\Ld^2_x}^2 \big) 
+\|\nabla\rho_{3,\e}\|_{\Ld^2_x}^2
+\big(\|\rho^\circ\|_{\Ld^\infty_{x}}^2+ \e^2 \mathcal C\big)\|  u_{3,\e}\|_{H^1_x}^2  \\ 
+  \e \|\nabla u_{1,\e}\|_{\Ld^2_x} \big(\|\rho_{3,\e} \|_{\Ld^4_x}^2 + \|g_{4,\e}\|_{\Ld^4_x \Ld^2_n}^2 \big)
 + \e \mathcal{C}\Vert   g_{4,\e} \Vert_{\Ld^2_{x,n}}^2.
\end{multline}
In order to estimate the $\Ld^4_x$ norms in the right-hand side, we appeal to interpolation and to the Sobolev inequality on the torus in the following form: in dimension $d\le3$, we have for any $\delta>0$,
\begin{align*}
\|\nabla u_{1,\e}\|_{\Ld^2_x}\|g_{4,\e}\|_{\Ld^4_x \Ld^2_n}^2
\quad & \leq \quad C\|\nabla u_{1,\e}\|_{\Ld^2_x}\|g_{4,\e}\|_{\Ld^2_{x,n}}^{2(1-\frac d4)}\|g_{4,\e}\|_{H^1_x\Ld^2_n}^\frac d2\\
\quad & \leq \quad \delta \|\nabla_x g_{4,\e}\|_{\Ld^2_{x,n}}^2+C\delta^{-\frac{d}{4-d}}\|\nabla u_{1,\e}\|_{\Ld^2_x}^\frac4{4-d}\|g_{4,\e}\|_{\Ld^2_{x,n}}^2 \\
& \qquad + C\|\nabla u_{1,\e}\|_{\Ld^2_x}\| g_{4,\e}\|_{\Ld^2_{x,n}}^2.
\end{align*}
Likewise, further noting that $\rho_{3,\e}$ is mean-free (given that the defining equations for $\rho_0,\rho_1,\rho_2$ ensure $\int_{\T^d}\rho_1=\int_{\T^d}\rho_2=0$ and $\int_{\T^d}\rho_\e=\int_{\T^d}\rho_0=1$), we have for all $\eta>0$,
\begin{eqnarray*}
    \|\nabla u_{1,\e}\|_{\Ld^2_x} \|\rho_{3,\e} \|_{\Ld^4_{x}}^2 
    & \le& C \|\nabla u_{1,\e}\|_{\Ld^2_x}\|\rho_{3,\e}\|_{\Ld^2_{x}}^{2(1-\frac d4)}\| \nabla\rho_{3,\e} \|_{\Ld^2_{x}}^\frac d2 \\
    & \le& \eta\| \nabla \rho_{3,\e}\|_{\Ld^2_{x}}^2+C \eta^{-\frac{d}{4-d}}\|\nabla u_{1,\e}\|_{\Ld^2_x}^\frac4{4-d}\|\rho_{3,\e}\|_{\Ld^2_{x}}^2.
\end{eqnarray*}
Inserting these bounds into~\eqref{eq:estim-g4eps-energy-0}, and choosing $\delta \simeq \Pe^{-1}$ and $\eta \simeq 1$, we get for $\e \lesssim 1$
\begin{multline} \label{est.g_4,eps}
\e\tfrac{\dd}{\dd t}\|g_{4,\e}\|_{\Ld^2_{x,n}}^2+\|\nabla_ng_{4,\e}\|_{\Ld^2_{x,n}}^2 \\
\,\lesssim \, \mathcal{C}\big(1+\|u_\e\|_{H^1_x}^2 + \| \rho_{3,\e}\|_{\Ld^2_x}^2 \big)
+ \|\nabla\rho_{3,\e}\|_{\Ld^2_x}^2
+ \big(\|\rho^\circ\|_{\Ld^\infty_{x}}^2+ \e^2 \mathcal C\big)\|  u_{3,\e}\|_{H^1_x}^2 \\
+ \e\|\nabla u_{1,\e}\|_{\Ld^2_x}^\frac4{4-d}\|\rho_{3,\e}\|_{\Ld^2_{x}}^2 
+ \e \mathcal{C} \Big(1+\|\nabla u_{1,\e}\|_{\Ld^2_x}^\frac4{4-d}  \Big)\|g_{4,\e}\|_{\Ld^2_{x,n}}^2.
\end{multline}

\medskip
\step4 Conclusion.\\
Combining estimates~\eqref{est.rho_3,eps} and~\eqref{est.g_4,eps} in such a way that the term~$\|\nabla \rho_{3,\e}\|_{\Ld^2_x}^2 $ in the right-hand side of~\eqref{est.g_4,eps} can be absorbed into the left-hand side of~\eqref{est.rho_3,eps}, we obtain
\begin{multline*}
\Pe\tfrac{\dd}{\dd t}\|\rho_{3,\e}\|_{\Ld^2_x}^2
+\e\tfrac{\dd}{\dd t}\|g_{4,\e}\|_{\Ld^2_{x,n}}^2
+\|\nabla\rho_{3,\e}\|_{\Ld^2_x}^2
+\|\nabla_ng_{4,\e}\|_{\Ld^2_{x,n}}^2\\
\,\lesssim \, \mathcal{C}\big(1+\|u_\e\|_{H^1_x}^2 + \| \rho_{3,\e}\|_{\Ld^2_x}^2 \big)
+\big((1+\Pe)^2\|\rho^\circ\|_{\Ld^\infty_x}^2+\e^2\Cc\big)\|u_{3,\e}\|_{H^1_x}^2
+\Cc\big(1+\|g_{3,\e}\|_{\Ld^2_{x,n}}^2\big) \\
+ \e\|\nabla u_{1,\e}\|_{\Ld^2_x}^\frac4{4-d}\|\rho_{3,\e}\|_{\Ld^2_{x}}^2
+ \e \mathcal{C} \Big(1+\|\nabla u_{1,\e}\|_{\Ld^2_x}^\frac4{4-d}  \Big)\|g_{4,\e}\|_{\Ld^2_{x,n}}^2.
\end{multline*}
Further using the fact that
\begin{align*}
\|u_\e\|_{H^1_x}^2 \, \lesssim \,  \mathcal{C}+ \e^6\|u_{3,\e}\|_{H^1_x}^2,
\qquad\|g_{3,\e}\|_{\Ld^2_{x,n}}^2\| \,  \lesssim \, \mathcal{C}+ \eps^2 \|g_{4,\e}\|_{\Ld^2_{x,n}}^2,
\end{align*}
the above reduces to
\begin{multline}\label{eq:estim-combined-rho3g4}
\Pe\tfrac{\dd}{\dd t}\|\rho_{3,\e}\|_{\Ld^2_x}^2
+\e\tfrac{\dd}{\dd t}\|g_{4,\e}\|_{\Ld^2_{x,n}}^2
+\|\nabla\rho_{3,\e}\|_{\Ld^2_x}^2
+\|\nabla_ng_{4,\e}\|_{\Ld^2_{x,n}}^2\\
\,\lesssim \, \mathcal{C}\big(1+ \| \rho_{3,\e}\|_{\Ld^2_x}^2 \big)
+\big((1+\Pe)^2\|\rho_0\|_{\Ld^\infty_x}^2+\e^2\Cc\big)\|u_{3,\e}\|_{H^1_x}^2\\
+ \e \mathcal{C} \Big(1+\|\nabla u_{1,\e}\|_{\Ld^2_x}^\frac4{4-d}  \Big)\Big(\|\rho_{3,\e}\|_{\Ld^2_{x}}^2+\|g_{4,\e}\|_{\Ld^2_{x,n}}^2\Big).
\end{multline}
We now aim to combine this with~\eqref{est.u_3,eps} to conclude the proof, and we separately consider the Stokes and Navier--Stokes cases, splitting the proof into two further substeps.

\medskip
\substep{4.1} Stokes case $\Rey=0$, $d\le3$.\\
In the Stokes case, the estimate~\eqref{est.u_3,eps} becomes, using Poincar\'e's inequality,
\begin{eqnarray}
\Vert u_{3,\e}\Vert_{H^1_x}^2
&\lesssim& \mathcal{C}+ \lambda^2\theta^2 \|g_{4,\e}\|_{\Ld^2_{x,n}}^2 +\mathcal{C}\big( \Vert \rho_{3,\e} \Vert_{\Ld^2_x}^2+\e^2\|g_{4,\e}\|_{\Ld^2_{x,n}}^2 \big)\label{est.u_3,eps-rep}\\
& \lesssim & \mathcal{C}+ \lambda^2\theta^2 \| \nabla_n g_{4,\e}\|_{\Ld^2_{x,n}}^2 +\mathcal{C}\big( \Vert \rho_{3,\e}\Vert_{\Ld^2_x}^2+\e^2\|g_{4,\e}\|_{\Ld^2_{x,n}}^2 \big).\nonumber
\end{eqnarray}
We now insert this estimate for $u_{3,\e}$ in the right-hand side of~\eqref{eq:estim-combined-rho3g4}.
Provided that
\begin{align*}
\lambda^2\theta^2 (1+\Pe)^2\big(\|\rho^\circ\|_{\Ld^\infty_{x}}^2+ \e^2 \Cc\big) \, \ll \,  1
\end{align*}
is small enough, we can absorb the term involving $\| \nabla_n g_{4,\e}\|_{\Ld^2_{x,n}}^2$ in the estimate for $u_{3,\e}$, and we end up with 
\begin{multline*}
\Pe \tfrac{\mathrm{d}}{\mathrm{d}t} \|\rho_{3,\e}\|_{\Ld^2_x}^2  + \e\tfrac{\dd}{\dd t}\|g_{4,\e}\|_{\Ld^2_{x,n}}^2+ \|\nabla\rho_{3,\e}\|_{\Ld^2_x}^2+\|\nabla_ng_{4,\e}\|_{\Ld^2_{x,n}}^2 \\ 
\,\lesssim\,\Cc
+\Big(\Cc+\e\|\nabla u_{1,\e}\|_{\Ld^2_x}^\frac4{4-d} \Big)\|\rho_{3,\e}\|_{\Ld^2_{x}}^2
+\e\Cc\Big(1+\|\nabla u_{1,\e}\|_{\Ld^2_x}^\frac4{4-d} \Big) \|g_{4,\e}\|_{\Ld^2_{x,n}}^2.
\end{multline*}
By Gr\"onwall's inequality, we deduce 
\begin{multline*}
\|\rho_{3,\e}\|_{\Ld^\infty_t\Ld^2_x}^2 + \|\rho_{3,\e}\|_{\Ld^2_t H^1_x}^2 + \e\|g_{4,\e}\|_{\Ld^\infty_t \Ld^2_{x,n}}^2+\|\nabla_ng_{4,\e}\|_{\Ld^2_{t,x,n}}^2\\
\,\le\,\Cc(t)\Big(1+\|\rho_{3,\e}|_{t=0}\|_{\Ld^2_x}^2+\e\|g_{4,\e}|_{t=0}\|_{\Ld^2_{x,n}}^2\Big) \exp\bigg(\Cc(t)\int_0^t\|\nabla u_{1,\e}\|_{\Ld^2_{x}}^\frac4{4-d}\bigg).
\end{multline*}
Now expanding $u_{1,\e} =  u_1+\e u_2+\e^2 u_{3,\e}$
and using once more~\eqref{est.u_3,eps-rep}, we get
\begin{multline*}
\|\rho_{3,\e}\|_{\Ld^\infty_t\Ld^2_x}^2 + \|\rho_{3,\e}\|_{\Ld^2_t H^1_x}^2 + \e\|g_{4,\e}\|_{\Ld^\infty_t \Ld^2_{x,n}}^2+\|\nabla_ng_{4,\e}\|_{\Ld^2_{t,x,n}}^2\\
\,\le\,\Cc(t)\Big(1+\|\rho_{3,\e}|_{t=0}\|_{\Ld^2_x}^2+\e\|g_{4,\e}|_{t=0}\|_{\Ld^2_{x,n}}^2\Big)
\exp\bigg(\e^4\Cc(t)\Big(\|\rho_{3,\e}\|_{\Ld^\infty_t\Ld^2_{x}}+\|g_{4,\e}\|_{\Ld^\infty_t\Ld^2_{x,n}}  \Big)^\frac4{4-d}\bigg).
\end{multline*}
The well-preparedness assumption~\eqref{eq:well-prepared^4} precisely ensures that the initial terms in the right-hand side are uniformly bounded.
Using again Poincar\'e's inequality on $\S^{d-1}$ and appealing to a standard continuity argument, the stated estimates on~$\rho_{3,\e}$ and~$g_{4,\e}$ follow. The stated estimate on~$u_{3,\e}$ is then deduced from~\eqref{est.u_3,eps-rep}.

\medskip
\substep{4.2} Navier--Stokes case $\Rey=0$, $d=2$.\\
In the Navier--Stokes case,
provided that
\begin{align*}
\lambda^2\theta^2 (1+\Pe)^2 \big(\|\rho^\circ\|_{\Ld^\infty_{x}}^2+ \e^2 \Cc\big) \, \ll \,  1
\end{align*}
is small enough,
combining~\eqref{est.u_3,eps} and~\eqref{eq:estim-combined-rho3g4} in such a way that the term~$\|\nabla u_{3,\e}\|_{\Ld^2_x}$ in the right-hand side of~\eqref{eq:estim-combined-rho3g4} can be absorbed into the corresponding dissipation term in~\eqref{est.u_3,eps}, while further absorbing the term $\|g_{4,\e}\|_{\Ld^2_{x,n}}\lesssim\|\nabla_ng_{4,\e}\|_{\Ld^2_{x,n}}$ in the right-hand side of~\eqref{est.u_3,eps} into the corresponding dissipation term in~\eqref{eq:estim-combined-rho3g4},
we obtain
\begin{multline*}
\tfrac{\dd}{\dd t}\|u_{3,\e}\|_{\Ld^2_x}^2+\Pe\tfrac{\dd}{\dd t}\|\rho_{3,\e}\|_{\Ld^2_x}^2
+\e\tfrac{\dd}{\dd t}\|g_{4,\e}\|_{\Ld^2_{x,n}}^2
+\|\nabla u_{3,\e}\|_{\Ld^2_x}^2
+\|\nabla\rho_{3,\e}\|_{\Ld^2_x}^2
+\tfrac12\|\nabla_ng_{4,\e}\|_{\Ld^2_{x,n}}^2\\
\,\lesssim \, 
\Cc\Big(1+\|u_{3,\e}\|_{\Ld^2_x}^2+\|\rho_{3,\e}\|_{\Ld^2_{x}}^2\Big)
+\e \Cc\|\nabla u_{1,\e}\|_{\Ld^2_x}^\frac4{4-d}\Big(\|\rho_{3,\e}\|_{\Ld^2_{x}}^2+\|g_{4,\e}\|_{\Ld^2_{x,n}}^2\Big),
\end{multline*}
and thus, by Gr\"onwall's inequality,
\begin{multline*}
\|u_{3,\e}\|_{\Ld^\infty_t\Ld^2_x}^2
+\|\rho_{3,\e}\|_{\Ld^\infty_t\Ld^2_x}^2
+\e\|g_{4,\e}\|_{\Ld^\infty_t\Ld^2_{x,n}}^2
+\|\nabla u_{3,\e}\|_{\Ld^2_{t,x}}^2
+\|\nabla\rho_{3,\e}\|_{\Ld^2_{t,x}}^2
+\|\nabla_ng_{4,\e}\|_{\Ld^2_{t,x,n}}^2\\
\,\lesssim \, 
\Cc(t)\Big(1+\|u_{3,\e}|_{t=0}\|_{\Ld^2_x}^2
+\|\rho_{3,\e}|_{t=0}\|_{\Ld^2_x}^2
+\e\|g_{4,\e}|_{t=0}\|_{\Ld^2_{x,n}}^2\Big)
\exp\Big(\Cc(t)\int_0^t\|\nabla u_{1,\e}\|_{\Ld^2_x}^{\frac{4}{4-d}}\Big).
\end{multline*}
In the 2D case, the exponent in the exponential reduces to $\frac{4}{4-d}=2$. Then expanding again $u_{1,\e}=u_1+\e u_2+\e^2u_{3,\e}$, we can appeal to a continuity argument and the conclusion follows by noting that the well-preparedness assumption~\eqref{eq:well-prepared^4} precisely ensure the uniform boundedness of the initial terms.\qed

\subsection{Proof of Proposition~\ref{lem:Bous}}\label{seclem:Bous}
The proof mainly consists in inserting the explicit expressions for $g_1$ and $g_2$ computed in Proposition~\ref{prop:first-terms} in terms of $(u_0,\rho_0)$ and $(u_1,\rho_1)$ inside the systems~\eqref{eq :(u_0,f_1)} and \eqref{eq :(u_1,f_2)}.
Recalling the definition of stress tensors in~\eqref{eq:eps-kin00}, inserting the explicit expressions for $g_1,g_2$ in Proposition~\ref{prop:first-terms}, and using the elementary integral identities in~\eqref{eq:explicit-integral-n}, we are led to
\begin{align*}
       \sigma_1[g_1] \,= \,\lambda \theta \int_{\Sp^{d-1}} (n \otimes n - \tfrac{1}{d} \Id) \,g_1 \, \dd n \, = \, \lambda \theta \tfrac{\omega_d}{d(d+2)}  \rho_0 \D(u_0),
\end{align*}
and 
\begin{align*}
    \sigma_2[ \rho_0, u_0]  \,= \,\lambda \int_{\Sp^{d-1}} (n \otimes n) \nabla u_0 (n \otimes n) \rho_0 \, \dd n \, = \, \lambda \tfrac{2 \omega_d}{ d(d+2)} \rho_0 \D(u_0).
\end{align*}
This proves that the couple $(u_0,\rho_0)$ solves the system~\eqref{eq:v0rho0} with parameters $\eta_0 = 1$ and $\eta_1= \lambda \tfrac{\omega_d (\theta +2)}{2d(d+2)}$. Now turning to $(u_1,\rho_1)$, we compute
\begin{multline*}
\sigma_1[ g_2]\,=\, \lambda\theta \eps \int_{\Sp^{d-1}}\big(n\otimes n-\tfrac1d\Id\big)g_2\, \dd n\\
\,=\, \lambda\theta  \tfrac{\omega_d}{d(d+2)}\big(\rho_1\!\D( u_0)+\rho_0\!\D( u_1)\big)
+\lambda\theta \tfrac{\omega_d}{d(d+2)}\bigg(\!-\tfrac1{4d}A_2'(u_0,\rho_0)+2\tfrac{d+2}{d(d+4)}\rho_0\!\D(u_0)^2\\
 +\tfrac{1}{d(d-1)}U_0^2\nabla^2\rho_0-\tfrac1{d(d+4)}\Tr(\D(u_0)^2)\Id\bigg).
\end{multline*}
We also compute
\begin{align*}
\sigma_2[\rho_0,\nabla u_1]
\,=\,\lambda\int_{\Sp^{d-1}}(n\otimes n)(\nabla u_1)(n\otimes n)\,\rho_0 \, \dd n\,=\,\lambda\tfrac{2\omega_d}{d(d+2)}\rho_0\!\D(u_1),
\end{align*}
as well as 
\begin{align*}
\sigma_2[\rho_1+g_1,\nabla u_0]  \,=\,\lambda\tfrac{2\omega_d}{d(d+2)} \rho_1\!\D( u_0)
+\lambda\tfrac{2\omega_d}{d(d+2)(d+4)}\rho_0\Big(2\D( u_0)^2+\tfrac12\Tr(\D(u_0)^2
)\Id\Big).
\end{align*}
Inserting these identities into~\eqref{eq :(u_1,f_2)}, we verify that $u_1$ satisfies the fluid equation in \eqref{eq:v1rho1} for some modified pressure field~$p_1$ and with the expected coefficients
\begin{align*}
    \gamma_1 = -\lambda\theta\tfrac{\omega_d}{4d^2(d+2)} \qquad \text{and} \qquad  \gamma_2 = \lambda\tfrac{\omega_d}{2d^2(d+4)}\big(\theta+\tfrac{2d}{d+2}\big).
\end{align*}
We turn to the derivation of the corresponding equation for $\rho_1$. By the defining equations for $\rho_1$ in \eqref{eq :(u_1,f_2)}, we can write
\begin{equation*}
(\partial_t- \tfrac 1 \Pe \Delta+ u_0\cdot\nabla)\rho_1 = - u_1 \cdot \nabla \rho_0 -\langle U_0n\cdot\nabla_xg_1\rangle.\\
\end{equation*}
Inserting the explicit expressions for $g_1$ in Proposition~\ref{prop:first-terms}, and using the above computations, we find
\begin{equation*}
\langle U_0n\cdot\nabla_xg_1\rangle\,=\,-\tfrac1{d-1}U_0^2\langle n\otimes n:\nabla^2\rho_0\rangle\,=\,-\tfrac1{d(d-1)}U_0^2\Delta\rho_0,
\end{equation*}
and the conclusion follows for $\rho_1$.
\qed

\section{Modelling and physical aspects of the results}\label{sec:Doi-physics}
{In this section, we review the different models considered in this work and their physical content. We start with the Doi--Saintillan--Shelley model~\eqref{eq:eps-kin00} for suspensions and we describe in particular the signification of the different parameters. Next, we turn to ordered fluid models: we recall their standard definition and the non-Newtonian properties that they describe, before turning to non-standard modifications of those models at finite P\'eclet number and for inhomogeneous suspensions, leading to the specific equations~\eqref{eq:second-order-fluid-intro-eqn} that we have derived in the hydrodynamic approximation. Finally, we compare the explicit non-Newtonian properties of the latter qualitatively with experimental data, and we discuss how predictions get modified in the next-order (third-order fluid) description.}

In the present paper we deal with the rigorous derivation of macroscopic models from  the Doi--Saintillan--Shelley model.

\subsection{Doi and Doi--Saintillan--Shelley models}\label{sec:nondimensionalization}
{As outlined at the beginning of the introduction, the Doi--Saintillan--Shelley model is a {\it multiscale} kinetic model for suspensions, which can be viewed as intermediate between {\it microscopic} fluid-particle models and non-Newtonian {\it macroscopic} fluid models. In microscopic models, particles are modelled as rigid bodies in the fluid evolving according to individual ODEs or SDEs coupled to a fluid PDE. In contrast, in so-called {\it multiscale} models, particles are only described via their density, which satisfies some kinetic equation and affects the fluid flow in some simplified way: in a nutshell, it essentially amounts to neglecting multi-body effects.}
The usual dimensional form of the Doi--Saintillan--Shelley model reads (see \cite[Chapter 8]{DoiEdwards88} and~\cite{Saintillan2018review})
\begin{equation} \label{Doi}
\left\{\begin{array}{l}
\rho_{\fl} (\partial_t+u\cdot\nabla)u-\mu_{\fl}\Delta u+\nabla p\,=\,h+\Div(\sigma_1[f])+\Div(\sigma_2[f,\nabla u]),\\[2mm]
\partial_t f+\Div_x\big((u+V_0n)f\big)+\Div_n\big(\pi_n^\bot(\nabla u)nf\big)\\
\hspace{5cm}=D_{\tr}\Div_x \big((\Id + n \otimes n) \nabla_x f\big) + D_{\ro}\Delta_nf,\\[2mm]
\Div(u)=0,
\end{array}\right.
\end{equation}
where elastic and viscous stresses are respectively given by 
\begin{eqnarray}
\sigma_1[f]&:=&\left(3 k_B \Theta+ \alpha \mu_{\fl} |V_0| \ell^2\right) \int_{\Sp^{d-1}}\big(n\otimes n-\tfrac1d\Id\big)f(\cdot,n)\, \dd n,\nonumber\\
\sigma_2[f,\nabla u]&:=&  \tfrac{1}{2}\zeta_{\ro} \int_{\Sp^{d-1}}(n\otimes n):\D(u)(n\otimes n)\,f(\cdot,n)\, \dd n,\label{sigma_2}
\end{eqnarray}
where $h$ is an internal force,
{and where the different parameters $\rho_{\fl},\mu_{\fl},D_{\tr},D_{\ro},\zeta_{\ro},V_0,\ell$ are all assumed to be constant.
In this kinetic description, $\rho_{\fl}$ stands for the fluid density, $\mu_{\fl}$ for the solvent viscosity, $D_{\tr}$ and $D_{\ro}$ for the translational and rotational diffusion coefficients, $\zeta_{\ro}$ for the rotational resistance coefficient, $V_0$ for the self-propulsion speed of the particles, and $\ell$ for the length of the particles.
The fluid flow is coupled to the kinetic equation for the particle density via the additional stresses~$\sigma_1$ (elastic stress) and~$\sigma_2$ (viscous stress), which make the fluid equations non-Newtonian.
We briefly describe the structure and physical origin of those contributions:}
\begin{enumerate}[---]
\item {The viscous stress $\sigma_2$ arises from the rigidity of suspended particles in the fluid flow: it is formally understood by homogenization of the solid phase, viewed as inclusions with infinite shear viscosity in the fluid.}
The above expression~\eqref{sigma_2} is an approximation for very elongated particles: for general axisymmetric particles, the viscous stress involves additional terms depending on the precise shape of the particles (see e.g.~\cite{HinchLeal72} for spheroids), but slender body theory indeed shows that it reduces to the above form in the limit of very elongated particles (see e.g.~\cite[Section 3.4]{KimKarilla13}).
\smallskip\item The elastic stress $\sigma_1$ contains a passive and an active contribution. The passive part, proportional to the the Boltzmann constant $k_B$ and to the absolute temperature $\Theta$, is created by the random torques that are responsible for the rotational Brownian motion of the particles. 
These torques indeed create stresses due to the rigidity and anisotropy of the particles. 
On the other hand, the active part arises directly from the swimming mechanism, which is encoded in the parameter $\alpha \in \R$: a so-called puller particle corresponds to $\alpha>0$, and a pusher particle corresponds to $\alpha<0$.
\end{enumerate}
For $V_0=0$, the model~\eqref{Doi} reduces to the classical Doi model for passive Brownian particles~\cite{DoiEdwards88}.
Finally, for very elongated particles in 3D, we also note that the translational and rotational diffusion and resistance coefficients are asymptotically given as follows (see~\cite[Section 5.15]{Dhont96}),
\begin{equation}\label{eq:param-elong3D}
\begin{array}{rllrrll}
    D_{\ro} &=& \tfrac{k_B \Theta}{\zeta_{\ro}} + D_{\act},&\qquad&  \zeta_{\ro} &=& \tfrac{\pi \mu_{\fl} \ell^3}{3 \log(\ell/a)} ,\\[2mm]
    D_{\tr} &=& \tfrac{k_B \Theta}{\zeta_{\tr}},  &&   \zeta_{\tr} &=& \tfrac{2 \pi \mu_{\fl} \ell}{\log(\ell/a)},
\end{array}
\end{equation}
where $a$ is the width of the particles and where $D_{\act}$ is some possible active contribution to the rotational diffusion (tumbling). {We briefly explain how our starting point~\eqref{eq:eps-kin00} in this work follows from the above usual form~\eqref{Doi} of the Doi--Saintillan--Shelley model up to some simplification and non-dimensionalization procedure, and we then briefly recall the history of such Doi-type models.

\subsubsection{Simplification and non-dimensionalization}\label{sec:DOI-nondim-lala}
Compared to~\eqref{Doi}, we make two minor simplifications in the model:}
\begin{enumerate}[---]
\item While the translational diffusion in~\eqref{Doi} is proportional to $\Id+n\otimes n$, hence is twice as strong in the direction of particle orientation as in the orthogonal directions, we choose to neglect this~$O(1)$ difference and assume that the diffusion is isotropic. This choice is for simplicity and does not change anything in the analysis of the model.
\smallskip\item It has been observed in the seminal work~\cite{BergBrown72} that, for the example of E.~coli bacteria and related microswimmers, the contribution of thermal rotation and active tumbling is of the same order, and we therefore set $D_{\act}=0$ for simplicity.
\end{enumerate}
After these two simplifications, we non-dimensionalize the model~\eqref{Doi} in terms of the typical speed~$u_0$ of the fluid, the typical macroscopic length scale~$L$, and the typical number~$N$ of particles in a cube of side length~$L$. We further rescale time according to the timescale of the fluid flow, that is, $T = L/u_0$. More precisely, we define
\begin{align*}
\hat u(t,x) := \tfrac1{u_0}u(T t, L x), && \hat f(t,x) := \tfrac{L^d}{N} f(T t, L x), && \hat h(t,x):=\tfrac{L^2}{\mu u_0}h(Tt,Lx).
\end{align*}
This leads to the following dimensionless model, dropping the hats for simplicity,
\begin{equation} \label{Doi.nondimensional}
\left\{\begin{array}{l}
\Rey(\partial_t+u\cdot\nabla)u-\Delta u+\nabla p\,=\,h+\Div(\sigma_1[f])+\Div(\sigma_2[f,\nabla u]),\\[2mm]
\partial_t f+\Div_x\big((u+U_0n)f\big)+\Div_n\big(\pi_n^\bot(\nabla u)nf\big)=\frac 1 {\operatorname{Pe}}\Delta_xf+\frac 1 {\operatorname{Wi}}\Delta_nf,\\[2mm]
\Div(u)=0,
\end{array}\right.
\end{equation}
where the dimensionless counterparts of the additional stresses $\sigma_1,\sigma_2$ now take the form
\begin{eqnarray*}
\sigma_1[f]&=&(6\tfrac{\lambda}{\operatorname{Wi}}+ \gamma) \int_{\Sp^{d-1}}\big(n\otimes n-\tfrac1d\Id\big)f(\cdot,n)\, \mathrm{d}n,\\
\sigma_2[f,\nabla u]&=& \lambda\int_{\Sp^{d-1}}(n\otimes n)(\nabla u)(n\otimes n)\,f (\cdot,n) \, \mathrm{d}n.
\end{eqnarray*}
Here, $\operatorname{Re},\operatorname{Pe},\operatorname{Wi}>0$ stand for the so-called Reynolds, P\'eclet, and Weissenberg numbers, $\lambda >0$ depends only on the shape and number density of the particles, and $\gamma\in\R$ accounts for the active contribution to the stress. 
More precisely, these parameters are given by
\begin{equation*}
\begin{array}{rllrrll}
\Rey &=& \tfrac{\rho_{\fl} u_0 L}{\mu_{\fl}},&\qquad&
\lambda&=& \tfrac{\zeta_{\ro} N}{2\mu_{\fl} L^d}, \\[2mm]
\Pe &=& \tfrac{u_0 L \zeta_{\tr}}{k_B \Theta},&&
\Wi&=&  \tfrac{u_0 \zeta_{\ro}}{k_B \Theta L}, \\[2mm]
\gamma &=& \alpha N  \tfrac{|U_0| \ell^2}{L^{d-1}},&&
U_0&=& \tfrac{V_0}{u_0},
\end{array}
\end{equation*}
and we briefly comment on their range and interpretation:
\begin{enumerate}[---]
\item The Weissenberg  number $\operatorname{Wi}$ is the ratio between convection and relaxation timescales. For the kinetic viscoelastic models under consideration, it coincides with the rotational P\'eclet number. {For elongated particles in 3D, as $\zeta_{\ro}$ is proportional to the cube of the particle length, cf.~\eqref{eq:param-elong3D}, the regime when $\Wi$ is of order $O(1)$ is very narrow, and lies under standard flow conditions at particle lengths of around $10$ microns.} In this work, we are interested in the derivation of hydrodynamic approximations in case of very small $\Wi\ll1$: this means that we have applications in mind where the particles have a length of a few microns and below, which is in particular the case for many types of bacteria. For notational convenience, we rename the Weissenberg number as
\[\e\,:=\,\Wi\ll1.\]
\item The (translational) P\'eclet number $\Pe$ is typically much larger than its rotational counterpart $\Wi$
since\footnote{The same holds with an additional logarithmic correction in 2D.}
\[\tfrac\Pe\Wi\, =\, \tfrac{\zeta_{\tr}L^2}{\zeta_{\ro}}\,=\,\tfrac{6 L^2}{\ell^2}\, \gg\, 1.\]
In fact, from the application-oriented perspective, it makes sense to consider $\operatorname{Pe} \gg 1$. However, due to well-posedness and stability issues, our analysis crucially relies on keeping $\operatorname{Pe}$ not too large, and we shall generally assume $\Pe\simeq1$. We could actually allow for $\Pe$ to slightly diverge, but more slowly than the inverse of the volume fraction of the particles, cf.~\eqref{eq:smallness-cond-RE}, {which is anyway incompatible with typical applications.}
\smallskip\item The shape parameter $\lambda$ is typically quite small as it is proportional to the volume fraction $NL^{-d}$ of the particles.
\smallskip\item The prefactor in the viscous stress $\sigma_1$ reads
\[6\tfrac\lambda{\operatorname{Wi}}+\gamma\,=\,\tfrac\lambda{\operatorname{Wi}}\Big(6+ 2\alpha \tfrac{|V_0| \ell^2 \mu_{\fl}}{k_B \Theta}\Big),\]
where the term $\frac{|V_0| \ell^2 \mu_{\fl}}{k_B \Theta}$ is of order $1-10$ for typical microswimmers like E.~coli bacteria~\cite{BergBrown72}. We shall set for abbreviation
\begin{align}\label{def:theta}
    \theta\,:=\,6+ 2\alpha \tfrac{|V_0| \ell^2 \mu_{\fl}}{k_B \Theta}.
\end{align}
Note in particular that for passive particles $V_0=0$ we find $\theta=6$.
\smallskip\item The self-propulsion speed $V_0$ of active particles is typically around $10$ micron per second.
This is so slow with respect to typical shear flows considered in experimental settings that the ratio $U_0=\frac{V_0}{u_0}$ is typically tiny.
However, extremely low shear rates leading to $u_0\sim V_0$ can possibly also be achieved, as for instance in the experiments reported by~\cite{lopez2015turning}. For that reason, we choose not to neglect $U_0$ in the equations and to keep track of its effects.
\end{enumerate}
{This non-dimensionalized model~\eqref{Doi.nondimensional} coincides with~\eqref{eq:eps-kin00} upon setting $\Wi = \e$ and $\theta$ as in~\eqref{def:theta}, and it constitutes our starting point in this work.

\subsubsection{Derivation of the Doi--Saintillan--Shelley model}
For completeness, we briefly review the history of the derivation of Doi-type models.}
{In the physics literature}, the systematic theoretical study of the effective rheology of suspensions has been initiated by Einstein~\cite{Ein06}, who found that passive non-Brownian spherical rigid particles effectively increase the fluid viscosity by~$\frac52 \phi \eta$, where $\phi$ is the volume fraction of the particles and $\eta$ is the viscosity of the solvent. Jeffery~\cite{Jeffery22} studied the analogous problem for ellipsoidal particles and found an increase of the viscosity depending locally on orientations of the particles. By slender-body theory, in the limit of very elongated particles, Jeffery's viscous stress exactly takes the form of~$\sigma_2$ in~\eqref{eq:eps-kin00}; see e.g.~\cite{Brenner74,KimKarilla13}.
Starting in the 1930s, there is a vast literature in physics on the rheology of suspensions of non-spherical rigid particles, see e.g.~\cite{Kuhn32, Eisenschitz33, Peterlin38, Burgers38}, but this early work was restricted to specific fluid flows like simple shear, and Brownian effects were neglected.
Brownian particles were first considered in~\cite{Shima40, KuhnKuhn45, RisemannKirkwood50, Saito51}, by means of different models and justifications, finally leading to the additional elastic stress $\sigma_1$ in~\eqref{eq:eps-kin00}.
These models were largely reviewed in~\cite{HinchLeal71, HinchLeal72, Brenner74}:
in particular, the multiscale kinetic model~\eqref{eq:eps-kin00} in the passive case ($U_0=0$, $\theta=6$) then entered textbooks such as~\cite{DoiEdwards88, Graham} {and became known in the mathematical literature as the Doi model.}
The extension to active suspensions has been proposed by Saintillan and Shelley~\cite{saintillan2008instabilities2}: by coarse-graining force dipoles exerted by the particles on the fluid, {they derived a further contribution of the elastic stress $\sigma_1$ due to particle activity;} see also~\cite{haines2008effective,haines2009three,saintillan2010dilute,potomkin2016effective,Degond-19}.

{All of the aforementioned derivations are of formal or phenomenological nature. On the mathematical side,} the derivation of the viscous stress $\sigma_2$ from microscopic models has received considerable attention in recent years. When the fluid is modeled by the Stokes equations and the particle distribution is given, the effective increase of the fluid viscosity is by now well understood~\cite{HainesMazzucato12,NiethammerSchubert19,Gerard-VaretHillairet19,HillairetWu19,DuerinckxGloria20,Gerard-Varet21, Gerard-VaretHoefer21,DuerinckxGloria21,Gerard-VaretMecherbet20,Duerinckx20,DuerinckxGloria23}; {see~\cite{DG-22rev} for a review.}
In a similar setting, for active suspensions, the elastic stress $\sigma_1$ has been derived in~\cite{Girodroux-Lavigne22,Bernou-Duerinckx-Gloria-22}.
Regarding the derivation of the passive part of the elastic stress~$\sigma_1$ for Brownian particles, we refer 
to~\cite{hofer2023derivation}, where the authors start from a simplified microscopic model where {the particle dynamics is given} by Brownian motion and not coupled to the fluid.
Yet, non-Newtonian effects originate from the retroaction of the fluid on the particles, that is, from coupling the particle dynamics to the fluid flow: beyond the derivation of $\sigma_1$ and $\sigma_2$ {in steady regimes,} the derivation of kinetic equations for the particle density is thus of key interest.
First steps in that direction
have been undertaken in~\cite{HoferSchubert21,hofer2024non, Duerinckx23},
{but we emphasize that a complete derivation is a very challenging open problem for which the existing rigorous results are still highly fragmentary.}

\subsection{Ordered fluid models} \label{sec:standard-2nd-order}
Various non-Newtonian fluid models have been considered in the literature, taking into account nonlinear and memory effects in different ways on the macroscopic scale.
{In this work, we focus on so-called ordered fluid models that are shown to naturally appear as hydrodynamic approximations of the Doi--Saintillan--Shelley theory.
In this section, we start by recalling the usual form of such models, we review their non-Newtonian properties, we motivate their non-standard modification at finite P\'eclet number and for inhomogeneous suspensions, and summarize our findings regarding well-posedness and stability issues. Finally, we also discuss the explicit non-Newtonian parameters~\eqref{eq:parameters-2nd} that we have derived in the hydrodynamic approximation, and we comment on the physical content of the next-order approximation.}

\subsubsection{Usual form of ordered fluid models}
Ordered fluid equations take the form 
\begin{equation} \label{ordered.fluid}
\left\{\begin{array}{l}
\Rey(\partial_t +u\cdot\nabla)u-\Div(\sigma)+\nabla p= h,\\[1mm]
\Div(u)=0,
\end{array}\right.
\end{equation}
{where the stress $\sigma$ is a function of the Rivlin--Ericksen tensors $\{A_k(u)\}_k$ associated with the fluid velocity $u$, which} are defined iteratively through
\begin{gather}
A_1(u) \,:=\, 2\D(u)\,=\,\nabla u + (\nabla u)^T,\nonumber\\ A_{n+1}(u) \,:=\, (\partial_t  + u \cdot \nabla) A_n(u) + (\nabla u)^T A_n(u) + A_n(u) \nabla u,\qquad n\ge1. \label{eq:RE-tensor}
\end{gather}
In other words, $A_{n+1}(u)$ is the so-called lower-convected derivative of~$A_n(u)$, which ensures frame indifference of the equations.\footnote{In the literature, {\it upper}-convected instead of lower-convected derivatives are sometimes used to define ordered fluids (see e.g.~\cite[Section 6]{BCAH87}). This is merely a choice of convention, as both lead to equivalent fluid equations (although the value and interpretation of parameters of course depends on the chosen convention). Indeed, the upper-convected derivative of $A_n(u)$ can be rewritten in terms of the lower-convected derivative as
\[(\partial_t+u\cdot\nabla)A_n(u)- A_n(u) (\nabla u)^T -(\nabla u) A_n(u)\,=\,A_{n+1}(u)-\big(A_1(u)A_n(u)+A_n(u)A_1(u)\big).\]\vspace{-0.3cm}}
While first-order fluids are of the form~\eqref{ordered.fluid} with constitutive law $\sigma=\eta_0A_1(u)$, thus coinciding with standard Newtonian fluids,
second-order fluids are of the form~\eqref{ordered.fluid} with
\begin{align} \label{sigma.first.order}
\sigma \,:=\, \eta_0 A_1(u) + \alpha_1 A_2(u) + \alpha_2 A_1(u)^2,
\end{align}
see e.g.~\cite[Eq.~(8.48)]{Boehme12}, and third-order fluids amount to
\begin{multline}\label{eq:3rd-order}
    \sigma \,:=\, \eta_0 A_1(u) + \alpha_1 A_2(u)+ \alpha_2 A_1(u)^2 \\
    + \beta_1 A_3(u) + \beta_2\big(A_1(u)A_2(u) + A_2(u)A_1(u)\big) + \beta_3A_1(u) \Tr\big(A_1(u)^2\big),
\end{multline}
for some coefficients $\eta_0, \alpha_1, \alpha_2,\beta_1,\beta_2,\beta_3 \in \R$.\footnote{In arbitrary space dimension $d$, a further term $\beta_4 A_1(u)^3$ should be included in general in the stress $\sigma$ of third-order fluids, but it is redundant in dimensions $d\le3$ as we then have $B^3=\frac 1 2 B\,\Tr (B^2)+\frac 1 3 \Tr (B^3)$ by the Cayley--Hamilton theorem for any symmetric trace-free matrix $B$.}
Since third- and higher-order fluid models are fairly complicated and involve a lot of terms, the most widely used ordered fluid model for viscoelastic fluids is the second-order model~\eqref{sigma.first.order}.
For more details on modelling aspects and applications, we refer to~\cite{Dunn-Rajagopal-95, Boehme12, BAH87, Joseph13, Phan-Thien-MaiDuy13}.

\subsubsection{Non-Newtonian properties of ordered fluids}\label{sec:non-newt-prop}
We recall the basic rheological properties of these ordered fluid models depending on the different parameter values; we focus here on the 3D setting.
The coefficient $\eta_0>0$ in~\eqref{sigma.first.order} and~\eqref{eq:3rd-order} is the zero-shear viscosity, as in the first-order (Newtonian) model, while other parameters account for various non-Newtonian behaviors as we briefly describe.
\begin{enumerate}[---]
\item {\it Shear-dependent viscosity:} In a simple shear flow $u_0(t,x)=\kappa x_2e_1$ with shear rate $\kappa>0$, the shear-dependent viscosity is defined as
\begin{equation} \label{definition shear-dependent viscosity}
    \eta(\kappa)\,:=\,\frac{\sigma_{12}}\kappa.
\end{equation}
Noting that
\begin{align*}
\qquad A_1(u_0) = \kappa \begin{pmatrix} 0 & 1 & 0 \\ 1 & 0 & 0 \\ 0 & 0 & 0 
\end{pmatrix},&& A_1(u_0)^2 = \kappa^2 \begin{pmatrix} 1 & 0 & 0 \\ 0 & 1 & 0 \\ 0 & 0 & 0 
\end{pmatrix},
&&A_2(u_0) = \kappa^2 \begin{pmatrix} 0 & 0 & 0 \\ 0 & 2 & 0 \\ 0 & 0 & 0 
\end{pmatrix},
\end{align*}
we compute for the second-order fluid that the shear-dependent viscosity simply coincides with the zero-shear viscosity, $\eta(\kappa)=\eta_0$. For the third-order model, in contrast, we find a nontrivial shear-dependent relation,
\[\eta(\kappa)=\eta_0+2\kappa^2(\beta_2+\beta_3).\]
Most real-life viscoelastic fluids, and in particular passive dilute suspensions, happen to be shear-thinning, meaning that the map $\kappa\mapsto\eta(\kappa)$ is decreasing: this holds for the third-order fluid model provided that the coefficients satisfy $\beta_2+\beta_3<0$.

\smallskip\item {\it Normal-stress differences:}
In a simple shear flow $u_0(t,x)=\kappa x_2 e_1$ with shear rate $\kappa \in \R$, non-Newtonian fluids typically display  non-zero normal stresses. {This is responsible for a number of phenomena, of which the rod-climbing effect is the best known; see e.g.~\cite[Chapter 17]{Joseph13} and~\cite[Chapter 6]{BAH87}.}
Normal-stress coefficients are defined as\footnote{Beware of different sign conventions for the normal-stress coefficients. We follow here the choice of~\cite[Chapter 2.2]{Boehme12}.}
\begin{align} \label{definition normal stress}
    \nu_{10} \,:=\, \frac {\sigma_{11} - \sigma_{22}}{\kappa^2}, &&  \nu_{20} \,:=\, \frac {\sigma_{22} - \sigma_{33}}{\kappa^2},
\end{align}
and thus, for second- and third-order models,
\begin{align}\label{normal.stress}
    \nu_{10} \,=\, - 2  \alpha_1, &&  \nu_{20} \,=\,  2  \alpha_1 + \alpha_2.
\end{align}
In other words, $\alpha_1,\alpha_2$ are related to normal-stress coefficients via
\begin{align} 
    \alpha_1 \,=\, - \tfrac 1 2 \nu_{10}, && \alpha_2 \,=\,  \nu_{10} + \nu_{20}.
\end{align}
For most real-life viscoelastic fluids, and in particular for polymer solutions, it is found experimentally that $\nu_{10} > 0$, $ \nu_{20} < 0$, and that $|\nu_{20}|$ is considerably smaller than $|\nu_{10}|$  (by a factor of around $10$, see e.g. \cite[Section 2.2]{Boehme12} and~\cite[Section 2.2]{Phan-Thien-MaiDuy13}),
which means in particular, in terms of second-order fluid coefficients,
\begin{align}\label{eq:cond-alpha-exp}
    \alpha_1 < 0, && \alpha_2 > 0.
\end{align}
\item {\it Elongational viscosity:}
The apparent viscosity of a non-Newtonian fluid may be completely different in an elongational flow. Given a uniaxial elongational flow $u_0(t,x)=\kappa(x_1e_1-\frac12(x_2e_2+x_3e_3))$ in the direction $e_1$, the elongational viscosity is defined as
\begin{align} \label{definitio elongational viscosity}
\eta_{E} \,:=\, \frac{\sigma_{11} - \frac12(\sigma_{22}+ \sigma_{33})}{\kappa}.
\end{align}
Noting that
\begin{align}
    A_1(u_0) = \kappa\begin{pmatrix} 2 & 0 & 0 \\ 0 & -1 & 0 \\ 0 & 0 & -1 
    \end{pmatrix}, && A_2(u_0) = A_1(u_0)^2 = \kappa^2\begin{pmatrix} 4 & 0 & 0 \\ 0 & 1 & 0 \\ 0 & 0 & 1 
    \end{pmatrix},
\end{align}
we compute for the second-order model,
\begin{align*}
\eta_{E} \,=\, 3\eta_0 + 3 \kappa (\alpha_1 + \alpha_2).
\end{align*}
For real-life viscoelastic fluids, it is typically observed that the elongational viscosity increases with the strain rate (so-called strain-thickening behavior), which holds provided that coefficients satisfy~$\alpha_1 + \alpha_2 > 0$.
\smallskip\item {\it Retardation phase shift in oscillatory flow:}
In a simple shear flow $u_0(t,x)=\kappa(t) x_2e_1$ with oscillatory shear rate $\kappa(t) = \sin t$, we compute
\[A_2(u_0) = \dot \kappa\begin{pmatrix} 0 & 1 & 0 \\ 1 & 0 & 0 \\ 0 & 0 & 0 
    \end{pmatrix} + \kappa^2 \begin{pmatrix} 0 & 0 & 0 \\ 0 & 2 & 0 \\ 0 & 0 & 0 
    \end{pmatrix},\]
which gives rise to a phase shift for $\sigma_{12}$ in the second-order model, in form of
\[\sigma_{12}\,=\,\eta_0\kappa+\alpha_1\dot\kappa\,=\,(\eta_0^2+\alpha_1^2)^{\frac12}\sin\big(t+\arctan(\tfrac{\alpha_1}{\eta_0})\big).\]
This constitutes another typical (time-dependent) non-Newtonian feature, in link with the dependence of the stress on the flow history.
\end{enumerate}

\begin{rem}[Connection to Oldroyd--B model]\label{rem:Oldroyd-lala}
The Oldroyd--B model is a special case of a simple fluid model, which is particularly popular as a formal exact closure of the kinetic Hookean dumbbell model.
It is characterized by the following constitutive equation for the stress tensor,
\begin{align}
    \sigma = 2\eta_s \D(u) + \eta_p \tau, && \tau + \Wi \left((\partial_t+u\cdot\nabla)\tau- \tau (\nabla u)^T -(\nabla u) \tau\right)\,=\,2\D (u).
\end{align}
By a formal expansion with respect to $\Wi\ll1$, this model is found to agree to order $\Wi^k$ with a $k$th-order fluid model with some specific choice of parameters. In particular, to order $\Wi^2$, we recover the second-order fluid model~\eqref{sigma.first.order} with $\alpha_2=-2\alpha_1$, which means in particular that the second normal-stress coefficient vanishes, $\nu_{20}=0$.
We emphasize that this is {\it not} the case for the second-order fluid model that we derive from the Doi--Saintillan--Shelley theory: the derived parameters in~\eqref{eq:parameters-2nd} yield in 3D
\[\tfrac{\alpha_2}{\alpha_1}\,=\,\tfrac{\gamma_2}{\gamma_1}\,=\,-\tfrac{10}7(1+\tfrac{6}{5\theta}),\]
which is indeed $\ne-2$ provided $\theta\ne3$. Recall that passive suspensions correspond to $\theta=6$. Hence, at small $\Wi$, the viscoelastic effects of suspensions of rigid Brownian particles are not well described by an Oldroyd--B model.
\end{rem}

\subsubsection{Modified ordered fluid models at $\Pe<\infty$}\label{sec:nonstandard-ordered-Pe}
{While an infinite P\'eclet number $\Pe=\infty$ is usually considered in applications, we have needed to include a nontrivial spatial diffusion in the Doi--Saintillan--Shelley model~\eqref{eq:eps-kin00} in this work for technical stability reasons. 
Now, at $\Pe<\infty$, the definition of Rivlin--Ericksen tensors~\eqref{eq:RE-tensor} then has to be adapted, leading to a natural modification of ordered fluid models for which we could surprisingly find no reference.
Note that this is reminiscent of the corresponding version of the Oldroyd--B model with stress diffusion that is often considered both for analytical and numerical studies~\cite{renardyThomases21,DebiecSuli23}.}
More precisely, Rivlin--Ericksen tensors are modified as follows,
\begin{gather}\label{eq:RE-tensor modified}
    A_1'(u)\,:=\,A_1(u)\,=\,2\D(u),\\[2mm]
    A_{n+1}'(u)\,:=\,(\partial_t-{\tfrac1 \Pe} \Delta+u\cdot\nabla)A_n'(u)+(\nabla u)^TA_n'(u)+A_n'(u)(\nabla u),\quad n\ge1,\nonumber
\end{gather}
hence in particular $A_2'(u):=A_2(u)-\frac1\Pe\Delta A_1(u)$.
The second-order fluid equations~\eqref{ordered.fluid}--\eqref{sigma.first.order} are then replaced by
\begin{equation} \label{ordered.fluid-Pe}
\left\{\begin{array}{l}
\Rey(\partial_t +u\cdot\nabla)u-\Div(\sigma)+\nabla p= h,\\[1mm]
\sigma=\eta_0A_1(u)+\alpha_1A_2'(u)+\alpha_2A_1(u)^2,\\[1mm]
\Div(u)=0.
\end{array}\right.
\end{equation}
For higher-order fluid models, in this diffusive setting $\Pe<\infty$, several additional tensors actually need to be further included in the stress beyond the above modified Rivlin--Ericksen tensors:
for the third-order model, instead of~\eqref{eq:3rd-order}, the stress rather needs to be chosen in general as
\begin{multline}\label{eq:3rdorder-peclet}
\sigma\,=\,\eta_0A_1(u)+\alpha_1A_2'(u)+\alpha_2A_1(u)^2
+\beta_1A_3'(u)+\beta_1'B_3'(u)\\
+\beta_2\big(A_1(u)A_2'(u)+A_2'(u)A_1(u)\big)+\beta_3 A_1(u)\Tr\big(A_1(u)^2\big),
\end{multline}
in terms of the following additional tensor,
\[B_3'(u)\,:=\,(\partial_t-\tfrac1{\Pe}\Delta+u\cdot\nabla)A_1(u)^2+\big((\nabla u)^TA_1(u)^2+A_1(u)^2(\nabla u)\big).\]
At $\Pe=\infty$, this additional tensor is redundant as it indeed reduces to $B_3'(u)=A_1(u)A_2(u)+A_2(u)A_1(u)-A_1(u)^3$, so we recover the usual form~\eqref{eq:3rd-order}.\footnote{Recall that, as in~\eqref{eq:3rd-order}, a further term $\beta_4A_1(u)^3$ should always be included in the stress in arbitrary space dimension $d$, but it reduces to $\frac12\beta_4A_1(u)\Tr(A_1(u)^2)$ in dimensions $d\le3$.}

\subsubsection{Models for inhomogeneous suspensions}\label{sec:ordered.nonstandard}
{Whether at $\Pe=\infty$ or $\Pe<\infty$,} the above formulations of ordered fluid models describe the behavior of {\it spatially homogeneous} suspensions, that is, the behavior of fluids with a {\it uniform density} of suspended particles.
This can be naturally generalized to an inhomogeneous setting to describe non-uniform particle suspensions.
{Surprisingly again,} we could not find any account of this generalization in the literature. We emphasize however that such an inhomogeneous setting
should also arise naturally from the small-$\Wi$ expansion of the inhomogeneous Oldroyd--B model, which is the formal exact closure of the kinetic Hookean dumbbell model.
Note that we are still considering an homogeneous solvent fluid.
At infinite P\'eclet number~$\Pe=\infty$, the adaptation is straightforward:
the fluid equations~\eqref{ordered.fluid} are simply coupled to a conservation equation for the particle density~$\rho:=\fint_{\Sp^{d-1}}f(\cdot,n)\,\dd n$,\footnote{The physical particle density is rather given by $x\mapsto\int_{\Sp^{d-1}}f(x,n)\,\dd n$, but for notational convenience we choose to normalize it by the area of~$\Sp^{d-1}$. In particular, we have $\rho\in\frac1{\omega_d}\Pc(\T^d)$.}
\[(\partial_t  + u \cdot \nabla)\rho= 0,\]
while the stress $\sigma$ is now a function both of $\rho$ and of the Rivlin--Ericksen tensors $\{A_k(u)\}_k$ where non-Newtonian corrections to the pure solvent viscosity~$\eta_{0}$ are taken proportional to the suspended particle density~$\rho$. More precisely, the inhomogeneous second-order fluid model takes the form
\begin{eqnarray} \label{ordered.fluid-inhom}
\left\{\begin{array}{l}
\Rey  (\partial_t +u\cdot\nabla)u-\Div(\sigma)+\nabla p= h,\\[1mm]
(\partial_t  + u \cdot \nabla)\rho= 0,\\[1mm]
\sigma \,=\, (\eta_{0}+\eta_1\rho) A_1(u) + \alpha_1 \rho A_2(u) + \alpha_2 \rho A_1(u)^2,\\[1mm]
\Div(u)=0,
\end{array}\right.
\end{eqnarray}
for some coefficients $\eta_1,\alpha_1,\alpha_2\in\R$. Inhomogeneous versions of higher-order fluid models are formulated similarly.

At finite P\'eclet number $\Pe<\infty$, on the other hand, the particle density $\rho$ is no longer simply transported by the fluid, but rather solves a transport-diffusion equation,
\[(\partial_t -\tfrac1{\Pe}\Delta+u\cdot\nabla) \rho= 0.\]
In this diffusive setting, the structure of ordered fluid models becomes slightly more complicated: due to diffusion, the transport-diffusion operator that appears in the Rivlin--Ericksen tensors~\eqref{eq:RE-tensor modified} at finite P\'eclet number does not commute with multiplication with the particle density $\rho$. The second-order fluid model then rather takes the form
\begin{eqnarray} \label{2ndordered.fluid-inhom-ex}
\left\{\begin{array}{l}
\Rey(\partial_t +u\cdot\nabla)u-\Div(\sigma)+\nabla p= h,\\[1mm]
(\partial_t-\tfrac1{\Pe}\Delta+ u \cdot \nabla) \rho= 0,\\[1mm]
\sigma\,=\,(\eta_0+\eta_1\rho)A_1(u)+\alpha_1A_2'(\rho,u)+\alpha_2\rho A_1(u)^2,\\[1mm]
\Div(u)=0,
\end{array}\right.
\end{eqnarray}
for some coefficients $\eta_1,\alpha_1,\alpha_2\in\R$, in terms of the modified inhomogeneous second-order Rivlin--Ericksen tensor
\begin{equation}\label{2ndordered.fluid-inhom-ex/def-A2'}
A_2'(\rho,u)\,:=\,(\partial_t-{\tfrac 1 \Pe}\Delta+ u\cdot\nabla)(\rho A_1( u))+\rho\big((\nabla u)^TA_1( u)+A_1( u)(\nabla u)\big).
\end{equation}
Indeed, due to the diffusion, the latter quantity does not reduce to the Rivlin--Ericksen tensor defined in~\eqref{eq:RE-tensor modified}: we have $A_2'(\rho,u)\ne\rho A_2'(u)$ in general along solutions --- in contrast with the case of infinite P\'eclet number.\footnote{At infinite P\'eclet number $\Pe=\infty$, as the particle density $\rho$ satisfies $(\partial_t+u\cdot\nabla)\rho=0$, we indeed obtain $(\partial_t+u\cdot\nabla)(\rho A_1(u))=\rho(\partial_t+u\cdot\nabla)A_1(u)$, so that the system~\eqref{2ndordered.fluid-inhom-ex} reduces to~\eqref{ordered.fluid-inhom}.}
{This non-standard modified second-order fluid system is exactly the one that we derive from the Doi--Saintillan--Shelley theory, cf.~\eqref{eq:second-order-fluid-intro-eqn}.
Similarly, corresponding modified models can be formulated to higher order:} the inhomogeneous third-order model amounts to~\eqref{2ndordered.fluid-inhom-ex} with stress
\begin{multline*}
\sigma\,=\,(\eta_0+\eta_1\rho)A_1(u)+\alpha_1A_2'(\rho,u)+\alpha_2A_1(u)^2
+\beta_1A_3'(\rho,u)+\beta_1'B_3'(\rho,u)\\
+\beta_2\big(A_1(u)A_2'(\rho,u)+A_2'(\rho,u)A_1(u)\big)+\beta_3\rho A_1(u)\Tr(A_1(u)^2),
\end{multline*}
in terms of the modified inhomogeneous third-order Rivlin--Ericksen tensor
\[A_3'(\rho,u)\,:=\,(\partial_t-\tfrac1{\Pe}\Delta+u\cdot\nabla)A_2'(\rho,u)+(\nabla u)^TA_2'(\rho,u)+A_2'(\rho,u)(\nabla u),\]
and in terms of the following additional quantity, which needs to be included similarly as in~\eqref{eq:3rdorder-peclet} at finite P\'eclet number,
\[B_3'(\rho,u)\,:=\,(\partial_t-\tfrac1{\Pe}\Delta+u\cdot\nabla)(\rho A_1(u)^2)+\rho\big((\nabla u)^TA_1(u)^2+A_1(u)^2(\nabla u)\big).\]

\subsubsection{Well-posedness of ordered fluid models}
We focus for shortness on the second-order fluid model.
There has actually been a fair amount of confusion on the relevant range of parameters $\alpha_1,\alpha_2$: the sign condition $\alpha_1<0$ in~\eqref{eq:cond-alpha-exp} is motivated by experimental normal-stress measurements, but it appears inconsistent with thermodynamics; see e.g.~\cite{Dunn-Rajagopal-95} for a detailed discussion from the physics perspective.
From the mathematical point of view, this inconsistency leads to ill-posedness issues. The matter was investigated by Galdi~\cite{Galdi-93}, who showed the following for the second-order fluid equations~\eqref{ordered.fluid}--\eqref{sigma.first.order} at infinite P\'eclet number $\Pe=\infty$:
\begin{enumerate}[---]
\item the local-in-time well-posedness holds whenever $1/\alpha_1>-\lambda_1$, where $\lambda_1$ stands for the Poincar\'e constant in~$\T^d$, which is quite consistent with the choice~\eqref{eq:cond-alpha-exp} (although the case of a negative $\alpha_1$ with a small absolute value is prohibited);
\smallskip\item the long-time well-posedness, as well as the stability of steady solutions, can only hold provided that $\alpha_1>0$.
\end{enumerate}
For the corresponding system~\eqref{ordered.fluid-Pe} with finite P\'eclet number $\Pe<\infty$, the situation is even worse: even local-in-time well-posedness actually fails whenever $\alpha_1<0$ because the equation then behaves like a backward heat equation.
Yet, {even though the kinetic Doi model itself is known to be thermodynamically consistent (see~\cite[Chapter 8]{DoiEdwards88}}), our analysis leads us to a second-order fluid with $\alpha_1<0$, 
and it is thus crucial to determine what meaning should be given to the model in that case.
In fact, in the small-$\Wi$ limit, $\e:=\Wi\ll1$ we find
\begin{equation}\label{eq:alpha-gamma}
(\alpha_1,\alpha_2)=(\e\gamma_1,\e\gamma_2),\qquad\text{for some $\gamma_1<0,~\gamma_2\in\R$.}
\end{equation}
In this perturbative setting $\e\ll1$,
although the equation is ill-posed for fixed $\e$,
there are several ways to rearrange the nonlinearity and define well-posed notions of approximate solutions that only satisfy the system~\eqref{ordered.fluid-Pe} up to a higher-order $O(\e^2)$ remainder. A hierarchical way has been formulated in Proposition \ref{prop:hier-inhom}. 
There are also non-hierarchical ways to perturbatively make sense of the ill-posed second-order fluid model~\eqref{ordered.fluid-Pe}, which may be more desirable in particular for stability issues. We address this aspect in Appendix~\ref{app.Bous}.

\subsubsection{Non-Newtonian properties of hydrodynamic approximation}\label{secrem:non-newtonian actif}
We comment on the rheological features of the second-order fluid system~\eqref{eq:second-order-fluid-intro-eqn} with explicit parameters~\eqref{eq:parameters-2nd}, {which we have derived from the Doi--Saintillan--Shelley theory.}
We briefly describe the resulting non-Newtonian properties and how they depend on microscopic features.
We focus here on the physically relevant 3D setting, and
{we point out that in the passive case the parameter values agree
with~\cite[Eqn~(41) and Table~2]{HinchLeal72} and \cite[Eqn~(7.4)]{Brenner74}.}

\begin{enumerate}[$\bullet$] 
\item \textit{Effective spatial diffusion:}\\The spatial diffusion $\frac1{\Pe}$ of the suspended particle density is enhanced by particle activity even at infinite P\'eclet number: it is replaced by $\frac1{\Pe}+\e \mu_0$ with $\mu_0=\frac16U_0^2$. This naturally follows from the coupling of particle swimming velocities with their rotational diffusion. This phenomenon of increased mixing has been observed in studies such as~\cite{saintillan2007orientational,Saintillan2018review}.

\smallskip \item \textit{Modified zero-shear viscosity:}\\The presence of suspended particles leads to a non-trivial contribution to the zero-shear viscosity: in the homogeneous setting $\bar\rho_\e\equiv\frac1{\omega_d}$, we obtain {a zero-shear viscosity
\begin{align*}
\tilde\eta_0\,:=\,\eta_0+\tfrac1{\omega_d}\eta_1=1+\tfrac{1}{30}\lambda(2+\theta).
\end{align*}
In particular, in case of passive particles, this zero-shear viscosity is always larger than the plain fluid viscosity, $\tilde\eta_0>1$, while particle activity can reverse this effect. For a precise description, first recall that $\theta=6+2\alpha \tfrac{|V_0|\ell\mu_{\fl}}{k_B \Theta}$, cf.~\eqref{def:theta}, where~$\alpha$ characterizes the swimming mechanism:
\begin{enumerate}[---]
\item for passive particles $\alpha=0$, the zero-shear viscosity is $\tilde\eta_0=1+\frac{4}{15}\lambda>1$;
\item for so-called puller particles $\alpha>0$, the zero-shear viscosity is even larger than for passive particles;
\item for so-called pusher particles $\alpha<0$, the zero-shear viscosity is smaller than for passive particles, and it can even be smaller than the plain fluid viscosity provided that the activity of the particles is strong enough: we find $\tilde\eta_0<1$ if $\alpha<-\frac{4k_B\Theta}{|V_0|\ell\mu_{\fl}}$. 
\end{enumerate}
This prediction is consistent with well-known experimental results, see e.g.~\cite{sokolov2009reduction,rafai2010effective, lopez2015turning}, and it has been largely confirmed
in the literature~\cite{haines2008effective,haines2009three,saintillan2010extensional,saintillan2010dilute,gyrya2011effective,alonso2016microfluidic}. In particular, for {E.\@ coli} bacteria (a typical pusher particle), we can assess the value of the parameter $\alpha$ using the experimental measurements performed in~\cite{drescher2011fluid}: this yields $\alpha<-\frac{4k_B\Theta}{|V_0|\ell\mu_{\fl}}$ and the experimental findings of~\cite{Clement-13} then confirm our prediction that the effective zero-shear viscosity is smaller than the plain fluid viscosity.}

\smallskip\item \textit{Normal-stress differences:}\\
In the homogeneous setting $\bar\rho_\e\equiv\frac1{\omega_d}$, we obtain the following values for first and second normal-stress coefficients, as defined in~\eqref{definition normal stress},
\begin{equation*}
\nu_{10} =\tfrac{\theta}{90}\e\lambda \qquad \text{and} \qquad \nu_{20} =\tfrac{3-\theta}{315}\e\lambda.
\end{equation*}
{For passive particles ($\theta=6$), we thus obtain $\nu_{10}>0$, $\nu_{20}<0$, and $\nu_{10}/|\nu_{20}| = 7$, which agrees with experiments as discussed in Section~\ref{sec:non-newt-prop}.} The amplitude of these normal-stress coefficients is even increased in case of puller particles.
In contrast, for pusher particles, normal-stress coefficients are reduced, and a very large activity could even result into opposite effects: we find $\nu_{10}<0$ if~$\theta<0$,
and $\nu_{20}>0$ if~$\theta<3$.
This behavior was also predicted in~\cite{saintillan2010dilute,potomkin2016effective}, but has yet to be experimentally verified.

\smallskip \item \textit{Elongational viscosity:}\\
In the homogeneous setting $\bar\rho_\e\equiv\frac1{\omega_d}$, in a uniaxial elongational flow in the direction~$e_1$, that is, $\bar u_\e=\kappa(x_1 e_1 - \tfrac{1}{2} (x_2 e_2 + x_3 e_3))$, {we obtain the following value for the elongational viscosity, as defined in~\eqref{definitio elongational viscosity},
\begin{equation*}
\eta_{E}=\big(3+\tfrac{1}{10}\lambda(2+\theta)\big)+\kappa\tfrac{1}{140}\e\lambda(4+\theta).
\end{equation*}
This shows that passive suspensions ($\theta=6$) lead to a strain-thickening behavior, which is even increased in case of puller particles. In contrast, for pusher particles, the strain-thickening behavior is reduced, and a very large activity could even result in the opposite effect: the system becomes strain-thinning if $\theta<-4$ (that is, $\alpha<-\tfrac{5k_B\Theta}{|V_0|\ell\mu_{\fl}}$).}
\qedhere
\end{enumerate}

\subsubsection{Next-order hydrodynamic approximation}\label{secrem:third-order}
{The small-$\Wi$ expansion performed in Theorem~\ref{th:main} is easily pursued to higher orders, but we stick to a formal discussion for shortness.
As we show in Appendix~\ref{app:3rd-order},} the next-order description of the suspended particle density involves additional nontrivial transport and anisotropic diffusion terms depending on the surrounding fluid flow: we find
\begin{equation}\label{eq:barrho/3rd}
(\partial_t+\bar u_\e\cdot\nabla)\bar\rho_\e-\Div\Big(\big({\tfrac 1 \Pe}+\e\mu_0+\e^2\mu_1\D(\bar u_\e)\big)\nabla\bar\rho_\e\Big)=\e^2\mu_2\,\Div(\bar\rho_\e\Delta\bar u_\e)+O(\e^3),
\end{equation}
with explicit coefficients
\begin{eqnarray*}
\mu_0&:=&\tfrac{U_0^2}{d(d-1)},\\
\mu_1&:=&\tfrac{(3d+1)U_0^2}{d(d-1)^2(d+2)},\\
\mu_2&:=&\tfrac{U_0^2}{2d(d-1)(d+2)}.
\end{eqnarray*}
{Those nontrivial terms, which differ from the inhomogeneous ordered models suggested in Section~\ref{sec:ordered.nonstandard}, do actually vanish for particles with vanishing swimming velocity ($\mu_1=\mu_2=0$ if $U_0=0$).
In any case, we note that homogeneous spatial densities are still stable for~\eqref{eq:barrho/3rd} to order~$O(\e^3)$,} and we shall henceforth restrict for simplicity to the homogeneous setting,
\[\bar\rho_\e=\tfrac1{\omega_d}+O(\e^3).\]
In addition, we shall focus on the case of infinite P\'eclet number and of vanishing particle swimming velocity,
\[\Pe=\infty,\qquad U_0=0,\]
as this choice strongly simplifies the macroscopic equations and as it seems anyhow to be the most relevant setting physically, cf.~Section~\ref{sec:nondimensionalization} (recall however that our rigorous results do not hold for $\Pe=\infty$).
In this setting,
{we formally derive} in Appendix~\ref{app:3rd-order} the following third-order fluid equations from the Doi--Saintillan--Shelley theory,
\begin{equation}\label{eq:baru/3rd}
\left\{\begin{array}{l}
\Rey(\partial_t+\bar u_\e\cdot\nabla)\bar u_\e-\Div(\bar\sigma_\e)
%-\e^2\kappa_4\nabla^2:(\bar\rho_\e\nabla^2\bar u_\e)
+\nabla\bar p_\e=h+O(\e^3),\\[1mm]
\Div(\bar u_\e)=0,
\end{array}\right.
\end{equation}
where the stress is given by
\begin{multline}\label{eq:baru/3rd-bis}
\bar\sigma_\e
=
(1+\eta_1)A_1(\bar u_\e)
+\e\gamma_1 A_2(\bar u_\e)
+\e\gamma_2 A_1(\bar u_\e)^2\\
\qquad \,  +\e^2\kappa_1 A_3(\bar u_\e)
+\e^2\kappa_2\big(A_1(\bar u_\e)A_2(\bar u_\e)
+A_2(\bar u_\e)A_1(\bar u_\e)\big)
+\e^2\kappa_3 A_1(\bar u_\e)\Tr(A_1(\bar u_\e)^2)
\end{multline}
with explicit coefficients
\begin{eqnarray*}
\eta_1&:=&\lambda\tfrac{1}{2d(d+2)}(\theta+2),\\
\gamma_1&:=&-\lambda\theta\tfrac{1}{4d^2(d+2)},\\
\gamma_2&:=&\lambda\tfrac{1}{2d^2(d+4)}(\theta+\tfrac{2d}{d+2}),\\
\kappa_1&:=&\lambda\theta\tfrac{1}{8d^3(d+2)},\\
\kappa_2&:=&-\lambda\tfrac{1}{8d^3(d+4)}(3\theta+\tfrac{2d}{d+2}),\\
\kappa_3&:=&\lambda\tfrac{1}{8d^3(d+2)^2(d+4)(d+6)}\Big(2d(3d^2+10d+6)+\theta(d+4)(3d^2+11d+12)\Big).
\end{eqnarray*}
These third-order fluid coefficients coincide in the passive case ($\theta=6$) with those computed by Brenner~\cite[Eq.~(7.4)]{Brenner74} (once the notation is properly identified).
Regarding the non-Newtonian phenomena discussed in Section~\ref{secrem:non-newtonian actif} above, the main observation is that this third-order fluid model further describes the expected shear-thinning behavior of suspensions. Indeed, the shear-dependent viscosity is given in 3D as follows, cf.~Section~\ref{sec:non-newt-prop},
\begin{align*}
\kappa~~\mapsto~~& 1+\eta_1+2\e^2(\kappa_2+\kappa_3)\kappa^2\\
\, \qquad \qquad =~&1+\lambda\tfrac{\theta+2}{30}
-\e^2\lambda\tfrac{19\theta-12}{18900}\kappa^2,
\end{align*}
which is decreasing in $\kappa$ if and only if $\theta>\frac{12}{19}$.
As expected, this shows that passive suspensions ($\theta=6$) lead to a shear-thinning behavior, which is even increased in case of puller particles. In contrast, for pusher particles, the shear-thinning behavior is reduced, and a very large activity could even result in the opposite effect: the system becomes shear-thickening if $\theta<\frac{12}{19}$ (that is, $\alpha<-\frac{51k_B\Theta}{19|V_0|\ell\mu_{\fl}}$).
This possible shear-thickening effect was indeed measured experimentally in~\cite{Clement-13,lopez2015turning} for suspensions of {E.~coli} bacteria (pusher-type microswimmers) with strong enough activity.
We also refer to~\cite{haines2009three,gyrya2011effective,potomkin2016effective} for analytical and numerical results showing the same effect.

{For active particles with $U_0\ne0$, we have already mentioned that the transport-diffusion equation for the particle density suggested in Section~\ref{sec:ordered.nonstandard} actually fails at order $O(\e^2)$ and is to be replaced by the much more complicated equation~\eqref{eq:barrho/3rd}.
In fact, we show in Appendix~\ref{app:3rd-order} that the effect of particle swimming is even worse: the resulting fluid equation for $\bar u_\e$ would differ from the third-order fluid model~\eqref{eq:baru/3rd}--\eqref{eq:baru/3rd-bis}} even at uniform particle density and at infinite P\'eclet number. In particular, an additional dispersive correction~$-\e^2\kappa_4\Delta^2\bar u_\e$ would appear in the fluid equations. We skip the details as the case $|U_0|\ll1$ seems to be the most physically relevant anyway.

{\section{Conclusion}\label{sec:Conclusion}

Starting from the Doi--Saintillan--Shelley system~\eqref{eq:eps-kin00}, which is a multiscale kinetic model for suspensions of (passive or active) Brownian rigid particles, we have rigorously derived in this work second-ordered fluid models in the regime of a small Weissenberg number $\Wi=\e \ll 1$. 
In the following concluding remarks, we highlight some important aspects of our results. 
\begin{enumerate}[$\bullet$]
\item Second-order fluid models are nonlinear macroscopic fluid equations. They have the benefit over multiscale models that their rheological properties are easier to ``read off'' from the parameters and therefore to compare to experiments. Moreover, since they are lower-dimensional models (the particle orientation is eliminated as a variable), they are more accessible to numerical simulations.
\smallskip\item The regime of a small Weissenberg number is generally characterized by a small relaxation timescale compared to the convection timescale. For the  Doi--Saintillan--Shelley model, the relaxation mechanism is orientational Brownian motion. Hence, the relaxation timescale is the inverse of the rotational diffusion constant.
As discussed in Section~\ref{sec:nondimensionalization}, the regime of a small Weissenberg number is therefore highly relevant for applications with very small particles.
\smallskip\item The regime of a small Weissenberg number corresponds to weak non-Newtonian effects, which are then naturally captured in form of an $\e$-expansion. Upon truncation, this formally leads to so-called ordered fluid models: the $k$-th-order fluid model amounts to truncating $O(\e^{k})$ effects.
In this work, we (mainly) focus on the first nontrivial term in the expansion, that is, second-order fluids.
\smallskip\item Second-order fluid models can exhibit the following well-known rheological properties: zero-shear viscosity different from the solvent viscosity, normal-stress differences, strain-thinning (or thickening), and retardation phase shift in oscillatory flows. Third-order fluids can further exhibit shear-thinning (or thickening). All those non-Newtonian properties are determined quantitatively by the parameters in ordered fluid models.
\smallskip\item Our rigorous derivation from Doi--Saintillan--Shelley theory comes together with explicit formulas for ordered fluid parameters, which
%findings for the parameters in the second (and third) order fluid models quantitatively
are in full accordance with previous formal computations for passive suspensions. Moreover, both for active and passive suspensions, our findings agree qualitatively with experimental rheological measurements. Also note that our results definitely rule out the validity of the Oldroyd--B model at order $O(\e)$ for suspensions, as it would wrongly predict a vanishing second normal-stress coefficient.
\smallskip\item Experimentally, most real-life viscoelastic fluids display
%(almost) all real-life non-Newtonian fluids show
a positive first normal-stress coefficient. However, this leads to second-order fluid models that are thermodynamically inconsistent --- and mathematically ill-posed. One might thus doubt whether those models have any reasonable meaning. Yet, in the regime of a small Weissenberg number, we propose a way to resolve this apparent issue: in a nutshell,
while second-order models neglect $O(\e^2)$ effects anyway, they can be rearranged up to $O(\e^2)$ errors into actual well-posed models.
%this does not affect the model in a physically relevant way.
Hence, while one might argue that there are no real-world fluids that are exactly described by second-order fluid models, the latter still constitute very effective and meaningful approximate models for $\eps \ll 1$.
\smallskip\item Our results are limited by the technical necessity to include spatial diffusion $\frac1\Pe>0$. Even though we could allow for a spatial diffusion that is rather weak at small particle volume fraction, cf.~\eqref{eq:smallness-cond-RE}, our requirements on $\Pe$ are outside of the range for typical applications (see discussion in Section~\ref{sec:DOI-nondim-lala}).
In addition, since our analysis requires strong solutions, we are restricted to the 2D Navier--Stokes case or the 3D Stokes case.
\end{enumerate}
}

\appendix
\section{Well-posedness of the Doi--Saintillan--Shelley system}\label{app:well-posedness}
This appendix is devoted to the proof of Proposition~\ref{prop:well-posedness} {on the well-posedness of the Cauchy problem~\eqref{eq:eps-kin00}--\eqref{eq:eps-kin-bis} for the Doi--Saintillan--Shelley system.} Let us assume that $h=0$ for simplicity, as well as $\mathrm{Pe}=\lambda=\theta=\e=1$, since those parameters play no role in the analysis.
We then omit the subscript $\e$ in the notation and we set $(u_\e,f_\e)\equiv (u,f)$.
We focus on the Stokes case $\Rey=0$ in dimension $d=3$, but standard adaptations of the proof allow to treat the 2D Navier--Stokes case without important additional difficulty.
For a given distribution function $f: \R^+ \times \T^3 \times \mathbb{S}^{2} \to \R^+$, we shall use the short-hand notation
\begin{align*}
\rho_f \, := \, \int_{\mathbb{S}^{2}} f(\cdot,n) \, \dd n.
\end{align*}
We split the proof into four main steps.

\medskip
\step1 A priori energy estimates: we show that a smooth solution $(u,f)$ of the system~\eqref{eq:eps-kin00} satisfies for all $t\ge0$,
\begin{eqnarray}
\hspace{-0.3cm}\|\nabla u\|_{\Ld^\infty_t\Ld^2_x}&\lesssim&e^{Ct}\| \rho_{f^\circ} \|_{\Ld^2_x},\label{ineq:energyU-Stokes}\\
\hspace{-0.3cm}\| \rho_f \|_{\Ld^\infty_t\Ld^2_x}+\| \nabla \rho_f  \|_{\Ld^2_{t}\Ld^2_x}&\lesssim&e^{Ct}\| \rho_{f^\circ} \|_{\Ld^2_x},\label{ineq:energyRhof}\\
\hspace{-0.3cm}\|f\|_{\Ld^\infty_t\Ld^2_{x,n}}+\|\nabla_xf\|_{\Ld^2_t\Ld^2_{x,n}}+\|\nabla_nf\|_{\Ld^2_t\Ld^2_{x,n}}&\lesssim&\exp\Big(te^{Ct}\big(1+\|\rho_{f^\circ}\|_{\Ld^2_x}^4\big)\Big)\|f^\circ\|_{\Ld^2_{x,n}}.\label{ineq:energyf}
\end{eqnarray}
On the one hand, testing the equation for the fluid velocity $u$ with $u$ itself, using the incompressibility constraint, and inserting the form of $\sigma_1,\sigma_2$, we find
\begin{equation*}
\| \nabla u \|_{\Ld^2_x}^2+ \int_{\T^3\times\Sp^{2}} (\nabla u : n \otimes n)^2 f\,=\,-\int_{\T^3\times\Sp^2} \nabla u: \big(n\otimes n - \tfrac13\Id\big)f,
\end{equation*}
and since the second left-hand side term is nonnegative,
\begin{equation}\label{eq:est-Doi-nabu-rhof}
\|\nabla u \|_{\Ld^2_x} \,\lesssim\, \| \rho_f \|_{\Ld^2_{x}}.
\end{equation}
On the other hand, testing the kinetic equation for the particle density $f$ with $f$ itself, using again the incompressibility constraint, and integrating by parts,
we find
\begin{eqnarray}
\tfrac{1}{2}\tfrac{\dd}{\dd t} \| f \|_{\Ld^2_{x,n}}^2+  \| \nabla_x f \|_{\Ld^2_{x,n}}^2+\| \nabla_n f \|_{\Ld^2_{x,n}}^2
& = & \int_{\T^3 \times \mathbb{S}^{2}} \nabla_n f \cdot \pi_n^\bot(\nabla u)nf \nonumber\\
& =& -\tfrac{1}{2}\int_{\T^3 \times \mathbb{S}^{2}}| f|^2 \Div_n\big( \pi_n^\bot(\nabla u)n\big) \nonumber\\
&\lesssim& \Vert \nabla u \Vert_{\Ld^2_x}\Vert f \Vert_{\Ld^4_x \Ld^2_n}^2.\label{eq:estim-energy-fdff}
\end{eqnarray}
By Ladyzhenskaya's inequality, we can estimate the last factor as
\begin{equation*}
    \Vert f \Vert_{\Ld^4_x \Ld^2_n}^2
    \,\lesssim\, \Vert f \Vert_{\Ld^2_{x,n}}^{\frac{1}{2}}\Vert f \Vert_{H^1_x \Ld^2_n}^{\frac{3}{2}}
    \,\lesssim\,  \Vert f \Vert_{\Ld^2_{x,n}}^2 + \Vert f \Vert_{\Ld^2_{x,n}}^{\frac{1}{2}}\Vert \nabla_x f \Vert_{\Ld^2_{x,n}}^{\frac{3}{2}}.
\end{equation*}
Inserting this into~\eqref{eq:estim-energy-fdff}, and appealing to Young's inequality to absorb the norm of $\nabla_xf$, we deduce
\begin{equation}\label{eq:est-Doi-fffff}
\tfrac{\dd}{\dd t} \Vert f \Vert_{\Ld^2_{x,n}}^2+  \| \nabla_x f \|_{\Ld^2_{x,n}}^2+\Vert \nabla_n f \Vert_{\Ld^2_{x,n}}^2 \, \lesssim \, \big(\| \nabla u\|_{\Ld^2_x} +\| \nabla u\|_{\Ld^2_x}^{4}\big)\| f \|^2_{\Ld^2_{x,n}}.
\end{equation}
Moreover, integrating the equation for $f$ with respect to the angular variable, we get
\begin{equation*}
\partial_t \rho_f + u \cdot \nabla \rho_f -\Delta \rho_f+  U_0\,\Div\Big( \int_{\mathbb{S}^{2}} n f\,\dd n\Big)\,=\,0,
\end{equation*}
and thus, testing this equation with $\rho_f$ itself and using the incompressibility constraint,
\begin{equation*}
\tfrac{1}{2}\tfrac{\dd}{\dd t} \Vert \rho_f \Vert_{\Ld^2_x}^2+ \Vert \nabla\rho_f  \Vert_{\Ld^2_x}^2 \, = \, U_0 \int_{\T^3\times\Sp^{2}} n f\cdot \nabla \rho_f \, \lesssim \, \int_{\T^3} \rho_f |\nabla \rho_f|,
\end{equation*}
hence
\begin{equation*}
\tfrac{\dd}{\dd t} \| \rho_f \|_{\Ld^2_x}^2+\| \nabla \rho_f  \|_{\Ld^2_x}^2 \, \lesssim \, \Vert \rho_f \Vert_{\Ld^2_x}^2.
\end{equation*}
By Gr\"onwall's inequality, this proves~\eqref{ineq:energyRhof}, and the claim~\eqref{ineq:energyU-Stokes} then follows after combination with~\eqref{eq:est-Doi-nabu-rhof}. Further combining with~\eqref{eq:est-Doi-fffff} and appealing again to Gr\"onwall's inequality, the claim~\eqref{ineq:energyf} also follows.

\medskip
\step2 Construction of approximate solutions.\\
In order to prove the existence of a weak solution for the system~\eqref{eq:eps-kin00},
we argue by means of a Galerkin approximation method.
More precisely, given $k \in \N$, we introduce the following orthogonal projection on $\Ld^2(\T^3)$,
\begin{align*}
P_k \,:\, \Ld^2(\T^3) \,\twoheadrightarrow\, F_k\,:=\, \big\{ u \in \Ld^2(\T^3) :\hat{u}(l)=0 \text{ for all $|l|>k$}\big\},
\end{align*}
where $\{\hat{u}(l)\}_{l\in\Z^3}$ stands for the Fourier coefficients of a periodic function $u\in\Ld^2(\T^3)$.
For all $u \in F_k$ and $s \geq 0$, we obviously have
\begin{align}\label{eq:Hs-bnd-truncPk}
\| u \|_{H^s_x} \,\le\, \langle k\rangle^s\| u \|_{\Ld^2_x}.
\end{align}
Given an initial condition $f^\circ\in H^1\cap\Pc(\T^3\cap\Sp^2)$, we shall consider the following approximate system,
\begin{equation}\label{system:regularized-plop}
\left\{\begin{array}{l}
-\Delta u_k+\nabla p_k
=P_k \Div(\sigma_1[f])+P_k\Div(\sigma_2[f_k,\nabla P_k u_k]),\\[1mm]
\partial_t f_k+\Div_x\big((u_k+U_0n)f_k\big)+\Div_n\big(\pi_n^\bot( \nabla P_ku_k)n f_k\big)=\Delta_xf_k+\Delta_nf_k, \\[1mm]
\Div(u_k)=0,
\quad\textstyle \int_{\T^3} u_k = 0,\\[1mm]
f_k|_{t=0}=f^\circ,
\end{array}\right.
\end{equation}
and we claim that this system admits a weak solution $(u_k,f_k)$ with
\begin{eqnarray}
u_k&\in& \Ld^\infty_\loc(\R^+;H^1(\T^3)^3),\label{system:regularized-plop-regul}\\
f_k&\in& \Ld^{\infty}_\loc(\R^+; \Ld^{2}\cap\Pc(\T^3\times\Sp^2)) \cap \Ld^{2}_\loc(\R^+; H^1(\T^3\times\Sp^2)).
\end{eqnarray}
To prove this, we argue by means of a Schauder fixed-point argument and we split the proof into four further substeps.

\medskip
\substep{2.1} Fixed-point problem.\\
Given~$T>0$, let
\[E\,:=\,\Big\{v\in \Ld^\infty(0,T; H^1(\T^3)^3):\Div(v)=0,~\textstyle\int_{\T^3}v=0\Big\},\]
and for all $v \in E$ define $\Lambda_k(v):=u$ as the unique solution of
\begin{equation}\label{eq:def-Lambda}
\left\{\begin{array}{l}
-\Delta u+\nabla p
=P_k \Div(\sigma_1[g])+P_k \Div(\sigma_2[g,\nabla P_k u]),\\[1mm]
\partial_t g+\Div_x\big((v+U_0n)g\big)+\Div_n\big(\pi_n^\bot( \nabla P_k v)n g\big)=\Delta_xg+\Delta_ng, \\[1mm]
\Div(u)=0,
\quad\textstyle \int_{\T^3}u= 0,\\[1mm]
g|_{t=0}=f^\circ.
\end{array}\right.
\end{equation}
Note that we use the velocity field~$u$ and not~$v$ in the viscous stress $\sigma_2$ in the equation for~$u$ so as to preserve its dissipative structure.
The existence of a weak solution $u_k$ for the approximate system~\eqref{system:regularized-plop} amounts to finding a fixed point $\Lambda_k(u_k)=u_k$.

Before going on with the fixed-point problem, we first check that the above system~\eqref{eq:def-Lambda} is indeed well-posed and defines a map $\Lambda_k:E\to E$.
First, given $v\in E$,
{by standard parabolic theory,}
the above kinetic equation for $g$ admits a unique weak solution
\[g\,\in\, \Ld^{\infty}(0,T; \Ld^{2}\cap\Pc(\T^3\times\Sp^2)) \cap \Ld^{2}(0,T; H^1(\T^3\times\Sp^2)).\]
Moreover, the a priori estimates of Step~1 can be repeated to the effect of
\begin{eqnarray}
\|\rho_g\|_{\Ld^\infty_T\Ld^2_x}+\|\nabla\rho_g\|_{\Ld^2_{T}\Ld^2_x}&\lesssim&e^{CT}\|\rho_{f^\circ}\|_{\Ld^2_x},\label{eq:apriori-est-ggg}\\
\|g\|_{\Ld^\infty_T\Ld^2_{x,n}}+\|\nabla_xg\|_{\Ld^2_{T}\Ld^2_{x,n}}+\|\nabla_ng\|_{\Ld^2_{T}\Ld^2_{x,n}}&\lesssim&\exp\Big(CT\big(1+\|\nabla v\|_{\Ld^\infty_T\Ld^2_x}^4\big)\Big)\|f^\circ\|_{\Ld^2_x}.\nonumber
\end{eqnarray}
To ensure the well-posedness of the elliptic equation for $u$ in~\eqref{eq:def-Lambda}, we first need to check that $\rho_g$ has some additional regularity.
Integrating the equation for $g$ with respect to the angular variable, we find the following equation for $\rho_g$,
\[\partial_t\rho_g+v\cdot\nabla\rho_g-\Delta\rho_g+U_0\,\Div\Big(\int_{\Sp^2}ng\,\dd n\Big)=0,\]
and thus, testing the equation with $\Delta\rho_g$, we find
\begin{eqnarray*}
\tfrac12\tfrac{\dd}{\dd t}\|\nabla\rho_g\|_{\Ld^2_x}^2+\|\Delta\rho_g\|_{\Ld^2_x}^2
&=&-\int_{\T^3}( v\cdot\nabla\rho_g)\Delta\rho_g-U_0\int_{\T^3}(\Delta\rho_g)\,\Div\Big(\int_{\Sp^2}ng\,\dd n\Big)\\
&\lesssim&\|v\|_{\Ld^4_x}\|\nabla\rho_g\|_{\Ld^4_x}\|\Delta\rho_g\|_{\Ld^2_x}+\|\Delta\rho_g\|_{\Ld^2_x}\|\nabla_xg\|_{\Ld^2_{x,n}}.
\end{eqnarray*}
By Ladyzhenskaya's and Poincar\'e's inequalities, the first right-hand side term can be bounded by
\begin{eqnarray*}
\|v\|_{\Ld^4_x}\|\nabla\rho_g\|_{\Ld^4_x}\|\Delta\rho_g\|_{\Ld^2_x}
&\lesssim&\|v\|_{\Ld^2_x}^\frac14\|\nabla v\|_{\Ld^2_x}^\frac34\|\nabla\rho_g\|_{\Ld^2_x}^\frac14\|\nabla^2\rho_g\|_{\Ld^2_x}^\frac74\\
&\lesssim&\|\nabla v\|_{\Ld^2_x}\|\nabla\rho_g\|_{\Ld^2_x}^\frac14\|\nabla^2\rho_g\|_{\Ld^2_x}^\frac74,
\end{eqnarray*}
and thus, inserting this into the above, recalling $\|\nabla^2\rho_g\|_{\Ld^2_x} \lesssim \|\Delta\rho_g\|_{\Ld^2_x}$, and using Young's inequality, we are led to
\begin{eqnarray*}
\tfrac{\dd}{\dd t}\|\nabla\rho_g\|_{\Ld^2_x}^2+\|\Delta\rho_g\|_{\Ld^2_x}^2
&\lesssim&\|\nabla v\|_{\Ld^2_x}^8\|\nabla\rho_g\|_{\Ld^2_x}^2+\|\nabla_xg\|_{\Ld^2_{x,n}}^2.
\end{eqnarray*}
By Gr\"onwall's inequality, this implies
\[\|\nabla\rho_g\|_{\Ld^\infty_T\Ld^2_x}
+\|\Delta\rho_g\|_{\Ld^2_T\Ld^2_x}
\,\lesssim\,\Big(\|\nabla\rho_{f^\circ}\|_{\Ld^2_x}+\|\nabla_xg\|_{\Ld^2_T\Ld^2_{x,n}}\Big)\exp\Big(C\int_0^t\|\nabla v\|_{\Ld^2_x}^8\Big).\]
Combined with~\eqref{eq:apriori-est-ggg}, this shows that $\rho_{g} \in \Ld^{\infty}(0,T; H^1(\T^3)) \cap \Ld^{2}(0,T; H^2(\T^3))$. In particular, we infer $\rho_{g}(t) \in \Ld^\infty(\T^3)$ for almost all $t \in [0,T]$.
With this additional regularity result for $\rho_g$, we can now finally ensure the well-posedness of the elliptic equation for $u$ in~\eqref{eq:def-Lambda}. Indeed, for all~$t\in[0,T]$, the weak formulation of this equation for $u(t)$ can be written as follows: for all~$w\in E$,
\begin{align*}
B_{t,k,g}(w,u(t))\,:=\,&\int_{\T^3} \nabla w:\nabla u(t)+\int_{\T^3\times\Sp^2}\big(n\otimes n:\nabla P_kw\big)\big(n\otimes n:\nabla P_ku(t)\big)\,g(t)\\
&\hspace{2cm}\,=\,-\int_{\T^3}\nabla w\cdot P_k\sigma_1[g(t)].
\end{align*}
For almost all $t$, as $\rho_{g}(t) \in \Ld^\infty(\T^3)$, we note that $B_{t,k,g}$ is a coercive continuous bilinear functional. As in addition we have $\|P_k\sigma_1[g(t)]\|_{\Ld^2_x}\lesssim\|\rho_g(t)\|_{\Ld^2_x}$, we can appeal to the Lax--Milgram theorem and deduce that there exists a unique solution $u\in\Ld^\infty(0,T;H^1(\T^3)^3)$ of~\eqref{eq:def-Lambda}. Further recalling~\eqref{eq:apriori-est-ggg}, it satisfies
\begin{equation}\label{eq:estim-u-constr-dis}
\|\nabla u\|_{\Ld^\infty_T\Ld^2_x}\,\lesssim\,\|\rho_g\|_{\Ld^\infty_T\Ld^2_x}\,\lesssim\,e^{CT}\|\rho_{f^\circ}\|_{\Ld^2_x}.
\end{equation}
Letting $\Lambda_k(v):=u$ be the solution of~\eqref{eq:def-Lambda}, this shows that we are indeed led to a well-defined map $\Lambda_k:E\to E$.

\medskip
\substep{2.2} Proof that the map $\Lambda_k:E\to E$ is compact.\\
Let $(v_r)_r$ be a bounded sequence in~$E$, and for all $r$ let $(u_r:=\Lambda_k v_r,g_r)$ be the corresponding solution of the system~\eqref{eq:def-Lambda} with $v$ replaced by $v_r$.
By~\eqref{eq:estim-u-constr-dis}, the sequence $(u_r)_r$ is bounded in $\Ld^{\infty}(0,T; H^1(\T^3)^3)$.
By the Aubin--Lions lemma, in order to prove that it is actually precompact in~$\Ld^{\infty}(0,T; H^1(\T^3)^3)$,
it suffices to check that $(u_r)_r$ is also bounded in $\Ld^{\infty}(0,T; H^2(\T^3)^3)$ and that $(\partial_tu_r)_r$ is bounded for instance in $\Ld^{8/3}(0,T; H^1(\T^3)^3)$.
On the one hand, noting that by definition we have $u_r=P_ku_r$,
and appealing to~\eqref{eq:Hs-bnd-truncPk},
we directly find
\[\| u_r \|_{H^2_x}\,=\,\| P_ku_r \|_{H^2_x}\,\lesssim_k\,\| u_r \|_{H^1_x},\]
which shows that $(u_r)_r$ is indeed also bounded in $\Ld^\infty(0,T;H^2(\T^3)^3)$.
We turn to the boundedness of time derivatives.
Taking the time derivative of the elliptic equation for $u_r$, we find
\begin{multline*}
-\Delta( \partial_t u_r)-P_k \Div(\sigma_2[ g_r,\nabla P_k (\partial_tu_r)]) + \nabla(\partial_tp_r) \\
\,= \,  P_k \Div(\sigma_1[\partial_t g_r])+P_k \Div(\sigma_2[\partial_t g_r,\nabla P_k  u_r]).
\end{multline*}
Arguing as for~\eqref{eq:estim-u-constr-dis}, we deduce
\begin{equation}\label{eq:estom-dtur}
\|\partial_t u_r\|_{H^1_x}\,\lesssim\,\|P_k\sigma_1[\partial_t g_r]\|_{\Ld^2_x}+\|P_k \sigma_2[\partial_t g_r,\nabla P_k u_r]\|_{\Ld^2_x}.
\end{equation}
To estimate the two right-hand side terms, we recall the definition of $\sigma_1,\sigma_2$, we insert the equation for $g_r$, we integrate by parts in the $n$-integrals, and we appeal to~\eqref{eq:Hs-bnd-truncPk} again. For instance, for the term involving $\sigma_1$, we find
\begin{eqnarray*}
\|P_k\sigma_1[\partial_t g_r]\|_{\Ld^2_x}
&=& \lambda\theta\Big\|P_k\int_{\Sp^{2}}\big(n\otimes n-\tfrac1d\Id\big)\Big(\Div_x\big((v_r+U_0n)g_r\big)\\
&&\hspace{1.5cm}+\Div_n\big(\pi_n^\bot( \nabla P_k v_r)n g_r\big)-\Delta_xg_r-\Delta_ng_r\Big)\,\dd n\Big\|_{\Ld^2_x}\\
&=& \lambda\theta\bigg\|
\sum_{i=1}^3P_k\nabla_{x_i}\int_{\Sp^{2}}\big(n\otimes n-\tfrac1d\Id\big)\big((v_r+U_0n)g_r\big)_i\,\dd n\\
&&\qquad-P_k\int_{\Sp^{2}}\nabla_n(n\otimes n)\pi_n^\bot( \nabla P_k v_r)n g_r\,\dd n\\
&&\qquad-P_k\Delta_x\int_{\Sp^{2}}\big(n\otimes n-\tfrac1d\Id\big)\,g_r\,\dd n
-P_k\int_{\Sp^{2}}\Delta_n(n\otimes n)\,g_r\,\dd n\bigg\|_{\Ld^2_x}\\
& \lesssim_k&
\Vert v_r \Vert_{\Ld^4_x}\|g_r\|_{\Ld^4_x\Ld^2_n}+\|g_r\|_{\Ld^2_x\Ld^2_n}.
\end{eqnarray*}
Noting that we have $\|P_k(h_1P_kh_2)\|_{\Ld^2_x}\le\|(P_{2k}h_1)(P_kh_2)\|_{\Ld^2_x}$ for all $h_1,h_2\in\Ld^1(\T^3)$,
we can argue similarly to estimate the second term in~\eqref{eq:estom-dtur}. Further using Jensen's inequality, we are led to
\begin{equation*}
\|\partial_t u_r\|_{H^1_x}\,\lesssim_k\,
(1+\|u_r\|_{\Ld^2_x})(1+\|v_r\|_{\Ld^4_x})\|g_r\|_{\Ld^4_x\Ld^2_n}.
\end{equation*}
By Ladyzhenskaya's inequality, this yields
\[\|\partial_tu_r\|_{H^1_x}\,\lesssim_k\,(1+\|u_r\|_{\Ld^2_x})(1+\|v_r\|_{H^1_x})\|g_r\|_{\Ld^2_{x,n}}^\frac14\| g_r\|_{H^1_{x}\Ld^2_n}^\frac34.\]
The a priori estimate~\eqref{eq:apriori-est-ggg} for $g=g_r$ then ensures that the sequence $(\partial_t u_r)_r$ is bounded in $\Ld^{8/3}(0,T;H^1(\T^3)^3)$,
and the compactness of $\Lambda_k$ follows.

\medskip
\substep{2.3} Proof that the map $\Lambda_k:E\to E$ is continuous.\\
Given a sequence $(v_r)_r$ that converges (strongly) to some $v$ in $E$,
we need to check that the image $u_r:=\Lambda_kv_r$ converges to~$\Lambda_kv$.
By compactness of $\Lambda_k$, we already know that up to a subsequence $u_r=\Lambda_kv_r$ converges (strongly) to some $w$ in $E$, and it remains to show that it satisfies $w=\Lambda_kv$. By definition of~$\Lambda_k$, we recall that $u_r$ satisfies the system~\eqref{eq:def-Lambda} with~$v$ replaced by~$v_r$, and we denote by~$g_r$ the corresponding density. Using the strong convergence of~$u_r$ and~$v_r$ along the extracted subsequence, recalling the a priori estimates for~$g_r$, cf.~\eqref{eq:apriori-est-ggg}, and using weak compactness for~$g_r$, we can pass to the limit in the weak formulation of this system. This shows that the limit~$w$ satisfies the system defining~$\Lambda_k v$, hence we have $w=\Lambda_kv$ by uniqueness.

\medskip
\substep{2.4} Application of Schauder's fixed-point theorem.\\
Recall the a priori estimate~\eqref{eq:estim-u-constr-dis}, that is,
\begin{equation*}
\| \Lambda_k v\|_{\Ld^\infty_TH^1_x}\,\le\,Ce^{CT}\|\rho_{f^\circ}\|_{\Ld^2_x}.
\end{equation*}
Considering the following non-empty closed convex subset of $E$,
\[K\,:=\,\big\{ u \in E : \| u \|_{\Ld^{\infty}_TH^1_x} \le Ce^{CT}\|\rho_{f^\circ}\|_{\Ld^2_x} \big\},\]
the restriction of $\Lambda_k$ defines a map $\Lambda_k|_K:K\to K$ that is compact and continuous by Steps~2.2 and~2.3.
By Schauder's theorem, this map must then admit a fixed point $u_k \in K$, that is, $u_k=\Lambda_k(u_k)$. Combined with the a priori estimates~\eqref{eq:apriori-est-ggg}, this concludes the proof of the existence of a weak solution of the approximate system~\eqref{system:regularized-plop} satisfying~\eqref{system:regularized-plop-regul}.

\medskip
\step{3} Existence of a weak solution for the system~\eqref{eq:eps-kin00}.\\
We argue by means of a Galerkin method based on the approximations defined in Step~2.
For that purpose, given an initial condition $f^\circ\in\Ld^2\cap\Pc(\T^3\times\Sp^2)$, we start by considering an approximating sequence of smooth initial conditions $(f^\circ_k)_k\subset C^\infty \cap\Pc(\T^3\times\Sp^2)$ that converges strongly to $f^\circ$ in $\Ld^2(\T^3\times\Sp^2)$.
By Step~2, for all $k$, we may then consider the solution $(u_k,f_k)$ of the approximate system~\eqref{system:regularized-plop} with initial condition~$f_k^\circ$, that is,
\begin{equation*}
\left\{\begin{array}{l}
-\Delta u_k+\nabla p_k
=P_k \Div(\sigma_1[f])+P_k \Div(\sigma_2[f_k,\nabla P_k  u_k]),\\[1mm]
\partial_t f_k+\Div_x\big((u_k+U_0n)f_k\big)+\Div_n\big(\pi_n^\bot( \nabla P_ku_k)n f_k\big)=\Delta_xf_k+\Delta_nf_k, \\[1mm]
\Div(u_k)=0,
\quad\textstyle \int_{\T^3} u_k = 0,\\[1mm]
f_k|_{t=0}=f_k^\circ,
\end{array}\right.
\end{equation*}
with
\begin{eqnarray*}
u_k &\in& \Ld^\infty(0,T; H^1(\T^3)^3),\\
f_k &\in& \Ld^{\infty}(0,T; \Ld^{2}\cap\Pc(\T^3 \times \mathbb{S}^2)) \cap \Ld^{2}(0,T; H^1(\T^3 \times \Sp^2)).
\end{eqnarray*} 
From Step~2, we further learn that $(u_k,f_k)$ is uniformly bounded in these spaces.
By weak compactness, we deduce that up to a subsequence $u_k$ converges weakly-* to some $u$ in $\Ld^\infty(0,T;H^1(\T^3))$, and that $f_k$ converges weakly-* to some $f$ in $\Ld^\infty(0,T;\Ld^2(\T^3\times\Sp^2))$ and weakly in $\Ld^2(0,T;H^1(\T^3\times\Sp^2))$.
In addition, examining the equation satisfied by $f_k$, a simple argument allows to check that~$\partial_tf_k$ is bounded e.g.\@ in $\Ld^\infty(0,T;W^{-2,1}(\T^3\times\Sp^2))$, so that the Aubin--Lions lemma further entails that $f_k$ converges strongly to $f$ in $\Ld^2(0,T;\Ld^2(\T^3)^3)$.
These convergences now allow to pass to the limit in the weak formulation of the system, and we easily conclude that the extracted limit $(u,f)$ satisfies the limiting system~\eqref{eq:eps-kin00}.

\medskip
\step{4} Weak-strong uniqueness for the system~\eqref{eq:eps-kin00}.\\
Let $(u,f)$ and $(u',f')$ be two weak solutions of~\eqref{eq:eps-kin00} with common initial condition
\[f|_{t=0}=f'|_{t=0}=f^{\circ},\]
and assume that $(u',f')$ further satisfies
\begin{eqnarray}
u' &\in& \Ld^2_{\loc}(\R^+; W^{1, \infty}(\T^3)^3),\label{Assump:weakstrong}\\
f' &\in& \Ld^{\infty}_{\loc}(\R^+; \Ld^{\infty}(\T^3\times\Sp^2)).\nonumber
\end{eqnarray}
Consider the differences $U=u-u'$ and $F=f-f'$. On the one hand, the equations for $u$ and $u'$ yield
\begin{eqnarray*}
- \Delta U+ \nabla P
&=& \Div(\sigma_1[F])+\Div(\sigma_2[f,\nabla u]-\sigma_2[f',\nabla u']) \\
& =& \Div(\sigma_1[F])+\Div(\sigma_2[f, \nabla U]+\sigma_2[F, \nabla u']).
\end{eqnarray*}
Testing this equation with $U$ itself, using the incompressibility constraint, and taking advantage of the additional dissipation given by $\sigma_2$, we get
\begin{eqnarray}
\|\nabla U\|_{\Ld^2_x}
&\le&\|\sigma_1[F]\|_{\Ld^2_x}+\|\sigma_2[F,\nabla u']\|_{\Ld^2_x}\nonumber\\
&\le&\big(1+\|\nabla u'\|_{\Ld^\infty_x}\big)\|F\|_{\Ld^2_{x,n}}.\label{elliptic:estimUniq}
\end{eqnarray}
On the other hand, the equations for $f$ and $f'$ yield
\begin{multline*}
\partial_tF+\Div_x\big((u+U_0n)F\big)+\Div_n\big(\pi_n^\bot(\nabla u)nF\big)-\Delta_xF-\Delta_nF\\
\,=\,
-\Div_x(Uf')
-\Div_n\big(\pi_n^\bot(\nabla U)nf'\big).
\end{multline*}
Testing this equation with $F$ itself, and integrating by parts, we get
\begin{eqnarray*}
\lefteqn{\tfrac12\tfrac{\dd}{\dd t}\Vert F \Vert^2_{\Ld^2_{x,n}} + \Vert \nabla_{x} F\Vert_{\Ld^2_{x,n}}^2+ \Vert \nabla_{n} F\Vert_{\Ld^2_{x,n}}^2}\\
&=&-\tfrac12\int_{\T^3\times\Sp^2}|F|^2\Div_n\big(\pi_n^\bot(\nabla u)n\big)+\int_{\T^3\times\Sp^2}\nabla_xF\cdot Uf'+\int_{\T^3\times\Sp^2}\nabla_nF\cdot \pi_n^\bot(\nabla U)nf'\\
&\lesssim&\|\nabla u\|_{\Ld^2_x}\|F\|_{\Ld^4_x\Ld^2_n}^2+\|\nabla_xF\|_{\Ld^2_{x,n}}\|U\|_{\Ld^2_x}\|f'\|_{\Ld^\infty_{x,n}}+\|\nabla_nF\|_{\Ld^2_{x,n}}\|\nabla U\|_{\Ld^2_x}\|f'\|_{\Ld^\infty_{x,n}}.
\end{eqnarray*}
By Ladyzhenskaya's inequality, the first right-hand side term can be bounded by
\begin{eqnarray*}
\|\nabla u\|_{\Ld^2_x}\|F\|_{\Ld^4_x\Ld^2_n}^2&\lesssim&\|\nabla u\|_{\Ld^2_x}\|F\|_{\Ld^2_{x,n}}^\frac12\|F\|_{H^1_x\Ld^2_n}^\frac32\\
&\lesssim&\|\nabla u\|_{\Ld^2_x}\|F\|_{\Ld^2_{x,n}}^2
+\|\nabla u\|_{\Ld^2_x}\|F\|_{\Ld^2_{x,n}}^\frac12\|\nabla_xF\|_{\Ld^2_{x,n}}^\frac32.
\end{eqnarray*}
Inserting this into the above and appealing to Young's and Poincar\'e's inequalities, we are led to
\begin{equation*}
\tfrac{\dd}{\dd t}\Vert F \Vert^2_{\Ld^2_{x,n}}
\,\lesssim\,\big(\|\nabla u\|_{\Ld^2_x}
+\|\nabla u\|_{\Ld^2_x}^4\big)\|F\|_{\Ld^2_{x,n}}^2
+\|\nabla U\|_{\Ld^2_x}^2\|f'\|_{\Ld^\infty_{x,n}}^2,
\end{equation*}
and thus, recalling~\eqref{elliptic:estimUniq},
\begin{equation*}
\tfrac{\dd}{\dd t}\Vert F \Vert^2_{\Ld^2_{x,n}}
\,\lesssim\,\Big(\big(1
+\|\nabla u\|_{\Ld^2_x}^4\big)+\big(1+\|\nabla u'\|_{\Ld^\infty_x}^2\big)\|f'\|_{\Ld^\infty_{x,n}}^2\Big)\|F\|_{\Ld^2_{x,n}}^2.
\end{equation*}
By the assumed regularity of $u,u',f'$, this implies $F=0$ by Gr\"onwall's inequality, hence also~$U=0$ by~\eqref{elliptic:estimUniq}.
\qed

\section{Perturbative well-posedness of ordered fluid equations}\label{app.Bous}
{This appendix is devoted to the perturbative well-posedness of the second-order fluid model~\eqref{eq:second-order-fluid-intro-eqn} that we derive from the Doi--Saintillan--Shelley model. We recall that this model is a priori ill-posed as it behaves like a backward heat equation, but we show that well-posed notions of {\it approximate} solutions can be defined.
We start below with the proof of Proposition~\ref{prop:hier-inhom} on so-called approximate hierarchical solutions, which are motivated by the formal $\e$-expansion.

Next, we show that
there are also alternative, non-hierarchical ways to perturbatively make sense of~\eqref{eq:second-order-fluid-intro-eqn}, which may be more desirable in particular for stability issues. For shortness, we shall focus here on {\it homogeneous} second-order fluids at $\Pe<\infty$ as defined in Section~\ref{sec:nonstandard-ordered-Pe}:
more precisely, instead of~\eqref{eq:second-order-fluid-intro-eqn}, we consider~\eqref{ordered.fluid-Pe}, that is,
\begin{equation}\label{ordered.fluid-Pe-RE}
\left\{\begin{array}{l} 
\Rey(\partial_t+\bar u_\e\cdot\nabla)\bar u_\e-\Div(\bar\sigma_\e)+\nabla\bar p_\e=h+O(\e^2),\\[1mm]
\bar\sigma_\e\,=\,\eta_0A_1(\bar u_\e)+\e\gamma_1A_2'(\bar u_\e)+\e\gamma_2A_1(\bar u_\e)^2,\\[1mm]
\Div(\bar u_\e)=0.
\end{array}\right.
\end{equation}
Comparing} to corresponding ill-posedness issues in the Boussinesq theory for water waves, we recall that there is a standard way to rearrange the ill-posed Boussinesq equation perturbatively and make it well-posed, see~\cite{Christov-Maugin-Velarde}: in a nutshell, the idea is to replace indefinite operators like $1+\e\Delta$ by corresponding positive operators like $(1-\e\Delta)^{-1}$ up to~$O(\e^2)$ errors.
We notice that a similar so-called {\it Boussinesq trick} can be used in the present situation as well: for {any} value of $\gamma_1,\gamma_2$, both at finite and infinite P\'eclet number, it leads us to a perturbative rearrangement of the second-order fluid equation that is well-posed and indeed equivalent to~\eqref{ordered.fluid-Pe-RE} up to~$O(\e^2)$ terms.
The so-defined solution is easily checked to differ from the corresponding hierarchical solution only by $O(\e^2)$.
We focus here on the relevant ill-posed case $\gamma_1\le0$.
Note that this procedure is easily extended to higher-order fluid equations (see also~\cite{Allaire-Briane-Vanninathan,duerinckx2023spectral} in a different context).

{\begin{prop}[Boussinesq-like solutions]\label{prop:Bous}
Consider the system~\eqref{ordered.fluid-Pe-RE} in the regime $\e\ll1$ with parameters $\eta_0, \Pe>0$, $\gamma_1\le0$, and $\gamma_2\in\R$.
\begin{enumerate}[(i)]
\item \emph{Stokes case $\Rey=0$, $d\le3$:}\\
Given $s>\frac d2$, $T_0>0$, and $h\in \Ld^\infty(0,T_0;H^{s+1}(\T^d)^d)\cap W^{1,\infty}(0,T_0;H^{s-1}(\T^d)^d)$, provided that $\e\ll1$ is small enough (depending on $s,T_0,h$ and on all parameters), the following nonlinear problem admits a unique solution $\bar u_{\e}\in\Ld^\infty(0,T_0;H^{s+1}(\T^d)^d)$,
\begin{equation}\label{eq:lim-St}
\qquad\left\{\begin{array}{l}
-\eta_0\Delta\bar u_\e+\nabla\bar p_\e=\big(1-\e\tfrac{\gamma_1}{\eta_0}(\partial_t-\tfrac1\Pe\Delta)\big)h+\e\Div(F_0(\bar u_\e)),\\[1mm]
\Div(\bar u_\e)=0,\quad \int_{\T^d} \bar u_\e=0,
\end{array}\right.
\end{equation}
where we have set for abbreviation
\begin{equation}\label{eq:def-F0hom}
\qquad F_0(u)\,:=\,\gamma_1(u\cdot\nabla)2\!\D(u)+\gamma_1\big((\nabla u)^T2\!\D(u)+2\!\D(u)(\nabla u)\big)+\gamma_2(2\!\D(u))^2.
\end{equation}
Moreover, the so-defined solution $\bar u_\e$ satisfies~\eqref{ordered.fluid-Pe-RE} with $\Rey=0$ for some controlled error term~$O(\e^2)$.

\smallskip\item \emph{Navier--Stokes case $\Rey=1$, $d=2$:}\\
Further assume $\eta_0\ge\frac1\Pe$.\footnote{For $\eta_0<\frac1\Pe$, the equations would need to be rearranged differently; we skip it here for shortness.}
Given $s>1$, $T_0>0$, $u^\circ\in H^s(\T^2)^2$ with \mbox{$\Div(u^\circ)=0$}, and $h\in \Ld^2(0,T_0;H^{s+1}(\T^2)^2)$, provided that $\e\ll1$ is small enough (depending on~$s,T_0,u^\circ,h$ and on all parameters), the following nonlinear problem admits a unique solution $\bar u_{\e}\in\Ld^\infty(0,T_0;H^{s}(\T^2)^2)\cap\Ld^2(0,T_0;H^{s+1}(\T^2)^2)$,
\begin{equation}\label{eq:lim-NS}
\qquad\left\{\begin{array}{l}
(\partial_t+\bar u_\e\cdot\nabla)\bar u_\e-\eta_0\Delta\bar u_{\e}-\e\gamma_1(\eta_0-\frac1{\Pe})\Delta^2\bar u_\e+\nabla\bar p_{\e}\\
\hspace{5cm}\,=\,(1+\e\gamma_1\Delta)h+\e\Div(F_1(\bar u_\e)),\\
\Div(\bar u_\e)=0,\\
\bar u_\e|_{t=0}=u^\circ,
\end{array}\right.
\end{equation}
where we have set for abbreviation
\begin{equation*}
\qquad F_1(u)\,:=\,2\gamma_1(\nabla u)^T(\nabla u)+\gamma_2(2\!\D(u))^2.
\end{equation*}
Moreover, the so-defined solution $\bar u_\e$ satisfies~\eqref{ordered.fluid-Pe-RE} with $\Rey=1$ with some controlled error term~$O(\e^2)$ and with initial condition $\bar u_\e|_{t=0}=u^\circ$.\qedhere
\end{enumerate}
\end{prop}}

\subsection{Proof of Proposition~\ref{prop:hier-inhom}}
Let $(u_0,\rho_0)$ and $(u_1,\rho_1)$ satisfy the systems~\eqref{eq:v0rho0} and~\eqref{eq:v1rho1}, respectively. Their superposition $(\bar u_\e,\bar \rho_\e)=(u_0+\e u_1,\rho_0 + \e \rho_1)$ is easily checked to satisfy the non-homogeneous second-order fluid equations~\eqref{eq:second-order-fluid-intro-eqn} in the form
\begin{equation*}
\left\{\begin{array}{l} 
\Rey(\partial_t+\bar u_\e\cdot\nabla)\bar u_\e-\Div(\bar\sigma_\e)+\nabla\bar p_\e=h+ \e^2 \mathcal{R}_\e,\\[1mm]
\partial_t\bar \rho_\e-(\tfrac1\Pe+ \e \mu_0)\Delta\bar \rho_\e+\bar u_\e\cdot\nabla\bar \rho_\e= \e^2 \mathcal{S}_\e,\\[1mm]
\bar\sigma_\e\,=\,(\eta_0+\eta_1\bar \rho_\e)A_1(\bar u_\e)+\e\gamma_1A_2'(\bar\rho_\e,\bar u_\e)+\e\gamma_2A_1(\bar u_\e)^2,\\[1mm]
\Div(\bar u_\e)=0,\\[1mm]
(\bar \rho_\e,\bar u_\e)|_{t=0}=(\rho^\circ,u^\circ),
\end{array}\right.
\end{equation*}
with explicit remainders given by 
\begin{eqnarray*}
\mathcal{R}_\e & :=& - \lambda \theta \tfrac{\omega_d}{d(d+2)} \Div (\rho_1  \D(u_1))\\
&&- \gamma_1 \tfrac{1}{\e} \Div \big( A_2'(\bar u_\e,\bar \rho_\e) - A_2'(u_0,\rho_0) \big) - \lambda \theta \tfrac{2 \omega_d}{d^2(d+4)} \tfrac{1}{\e} \Div \big( \bar \rho_\e \D(\bar u_\e)^2 - \rho_0 \D(u_0)^2 \big)\\
&& - \lambda \tfrac{2 \omega_d}{d(d+2)(d+4)} \Div \Big( \rho_0 \big( \D(\bar u_\e) \D(u_1) + \D(u_1) \D(\bar u_\e) \big) +  2 \rho_1 \D(\bar u_\e)^2 \Big)\\
&& + \Rey (u_1 \cdot \nabla) u_1 - \Div (\sigma_2[\rho_1+g_1,\nabla u_1]),\\
\mathcal{S}_\e &:=& - \tfrac{1}{d(d-1)} U_0^2 \Delta \rho_1 + u_1 \cdot \nabla \rho_1.
\end{eqnarray*}
In the proof of Proposition~\ref{prop:main}, we recall that we have actually shown that after elimination of $g_1,g_2$ the solutions $(u_0,\rho_1)$ and $(u_1,\rho_2)$ of the hierarchy~\eqref{eq :(u_0,f_1)}--\eqref{eq :(u_1,f_2)} precisely satisfy equations~\eqref{eq:v0rho0}--\eqref{eq:v1rho1} (for a specific set of parameters).
Hence, the well-posedness of hierarchical solutions of~\eqref{eq:v0rho0}--\eqref{eq:v1rho1} follows from the proof of Proposition~\ref{prop:first-terms}, which indeed concerns the well-posedness of the hierarchy~\eqref{eq :(u_0,f_1)}--\eqref{eq :(u_1,f_2)} in the Stokes case $\Rey=0$. The proof is analogous in the 2D Navier--Stokes case and is skipped for conciseness.
\qed

\subsection{Proof of Proposition~\ref{prop:Bous} in Stokes case}
Let $\Rey=0$ and $d\le3$.
We split the proof into three steps, first deriving the suitable reformulation of the equations, and then proving the existence and uniqueness of smooth solutions.

\medskip
\step1 Reformulation.\\
Equation~\eqref{eq:lim-St} implies in particular
\[h\,=\,-\eta_0\Delta\bar u_\e-\e\Div(F_0(\bar u_\e))+\e\tfrac{\gamma_1}{\eta_0}(\partial_t-\tfrac1\Pe)h+\nabla\bar p_\e,\]
and therefore, using this to replace $\e\frac{\gamma_1}{\eta_0}(\partial_t-\frac1{\Pe})h$ in the right-hand side of~\eqref{eq:lim-St},
\begin{multline*}
-\eta_0\Delta\bar u_\e+\nabla\bar p_\e=h+\e\Div(F_0(\bar u_\e))\\
-\e\tfrac{\gamma_1}{\eta_0}(\partial_t-\tfrac1\Pe)\Big(-\eta_0\Delta\bar u_\e+\nabla\bar p_\e-\e\Div(F_0(\bar u_\e))+\e\tfrac{\gamma_1}{\eta_0}(\partial_t-\tfrac1\Pe)h\Big).
\end{multline*}
Recalling the definition of $F_0$, cf.~\eqref{eq:def-F0hom}, and reorganizing $O(\e)$ terms, we precisely deduce the second-order fluid equation
\begin{equation*}
\left\{\begin{array}{l}
-\Div(\bar \sigma_\e)+\nabla\bar p_\e'=h+\e^2R_\e,\\[1mm]
\bar\sigma_\e=\eta_0A_1(\bar u_\e)+\e\gamma_1A_2'(\bar u_\e)+\e\gamma_2A_1(\bar u_\e)^2,
\end{array}\right.
\end{equation*}
for some modified pressure field $\bar p_\e'$ and some remainder term
\[R_\e\,:=\,\tfrac{\gamma_1}{\eta_0}(\partial_t-\tfrac1\Pe)\Big(\Div(F_0(\bar u_\e))-\tfrac{\gamma_1}{\eta_0}(\partial_t-\tfrac1\Pe)h\Big).\]

\medskip
\step2 Existence.\\
Let $s>\frac d2$ be fixed.
We proceed by an iterative scheme.
We set $\bar u_0:=0$ and for all $n\ge0$, given $\bar u_n\in H^{s+1}(\T^d)^d$, we define~$\bar u_{n+1}\in H^{s+1}(\T^d)^d$ as the unique solution of the linear problem
\begin{equation}\label{eq:iterative-scheme-St}
\left\{\begin{array}{l}
-\eta_0\Delta\bar u_{n+1}-\e\gamma_1\Div\big((\bar u_n\cdot\nabla)2\!\D(\bar u_{n+1})\big)+\nabla \bar p_{n+1}\\
\hspace{5cm}=\big(1-\e\tfrac{\gamma_1}{\eta_0}(\partial_t-\Delta)\big)h+\e\Div(G_0(\bar u_n)),\\
\Div(\bar u_{n+1})=0,
\end{array}\right.
\end{equation}
where we have set
\[G_0(u)\,:=\,\gamma_1\big((\nabla u)^T2\!\D(u)+2\!\D(u)(\nabla u)\big)+\gamma_2(2\!\D(u))^2.\]
We split the proof into two further substeps.

\medskip
\substep{2.1} Sobolev a priori estimates.\\
Applying $\langle\nabla\rangle^s=(1-\Delta)^{s/2}$ to both sides of the above equation, testing it with $\langle\nabla\rangle^s\bar u_{n+1}$, and using the incompressibility constraint, we find
\begin{multline*}
\eta_0\int_{\T^d}|\nabla\langle\nabla\rangle^s \bar u_{n+1}|^2
\,=\,\int_{\T^d}\langle\nabla\rangle^s\bar u_{n+1}\cdot\big(1-\e\tfrac{\gamma_1}{\eta_0}(\partial_t-\Delta)\big)\langle\nabla\rangle^sh\\
-\e\int_{\T^d}\nabla\langle\nabla\rangle^s\bar u_{n+1}: \langle\nabla\rangle^sG_0(\bar u_n)
-2\e\gamma_1\int_{\T^d}\D(\langle\nabla\rangle^s\bar u_{n+1}):[\langle\nabla\rangle^s,\bar u_n\cdot\nabla]\!\D(\bar u_{n+1}),
\end{multline*}
hence, by the Cauchy--Schwarz inequality,
\begin{multline}\label{eq:estim-Un+1}
\|\nabla\bar u_{n+1}\|_{H^s_x}
\,\lesssim\,\|h\|_{H^{s-1}_x}+\e\|(\partial_t-\Delta)h\|_{H^{s-1}_x}+\e\|G_0(\bar u_n)\|_{H^s_x}\\
+\e \|[\langle\nabla\rangle^s,\bar u_n\cdot\nabla]\!\D(\bar u_{n+1})\|_{\Ld^2_x}.
\end{multline}
To estimate the last term, we use the following form of the Kato--Ponce commutator estimate~\cite[Lemma~X1]{KatoPonce}: for all $u,v\in C^\infty(\T^d)^d$ and all $p,q\in[2,\infty]$ with $\frac1p+\frac1q=\frac12$, we have
\begin{equation}\label{eq:KatoPonce-est}
\|[\langle\nabla\rangle^s,u\cdot\nabla]\!\D(v)\|_{\Ld^2_x}\,\lesssim_s\,\|u\|_{W^{s,q}_x}\|v\|_{W^{2,p}_x}+\|u\|_{W^{1,\infty}_x}\|v\|_{H^{s+1}_x},
\end{equation}
and thus, properly choosing $p,q$ and appealing to the Sobolev embedding with $s>\frac d2$, we deduce
\begin{equation*}
\|[\langle\nabla\rangle^s,u\cdot\nabla]\!\D(v)\|_{\Ld^2_x}\,\lesssim_s\,\|u\|_{H^{s+1}_x}\|v\|_{H^{s+1}_x}.
\end{equation*}
Using this to estimate the last right-hand side term in~\eqref{eq:estim-Un+1},
inserting the definition of $G_0$, and further appealing to the Sobolev embedding with $s>\frac d2$, we are led to
\begin{equation*}
\|\bar u_{n+1}\|_{H^{s+1}_x}
\,\lesssim_s\,\|h\|_{H^{s-1}_x}+\e\|(\partial_t-\Delta)h\|_{H^{s-1}_x}+\e\|\bar u_n\|_{H^{s+1}_x}^2
+\e \|\bar u_n\|_{H^{s+1}_x}\|\bar u_{n+1}\|_{H^{s+1}_x}.
\end{equation*}
Provided that $\|\bar u_n\|_{H^{s+1}_x}\le C_0$ and that $\e C_0\ll_s1$ is small enough, the last right-hand side term can be absorbed and we are led to
\begin{equation*}
\|\bar u_{n+1}\|_{H^{s+1}_x}
\,\lesssim_s\,\|h\|_{H^{s-1}_x}+\e\|(\partial_t-\Delta)h\|_{H^{s-1}_x}+\e C_0^2.
\end{equation*}
Choosing $C_0:=C_s(\|h\|_{H^{s-1}_x}+\e \|(\partial_t-\Delta)h\|_{H^{s-1}})$ for some $C_s\gg_s1$ large enough, we then deduce, provided that $\e C_0\ll_s1$ is small enough,
\[\|\bar u_n\|_{H^{s+1}_x}\,\le\, C_0\qquad\implies\qquad\|\bar u_{n+1}\|_{H^{s+1}_x}\,\le\,C_0.\]
Recalling the choice $\bar u_0=0$, this proves by induction that we have for all $n\ge0$, provided that $\e\ll_{s,h}1$ is small enough,
\begin{equation}\label{eq:Sobolev-estim}
\|\bar u_{n}\|_{H^{s+1}_x}\,\le\,C_0\,\simeq_s\,\|h\|_{H^{s-1}_x}+\e \|(\partial_t-\Delta)h\|_{H^{s-1}_x}.
\end{equation}

\medskip
\substep{2.2} Contraction.\\
The difference $\bar u_{n+1}-\bar u_n$ satisfies
\begin{multline*}
-\eta_0\Delta (\bar u_{n+1}-\bar u_n)+\nabla(\bar p_{n+1}-\bar p_n)
=\e\gamma_1\Div\big((\bar u_{n}\cdot\nabla)2\!\D(\bar u_{n+1}-\bar u_n)\big)\\
+\e\gamma_1\Div\big(((\bar u_{n}-\bar u_{n-1})\cdot\nabla)2\!\D(\bar u_n)\big)
+\e\Div(G_0(\bar u_{n})-G_0(\bar u_{n-1})).
\end{multline*}
Testing this equation with $\bar u_{n+1}-\bar u_n$ and using the incompressibility constraint, we get similarly as in~\eqref{eq:estim-Un+1} above,
\begin{equation*}
\|\nabla(\bar u_{n+1}-\bar u_n)\|_{\Ld^2_x}
\,\lesssim\,\e\|(\bar u_{n}-\bar u_{n-1})\cdot\nabla\D(\bar u_n)\|_{\Ld^2_x}
+\e\|G_0(\bar u_{n})-G_0(\bar u_{n-1})\|_{\Ld^2_x},
\end{equation*}
and thus, inserting the definition of $G_0$ and appealing to the Sobolev embedding with~$s>\frac d2$,
\begin{equation*}
\|\nabla(\bar u_{n+1}-\bar u_n)\|_{\Ld^2_x}
\,\lesssim_s\,\e\|(\bar u_{n-1},\bar u_{n})\|_{H^{s+1}_x}\|\nabla(\bar u_{n}-\bar u_{n-1})\|_{\Ld^2_x}.
\end{equation*}
Provided that $\e\ll_{s,h}1$ small enough, inserting the a priori estimate~\eqref{eq:Sobolev-estim}, we deduce
\begin{equation*}
\|\nabla(\bar u_{n+1}-\bar u_n)\|_{\Ld^2_x}
\,\le\,\tfrac12\|\nabla(\bar u_{n}-\bar u_{n-1})\|_{\Ld^2_x}.
\end{equation*}
Together with the a priori bound~\eqref{eq:Sobolev-estim}, this contraction estimate entails that the sequence $(\bar u_n)_n$ converges weakly in $H^{s+1}(\T^d)^d$ to some limit $\bar u_\e\in H^{s+1}(\T^d)^d$. By the Rellich theorem, recalling $s>\frac d2$, this allows to pass to the limit in~\eqref{eq:iterative-scheme-St} and to deduce that the limit $\bar u_\e$ is a solution of the nonlinear equation~\eqref{eq:lim-St}. In addition, it automatically satisfies the a priori estimate
\begin{equation}\label{eq:apriori-St}
\|\bar u_\e\|_{H^{s+1}_x}\,\lesssim_s\,\|h\|_{H^{s-1}_x}+\e \|(\partial_t-\Delta)h\|_{H^{s-1}_x}.
\end{equation}

\medskip
\step3 Uniqueness.\\
Let $s>\frac d2$ be fixed.
Let $\bar u_\e,\bar v_\e\in H^{s+1}(\R^d)^d$ be two solutions of~\eqref{eq:lim-St}. Their difference then satisfies
\begin{multline*}
-\eta_0\Delta (\bar u_\e-\bar v_\e)+\nabla \bar p_\e
=\e\gamma_1\Div\big((\bar u_\e\cdot\nabla)2\!\D(\bar u_\e-\bar v_\e)\big)\\
+\e\gamma_1\Div\big(((\bar u_\e-\bar v_\e)\cdot\nabla)2\!\D(\bar v_\e)\big)
+\e\Div(G_0(\bar u_\e)-G_0(\bar v_\e)).
\end{multline*}
Arguing as in Step~2.2 above, we easily deduce
\[\|\nabla(\bar u_\e-\bar v_\e)\|_{\Ld^2_x}^2\,\lesssim\,\e\|(\bar u_\e,\bar v_\e)\|_{H^{s+1}_x}\|\nabla(\bar u_\e-\bar v_\e)\|_{\Ld^2_x}^2.\]
By the a priori estimate~\eqref{eq:apriori-St} for $\bar u_\e,\bar v_\e$, we deduce $\bar u_\e=\bar v_\e$ provided that $\e\ll_{s,h}1$ is small enough.\qed

\subsection{Proof of Proposition~\ref{prop:Bous} in Navier--Stokes case}
Let $\Rey=1$ and $d=2$.
We skip the proof of well-posedness, which is a straightforward modification of the above proof in the corresponding Stokes case.
It remains to derive the suitable reformulation of the equations.
For that purpose, first note that equation~\eqref{eq:lim-NS} can be rewritten as
\begin{multline}\label{eq:reform-lim-NS}
(\partial_t+\bar u_\e\cdot\nabla)\bar u_\e-\eta_0\Delta\bar u_\e
-\e\gamma_1(\eta_0-\tfrac1\Pe)\Delta^2\bar u_\e
+\nabla\bar p_\e\\
=h+\e\gamma_1\Div(2\D(h))
+\e\Div(F_0(\bar u_\e))-\e\gamma_1\Div\big(2\D((\bar u_\e\cdot\nabla)\bar u_\e)\big),
\end{multline}
where $F_0$ is defined in~\eqref{eq:def-F0hom}.
This yields in particular the following relation,
\begin{multline*}
h\,=\,(\partial_t+\bar u_\e\cdot\nabla)\bar u_\e-\eta_0\Delta\bar u_\e
+\nabla\bar p_\e\\
-\e\gamma_1(\eta_0-\tfrac1\Pe)\Delta^2\bar u_\e
-\e\gamma_1\Div(2\D( h))
-\e\Div(F_0(\bar u_\e))+\e\gamma_1\Div\big(2\D((\bar u_\e\cdot\nabla)\bar u_\e)\big).
\end{multline*}
Using this to replace $\e\gamma_1\Div(2\D(h))$ in the right-hand side of~\eqref{eq:reform-lim-NS}, we find after straightforward simplifications
\begin{multline*}
(\partial_t+\bar u_\e\cdot\nabla)\bar u_\e-\eta_0\Delta\bar u_\e
+\nabla\bar p_\e
=h+\e\Div(F_0(\bar u_\e))\\
+\e\gamma_1\Div\Big[2\D\Big((\partial_t-\tfrac1\Pe\Delta)\bar u_\e
+\nabla\hat p_\e
-\e\gamma_1(\eta_0-\tfrac1\Pe)\Delta^2\bar u_\e\\
-\e\gamma_1\Div(2\D(h))
-\e\Div(F_0(\bar u_\e))+\e\gamma_1\Div\big(2\D((\bar u_\e\cdot\nabla)\bar u_\e)\big)\Big)\Big].
\end{multline*}
Recalling the definition of $F_0$, cf.~\eqref{eq:def-F0hom}, and reorganizing $O(\e)$ terms, we precisely deduce the second-order fluid equation
\begin{equation*}
\left\{\begin{array}{l}
(\partial_t+\bar u_\e\cdot\nabla)\bar u_\e-\Div(\bar \sigma_\e)+\nabla\bar p_\e'=h+\e^2R_\e,\\
\bar \sigma_\e=\eta_0A_1(\bar u_\e)+\e\gamma_1A_2'(\bar u_\e)+\e\gamma_2A_1(\bar u_\e)^2,
\end{array}\right.
\end{equation*}
for some modified pressure field $P_\e'$ and some remainder term
\begin{multline*}
R_\e\,:=\,-\gamma_1\Div\Big[2\D\Big(\gamma_1(\eta_0-\tfrac1\Pe)\Delta^2\bar u_\e
+\gamma_1\Div(2\D(h))
+\Div(F_0(\bar u_\e))\\
-\gamma_1 \Div\big(2\D((\bar u_\e\cdot\nabla)\bar u_\e)\big)\Big)\Big].
\end{multline*}
This concludes the proof.
\qed

\section{Derivation of third-order fluid equations}\label{app:3rd-order}
In this section, we extend Proposition~\ref{lem:Bous}
and derive the corresponding equations to next order as stated in Section~\ref{secrem:third-order} --- we keep the derivation formal and skip detailed error estimates for shortness.
By definition of $u_0,u_1,u_2,\rho_0,\rho_1,\rho_2$ in Proposition~\ref{prop:first-terms}, we first note that $(\bar u_\e,\bar\rho_\e)$ now satisfies
\begin{equation}\label{eq:baru-barrho-2}
\left\{\begin{array}{l}
\Rey(\partial_t +\bar u_\e\cdot\nabla)\bar u_\e-\Delta\bar u_\e+\nabla\bar p_\e=h+\Div(\sigma_1[g_1+\e g_2+\e^2g_3])\\
\hspace{5cm}+\Div(\sigma_2[\bar\rho_\e+\e g_1+\e^2g_2,\nabla\bar u_\e])+O(\e^3),\\
(\partial_t-{\tfrac 1 \Pe}\Delta+\bar u_\e\cdot\nabla)\bar\rho_\e=-\langle U_0n\cdot\nabla_x(\e g_1+\e^2g_2)\rangle+O(\e^3),
\end{array}\right.
\end{equation}
and it remains to evaluate the different contributions of $g_1,g_2,g_3$ in the right-hand side.

\subsection{Equation for particle density}
Using $n=-\tfrac1{d-1}\Delta_n n$, we have
\begin{eqnarray*}
\langle U_0n\cdot\nabla_x(\e g_1+\e^2g_2)\rangle
&=&\e U_0\,\Div\Big(\fint_{\Sp^{d-1}}n(g_1+\e g_2)\,\dd n\Big)
\\
&=&-\e\tfrac{U_0}{d-1} \,\Div\Big(\fint_{\Sp^{d-1}}n\Delta_n(g_1+\e g_2)\,\dd n\Big).
\end{eqnarray*}
Using the defining equations for $g_1,g_2$, we find
\begin{align*}
\Delta_n(g_1+\e g_2) & \, = \,  U_0n\cdot\nabla_x\rho_0+\Div_n\big(\pi_n^\bot(\nabla u_0)n\rho_0\big)\\
& \quad +\e(\partial_t-\tfrac1 \Pe \Delta_x+u_0\cdot\nabla_x)g_1
+\e P_1^\bot\big(U_0n\cdot\nabla_x(\rho_1+g_1)\big)\\
& \quad  +\e\Div_n\big(\pi_n^\bot(\nabla u_1)n\rho_0+\pi_n^\bot(\nabla u_0)n(\rho_1+g_1)\big).
\end{align*}
Inserting this identity into the above, as well as the explicit expression for $g_1$ in Proposition~\ref{prop:first-terms}, then recalling~\eqref{spherical.divergence} and~\eqref{cancel-u_0}, and using integral computations~\eqref{eq:explicit-integral-n},
\begin{align*}
\langle U_0n\cdot\nabla_x(\e g_1+\e^2g_2)\rangle
& \, =  -\e\tfrac{U_0^2}{d(d-1)} \,\Delta\bar\rho_\e\\
& \quad \ -\e^2\tfrac{2U_0^2}{d(d-1)(d+2)} \,\Div\Big(\tfrac{d+1}{d-1}\D(\bar u_\e)\nabla\bar\rho_\e+\tfrac{1}2\Div(\bar\rho_\e\D(\bar u_\e))\Big)
+O(\e^3).
\end{align*}
Inserting this computation into~\eqref{eq:baru-barrho-2}, we precisely get the claimed equation~\eqref{eq:barrho/3rd}.

\subsection{Equation for fluid velocity}
As the obtained equation~\eqref{eq:barrho/3rd} for $\bar\rho_\e$ shows that homogeneous spatial densities are stable to order $O(\e^3)$, we can focus for simplicity on the homogeneous setting,
\[\bar\rho_\e\,\equiv\,\tfrac1{\omega_d}+O(\e^3).\]
In addition, we shall focus on the case of infinite P\'eclet number and of vanishing particle swimming velocity,
\[\Pe=\infty,\qquad U_0=0,\]
which substantially simplifies the macroscopic fluid equations (recall however that our rigorous results do not hold for $\Pe=\infty$).
We start by computing the contribution of the elastic stress $\sigma_1$ in~\eqref{eq:baru-barrho-2}.
Using that
\[n\otimes n-\tfrac1d\Id=-\tfrac1{2d}\Delta_n(n\otimes n-\tfrac1d\Id),\]
we have by definition
\begin{eqnarray*}
\sigma_1[g_1+\e g_2+\e^2g_3]&=&\lambda\theta\int_{\Sp^{d-1}}(n\otimes n-\tfrac1d\Id)\,(g_1+\e g_2+\e^2g_3)(\cdot,n)\,\dd n\\
&=&-\lambda\theta\tfrac{1}{2d}\int_{\Sp^{d-1}}(n\otimes n-\tfrac1d\Id)\,\Delta_n(g_1+\e g_2+\e^2g_3)(\cdot,n)\,\dd n.
\end{eqnarray*}
Using the defining equations for $g_1,g_2,g_3$ with $U_0=0$, in form of
\begin{multline*}
\Delta_n(g_1+\e g_2+\e^2g_3)=
\tfrac{1}{\omega_d}\Div_n(\pi_n^\bot(\nabla\bar u_\e)n)
+\e(\partial_t+\bar u_\e\cdot\nabla_x)(g_1+\e g_2)\\
+\e\Div_n(\pi_n^\bot(\nabla\bar u_\e)n(g_1+\e g_2))+O(\e^3),
\end{multline*}
the above becomes, after straightforward simplifications and integrations by parts,
\begin{align*}
&\sigma_1[g_1+\e g_2+\e^2g_3]_{ij} \\
& =
\tfrac{1}{2 \omega_d}\lambda \theta\D(\bar u_\e):\int_{\Sp^{d-1}}(n\otimes n)(n\otimes n-\tfrac1d\Id)_{ij}\,\dd n\\
& \quad -\e\lambda\theta\tfrac{1}{2d}(\partial_t+\bar u_\e\cdot\nabla)\int_{\Sp^{d-1}}(n\otimes n-\tfrac1d\Id)_{ij}(g_1+\e g_2)(\cdot,n)\,\dd n\\
& \quad +\e\lambda\theta\tfrac{1}{2d}(\nabla\bar u_\e)_{kl}\int_{\Sp^{d-1}}(\delta_{ik} n_jn_l+\delta_{jk}n_in_l-2n_in_jn_kn_l)\,(g_1+\e g_2)(\cdot,n)\,\dd n \\
& \quad +O(\e^3).
\end{align*}
Now inserting the explicit expressions for $g_1,g_2$ obtained in Proposition~\ref{prop:first-terms}, in form of
\begin{multline}\label{eq:repeat-g1g2}
(g_1+\e g_2)(\cdot,n)\,=\,
\tfrac1{2 \omega_d}(n\otimes n-\tfrac1d\Id):\Big(\D(\bar u_\e)-\e\tfrac1{4d} A_2(\bar u_\e)+\e\tfrac1d\D(\bar u_\e)^2\Big)\\
+\e\tfrac{1}{8\omega_d}\Big(\big((n\otimes n):\D(\bar u_\e)\big)^2-\tfrac2{d(d+2)}\Tr(\D(\bar u_\e)^2)\Big)
+O(\e^2),
\end{multline}
using the explicit integrals~\eqref{eq:explicit-integral-n}, and using that $\Tr(A_2(u))=4\Tr(\D(u)^2)$,
we are led to
\begin{align*}
&\sigma_1[g_1+\e g_2+\e^2g_3]
\\
& \qquad =\lambda\theta\tfrac{1}{2d(d+2)}A_1(\bar u_\e)
-\e\lambda\theta\tfrac{1}{d^2(d+2)(d+4)}A_1(\bar u_\e)^2\\
& \qquad \quad -\e\lambda\theta\tfrac{1}{4d^2(d+2)}\Big(
(\partial_t+\bar u_\e\cdot\nabla)A_1(\bar u_\e)
-(\nabla\bar u_\e)A_1(\bar u_\e)
-A_1(\bar u_\e)(\nabla\bar u_\e)^T\Big)\\
& \qquad \quad +\e^2\lambda\theta\tfrac{1}{8d^3(d+2)}\Big(
(\partial_t+\bar u_\e\cdot\nabla) A_2(\bar u_\e)
-(\nabla\bar u_\e)A_2(\bar u_\e)
-A_2(\bar u_\e)(\nabla\bar u_\e)^T\Big)\\
& \qquad \quad -\e^2\lambda\theta\tfrac{1}{4d^3(d+4)}\Big(
(\partial_t+\bar u_\e\cdot\nabla)A_1(\bar u_\e)^2
-(\nabla\bar u_\e)A_1(\bar u_\e)^2
-A_1(\bar u_\e)^2(\nabla\bar u_\e)^T\Big)\\
& \qquad \quad +\e^2\lambda\theta\tfrac{1}{4d^3(d+2)(d+4)}\Big(A_1(\bar u_\e)A_2(\bar u_\e)+A_2(\bar u_\e)A_1(\bar u_\e)\Big)\\
& \qquad \quad -\e^2\lambda\theta\tfrac{(5d+12)}{4d^3(d+2)(d+4)(d+6)}A_1(\bar u_\e)^3\\
& \qquad \quad -\e^2\lambda\theta\tfrac{(d+3)}{4d^2(d+2)^2(d+4)(d+6)}A_1(\bar u_\e)\Tr(A_1(\bar u_\e)^2)
+R_\e\Id+O(\e^3),
\end{align*}
for some scalar field $R_\e$ that will be absorbed in the pressure field.
Now recalling the definition of Rivlin--Ericksen tensors, cf.~\eqref{eq:RE-tensor}, and appealing to the following matrix identities, for any matrix $B$,
\begin{gather*}
-(\nabla u)B-B(\nabla u)^T\,=\,
(\nabla u)^TB+B(\nabla u)
-\big(A_1(u)B+BA_1(u)\big),\\
(\partial_t+u\cdot\nabla)A_1(u)^2-(\nabla u)A_1(u)^2-A_1(u)^2(\nabla u)^T
\,=\,A_1(u)A_2(u)+A_2(u)A_1(u)-3A_1(u)^3,
\end{gather*}
we are led to
\begin{align*}
&\sigma_1[g_1+\e g_2+\e^2g_3] \\
&=
\lambda\theta\tfrac{1}{2d(d+2)}A_1(\bar u_\e)
-\e\lambda\theta\tfrac{1}{4d^2(d+2)}A_2(\bar u_\e)
+\e\lambda\theta\tfrac{1}{2d^2(d+4)}A_1(\bar u_\e)^2\\
& \quad +\e^2\lambda\theta\tfrac{1}{8d^3(d+2)}A_3(\bar u_\e)
-\e^2\lambda\theta\tfrac{3}{8d^3(d+4)}
\Big(A_1(\bar u_\e)A_2(\bar u_\e)
+A_2(\bar u_\e)A_1(\bar u_\e)\Big)\\
& \quad +\e^2\lambda\theta\tfrac{(3d^2+19d+24)}{4d^3(d+2)(d+4)(d+6)}A_1(\bar u_\e)^3
-\e^2\lambda\theta\tfrac{(d+3)}{4d^2(d+2)^2(d+4)(d+6)}A_1(\bar u_\e)\Tr(A_1(\bar u_\e)^2)\\
& \quad +R_\e\Id+O(\e^3).
\end{align*}
Further recalling that $A_1(u)^3=\tfrac12A_1(u)\Tr(A_1(u)^2)+\tfrac13\Tr(A_1(u)^3)$ in dimension $d\le 3$ by the Cayley--Hamilton theorem with $\Tr(A_1(u))=0$,
we obtain
\begin{align}\label{eq:sigma1-expand3}
\begin{split}
&\sigma_1[g_1+\e g_2+\e^2g_3] \\
&=
\lambda\theta\tfrac{1}{2d(d+2)}A_1(\bar u_\e)
-\e\lambda\theta\tfrac{1}{4d^2(d+2)}A_2(\bar u_\e)
+\e\lambda\theta\tfrac{1}{2d^2(d+4)}A_1(\bar u_\e)^2\\
& \quad +\e^2\lambda\theta\tfrac{1}{8d^3(d+2)}A_3(\bar u_\e)
-\e^2\lambda\theta\tfrac{3}{8d^3(d+4)}
\Big(A_1(\bar u_\e)A_2(\bar u_\e)
+A_2(\bar u_\e)A_1(\bar u_\e)\Big)\\
& \quad +\e^2\lambda\theta\tfrac{3d^2+11d+12}{8d^3(d+2)^2(d+6)}A_1(\bar u_\e)\Tr(A_1(\bar u_\e)^2)
+R_\e\Id+O(\e^3),
\end{split}
\end{align}
up to modifying the scalar field $R_\e$.

Next, we turn to the computation of the contribution of the viscous stress $\sigma_2$ in~\eqref{eq:baru-barrho-2}. We have by definition
\begin{equation*}
\sigma_2[\tfrac{1}{\omega_d}+\e g_1+\e^2g_2,\nabla\bar u_\e]
\,=\,\lambda\int_{\Sp^{d-1}}(n\otimes n)(\nabla\bar u_\e)(n\otimes n)\,(\tfrac{1}{\omega_d}+\e g_1+\e^2g_2)\, \dd n,
\end{equation*}
and thus, using the explicit integrals~\eqref{eq:explicit-integral-n},
\begin{multline*}
\sigma_2[\tfrac{1}{\omega_d}+\e g_1+\e^2g_2,\nabla\bar u_\e]\\
\,=\,\lambda\tfrac{2}{d(d+2)} \D(\bar u_\e)
+\e\lambda\int_{\Sp^{d-1}}(n\otimes n)(\nabla\bar u_\e)(n\otimes n)\,(g_1+\e g_2)\, \dd n+R_\e\Id,
\end{multline*}
where again $R_\e$ stands for some scalar field that can be absorbed in the pressure field.
Inserting the explicit expressions for $g_1,g_2$ obtained in Proposition~\ref{prop:first-terms}, in form of~\eqref{eq:repeat-g1g2}, and computing integrals as above, we easily find
\begin{align*}
\sigma_2[\tfrac{1}{\omega_d}+\e g_1+\e^2g_2,\nabla\bar u_\e]
& = \lambda\tfrac{1}{d(d+2)} A_1(\bar u_\e)
+\e\lambda\tfrac{1}{d(d+2)(d+4)}A_1(\bar u_\e)^2\\
&\quad -\e^2\lambda\tfrac{1}{4d^2(d+2)(d+4)}\big(A_1(\bar u_\e)A_2(\bar u_\e)+A_2(\bar u_\e)A_1(\bar u_\e)\big)\\
& \quad +\e^2\lambda\tfrac{(5d+12)}{4d^2(d+2)(d+4)(d+6)}A_1(\bar u_\e)^3\\
& \quad +\e^2\lambda\tfrac{(d^2-2d-12)}{8d^2(d+2)^2(d+4)(d+6)}A_1(\bar u_\e)\Tr(A_1(\bar u_\e)^2)
+O(\e^3)+R_\e\Id.
\end{align*}
Using again $A_1(u)^3=\tfrac12A_1(u)\Tr(A_1(u)^2)+\tfrac13\Tr(A_1(u)^3)$ in dimension $d\le 3$, this becomes
\begin{align*}
\sigma_2[\tfrac{1}{\omega_d}+\e g_1+\e^2g_2,\nabla\bar u_\e]
&=\lambda\tfrac{1}{d(d+2)} A_1(\bar u_\e)
+\e\lambda\tfrac{1}{d(d+2)(d+4)}A_1(\bar u_\e)^2\\
& \quad -\e^2\lambda\tfrac{1}{4d^2(d+2)(d+4)}\big(A_1(\bar u_\e)A_2(\bar u_\e)+A_2(\bar u_\e)A_1(\bar u_\e)\big)\\
& \quad +\e^2\lambda\tfrac{(3d^2+10d+6)}{4d^2(d+2)^2(d+4)(d+6)}A_1(\bar u_\e)\Tr(A_1(\bar u_\e)^2)
+R_\e\Id
+O(\e^3).
\end{align*}
Inserting this together with~\eqref{eq:sigma1-expand3} back into~\eqref{eq:baru-barrho-2}, we precisely get the claimed third-order fluid equation~\eqref{eq:baru/3rd}.
\qed

\section*{Acknowledgements}
We thank David G\'erard-Varet for related discussions, we thank Dennis Trautwein for pointing out references on viscoelastic models with stress diffusion, and we thank Leonid Berlyand and Spencer Dang for references on the shear-thinning/thickening behavior of active suspensions. We also thank the referees for useful comments that have led to a strong improvement of this work.

M.D. acknowledges financial support from the F.R.S.-FNRS, as well as from the European Union (ERC, PASTIS, Grant Agreement n$^\circ$101075879).\footnote{{Views and opinions expressed are however those of the authors only and do not necessarily reflect those of the European Union or the European Research Council Executive Agency. Neither the European Union nor the granting authority can be held responsible for them.}}
L.E. acknowledges financial support from the ESPRC STABLE-CHAOS.
R.H. has been supported  by the German National Academy of Science Leopoldina, grant LPDS 2020-10.

\subsection*{Data Availability Statement}
Data sharing is not applicable to this article as no datasets were
generated or analyzed during the study.

\subsection*{Conflicts of interests/Competing interests.} The authors have no relevant financial or non-financial interests to disclose. The authors have no competing interests to declare that are relevant to the content of this article.

\end{document}